\documentclass[preprint,review, 3p, 11pt]{elsarticle}

\usepackage{amssymb}
\usepackage{amsmath}
\usepackage{amsthm}

\usepackage{algorithm}
\usepackage{algpseudocode}
\usepackage{tikz}
\usepackage{natbib}
\usepackage{amsfonts, amsmath, amssymb,dsfont}
\usepackage{algorithm}
\usepackage{comment}
\usepackage{tikz-cd}

\newcommand{\mcl}[1]{\mathcal{#1}}
\usepackage{thmtools}
\usepackage{thm-restate}

\usepackage{hyperref}

\usepackage{cleveref}

 \numberwithin{dummy}{section}
\newtheorem{theorem}{Theorem}
\newtheorem{proposition}[theorem]{Proposition}
\newtheorem{corollary}[theorem]{Corollary}
\newtheorem{lemma}[theorem]{Lemma}
\newtheorem{definition}[theorem]{Definition}
\newtheorem{assumption}{Assumption}
\newtheorem{conjecture}{Conjecture}
\newtheorem{remark}{Remark}

\makeatother
\newcommand{\redtext}[1]{{\leavevmode\color{red}#1}}
\newcommand{\bluetext}[1]{{\leavevmode\color{blue}#1}}
\newcommand{\greentext}[1]{{\leavevmode\color{green!50!black!100}#1}}
\newcommand{\sspace}{\mcl{S}^{(N)}}
\newcommand{\sspaceN}[1]{\mcl{S}{(#1)}}

\newcommand{\pspaceN}[1]{\mcl{P}^{(#1)}}
\newcommand{\mcY}[1]{\mcl{Y}(#1)}
\newcommand{\mcU}[1]{\mcl{U}(#1)}

\newcommand{\Pmat}[2]{\mathbf{P}^{(#1)}_{#2}} 
\newcommand{\PmatN}{\mathbf{P}}
\newcommand{\PmatNN}[1]{\mathbf{P}_{#1}}
\newcommand{\budget}{\alpha}
\newcommand{\budgetK}{\mathbf{b}^{E}}
\newcommand{\budgetk}[1]{\mathbf{b}^{(#1)}}
\newcommand{\constrWMDP}[3]{d^{#1}(#2, #3)}
\newcommand{\Rvec}[1]{\mathbf{r}^{(#1)}}
\newcommand{\RCvec}[1]{\Tilde{\mathbf{r}}^{(#1)}}
\newcommand{\RvecN}{\mathbf{r}}
\newcommand{\rew}[3]{r^{(#1)}(#2, #3)}
\newcommand{\syncconst}[1]{\rho_{#1}}

\newcommand{\bsN}{\mathbf{s}^{(N)}}
\newcommand{\bsn}[1]{s_{#1}}
\newcommand{\bsNN}{\mathbf{\hat{s}}^{(N)}}
\newcommand{\bsnN}[1]{\hat{s}_{#1}}
\newcommand{\bs}{\mathbf{s}}
\newcommand{\bSN}{\mathbf{S}^{(N)}}
\newcommand{\bSNN}{\mathbf{S}}

\newcommand{\baN}{\mathbf{a}^{(N)}}
\newcommand{\ban}[1]{a_{#1}}
\newcommand{\ba}{\mathbf{a}}
\newcommand{\bAN}{\mathbf{A}^{(N)}}
\newcommand{\bANN}{\mathbf{A}}

\newcommand{\VpolavgN}[3]{\mathbf{V}_{#1}^{#2}(#3)}
\newcommand{\VpolavgNnos}[2]{\mathbf{V}_{#1}^{#2}}
\newcommand{\VoptN}[1]{\mathbf{V}_{\text{OPT}}^{#1}}
\newcommand{\VoptNoN}{\mathbf{V}_{\text{OPT}}}

\newcommand{\VLPN}{\mathbf{V}_{\text{LP}}}
\newcommand{\VLPF}{\mathbf{V}_{\text{finite}}}
\newcommand{\VoptWMDP}[1]{\mathbf{V}_{\text{OPT-WMDP}}^{#1}}

\newcommand{\THor}{\tau}
\newcommand{\VLPROB}[2]{\mathbf{W}^{\text{FL}}\left(#1, #2\right)}
\newcommand{\VLPROBn}[4]{\Tilde{\mathbf{W}}_{#3}^{#4}\left(#1, #2\right)}
\newcommand{\Qval}[3]{\mathbf{Q}^{\text{FL}}_{\mu}(#1, #2, #3)}
\newcommand{\LrelT}[3]{\mathbf{L}^{#3}(#2, #1)}

\newcommand{\bx}{\mathbf{x}}
\newcommand{\bu}{\mathbf{u}}
\newcommand{\bualg}[1]{\bu_{\text{A2}(#1)}}

\newcommand{\bxs}[1]{\mathbf{x}(#1)}
\newcommand{\bxn}[2]{\mathbf{x}^{#2}(#1)}
\newcommand{\bhx}{\hat{\mathbf{x}}}
\newcommand{\bAalg}[1]{\mathbf{A}_{\text{A2}}(#1)}
\newcommand{\state}{\omega}

\newcommand{\bYn}[4]{y^{(#1)}_{#2, #3}(#4)} 
\newcommand{\bYnf}[3]{y^{(#1)}_{#2, #3}} 
\newcommand{\bYVn}[2]{\mathbf{y}^{#1}(#2)}
\newcommand{\bYVZn}[2]{\mathbf{y}^{(#1)}_{\cdot, 0}(#2)}
\newcommand{\bYVUn}[2]{\mathbf{y}^{(#1)}_{\cdot, 1}(#2)}
\newcommand{\bYNt}[1]{\mathbf{y}(#1)}
\newcommand{\bYVUN}[1]{\mathbf{y}_{\cdot, 1}(#1)}
\newcommand{\bYVZN}[1]{\mathbf{y}_{\cdot, 0}(#1)}
\newcommand{\bYVN}[2]{\mathbf{y}({#1})_{#2}}

\newcommand{\bYFLTn}[2]{\mathbf{y}^{\text{FL}, (#2)}(#1)}
\newcommand{\bYFLT}[1]{\mathbf{y}^{\text{FL}}(#1)}
\newcommand{\bYFLU}[1]{\mathbf{y}^{\text{FL}}_{\cdot, 1}(#1)}

\newcommand{\bYFLZn}[2]{\mathbf{y}_{#2}^{\text{FL}}_{\cdot, 0}(#1)}
\newcommand{\bYFLUn}[3]{\mathbf{y}_{#2, 1}^{\text{FL},#1}(#3)} %

\newcommand{\bYFLAn}[3]{\mathbf{y}_{\cdot, #2}^{\text{FL},#1}(#2)}
\newcommand{\ufl}[1]{\mathbf{u}^{\text{FL}}(#1)}

\newcommand{\bYWMDPn}[2]{\mathbf{y}_{\cdot,a}^{#2}(#1)}
\newcommand{\bYWMDPN}[1]{\mathbf{y}_{\cdot,a}^{#1}}
\newcommand{\VLWMDP}[2]{\mathbf{W}^{\text{FL}}_{\text{MDP}}\left(#1, #2\right)}
\newcommand{\VLPinfWMDP}[1]{\mathbf{V}_{\text{LP-} #1}^{\text{WMDP}}}

\newcommand{\gstar}{\mathbf{g}^{\star}}
\newcommand{\hstar}[1]{h^{\star}(#1)}

\newcommand{\bYNn}[1]{\mathbf{y}^{(#1)}}
\newcommand{\bYN}{\mathbf{y}}
\newcommand{\bYZN}{\mathbf{y}_{\cdot, 0}}
\newcommand{\bYUN}{\mathbf{y}_{\cdot, 1}}

\newcommand{\bYstarN}{\mathbf{y}^{\star}}
\newcommand{\bYstar}[1]{\mathbf{y}^{\star, (#1)}}
\newcommand{\bYZstarN}{\mathbf{y}^{\star}_{\cdot, 0}}
\newcommand{\bYUstarN}{\mathbf{y}^{\star}_{\cdot, 1}}

\newcommand{\bYinf}[1]{\mathbf{y}^{#1}}
\newcommand{\bYZinf}[1]{\mathbf{y}^{#1}_{\cdot, 0}}
\newcommand{\bYUinf}[1]{\mathbf{y}^{#1}_{\cdot, 1}}
\newcommand{\blaminf}{\mathbf{\lambda}_{\infty}}

\newcommand{\laginf}{\mu}

\newcommand{\blamT}{\mathbf{\lambda}}
\newcommand{\blam}[1]{\lambda(#1)}
\newcommand{\blamstarT}{\mathbf{\lambda}^{\star}}
\newcommand{\blamstar}[1]{\lambda^{\star}(#1)}

\newcommand{\indx}[2]{\mathbf{i}_{#1}(#2)}
\newcommand{\indxlp}[1]{\mathbf{i}_{#1}^{\text{lp}}}
\newcommand{\Windx}[2]{\mathbf{i}^{\text{WI}}_{#1}(#2)}

\newcommand{\E}{\mathbb{E}}
\newcommand{\ProbP}{\mathbb{P}}
\newcommand{\Prob}[1]{\ProbP\left(#1\right)}
\newcommand{\Real}[1]{\mathbb{R}^{#1}}
\newcommand{\Onne}{\mathbf{1}}

\newcommand{\floor}[1]{\lfloor #1\rfloor}

\newcommand{\costt}[2]{c(#1, #2)}
\newcommand{\rotcost}[2]{\Tilde{#1}(#2)}
\newcommand{\CostT}[2]{C_{#2}(#1)}

\newcommand{\stor}[1]{\psi(#1)}

\newcommand{\bxc}[1]{\{\bx^{(1)}, \bx^{(2)} \dots \bx^{(#1)} \dots \bx^{(N)}\}}
\newcommand{\bxcn}[1]{\{\bx^{(1)}, \bx^{(2)} \dots \Tilde{\bx}^{(#1)} \dots \bx^{(N)}\}}

\newcommand{\gstarWMDP}{\gstar_{\text{WMDP}}}
\newcommand{\hstarWMDP}[1]{h^{\star}_{\text{WMSP}}(#1)}

\journal{European Journal of Operations Research}

\newcommand\GN[1]{{\color{green!50!black!90}\textbf{Nicolas:}#1}}

\begin{document}
\begin{frontmatter}
\title{Model Predictive Control is almost Optimal for Heterogeneous Restless Multi-armed Bandits}
\author[1]{Dheeraj Narasimha\corref{cor1}%
\fnref{fn1}}
\ead{dheeraj.narasimha@gmail.com}

\affiliation[1]{organization={LAAS CNRS},
            addressline={ 7 Av. du Colonel Roche}, 
            city={Toulouse},
            postcode={31500}, 
            country={France}}
\author[2]{Nicolas Gast \fnref{fn2}}

\affiliation[2]{
organization = {INP Grenoble, INRIA, CNRS},
addressline = {Batiment IMAG
150 place du Torrent},
city = {St Martin d'Heres},
postcode={38400},
country={France}
}

\begin{abstract}
We consider a general infinite horizon Heterogeneous Restless multi-armed Bandit (RMAB). Heterogeneity is a fundamental problem for many real-world systems largely because it resists many concentration arguments. 
In this paper, we assume that each of the $N$ arms can have different model parameters. 
Model predictive control is a well-known control strategy that repeatedly solves a finite-horizon optimization problem of length $\THor$ to produce a policy that can be applied to an infinite-horizon setting. In this paper, we adopt this approach by repeatedly solving a finite linear program, yielding what we call the LP-update policy for the infinite-horizon problem. Under a mild assumption of uniform ergodicity, we show an $\mcl{O}\left(\sqrt{1/N}\right)$ suboptimality gap on this well-known algorithm that works very well in practice. In addition to the LP-update policy we are able to derive a finite-horizon policy (LP-update with recomputation) that segments the infinite time horizon into finite horizon problems that allow us to explicitly connect the length of computation time to the acceptable error tolerance.
Our simulations demonstrate that our algorithm works extremely well even when this finite-horizon, $\tau$, is very small (in our case $5$), which makes it computationally efficient.  Our theoretical results draw on techniques from the model predictive control literature by invoking the concept of \emph{dissipativity} and generalize quite easily to the more general weakly coupled heterogeneous Markov Decision Process setting. In addition, we draw a parallel between our own policy and the LP-index policy by showing that the LP-index policy corresponds to $\tau=1$.  
\end{abstract}



\end{frontmatter}

\section{Introduction and Related Works}

We consider the sequential decision-making problem known as the \textbf{Restless Multi-Armed Bandit (RMAB)} over an infinite discrete-time horizon. Our problem consists of \( N \) weakly coupled arms, each arm is treated as being statistically different from the others. At any given time step, the decision-maker must select a subset of no more than $N \alpha$ arms to pull which we call the \emph{budget constraint}. Each arm inhabits its own state space and may evolve differently depending on the selected action. Once the subset is chosen, each arm yields a reward based on its state and action, and transitions independently to a new state according to a state-action dependent transition kernel. Hence, the arms are ``weakly coupled'', interacting only through shared constraints and our selection policy. Both the transition kernels and reward functions of all arms are known beforehand. The objective is to design a policy that maps joint arm states to actions to maximize the long-run average reward.

Sequential decision-making problems with weakly coupled constraints arise in many domains, including communication systems \cite{Borkar2018}, queueing theory \cite{nino2002dynamic}, web crawling \cite{MearaWebcrawl01}, and scheduling \cite{dance2019optimal, veatch1996scheduling}, among others. A comprehensive overview of RMABs and their applications is presented in \cite{NinoMora23}. The RMAB problem is PSPACE-hard, as shown in \cite{PT99}, motivating the study of algorithms that are optimal as \( N \to \infty \), that is, asymptotically optimal.
\subsection{A history of Index Policies.}
The RMAB problem was first introduced by Whittle in his seminal work \cite{Wh88}, where he conjectured that under a condition called \emph{indexability}, a priority-based policy—\textbf{Whittle's index policy}—would be asymptotically optimal. This conjecture was later verified in the homogeneous setting (i.e., statistically identical arms) under a \emph{uniform global attractor condition}(UGAP) in \cite{WW90}, who also demonstrated counterexamples when this condition fails. Subsequently, many lines of research in the area have sought to relax both the attractor and indexability conditions, particularly in the homogeneous setting.
\subsection{Index Policies, Asymptotic Optimality and the UGAP assumption.}
A generalization of Whittle's index known as the \textbf{LP-priority index policy} was introduced in \cite{verloop2016asymptotically}, removing the need for indexability. While this work does not strictly require homogeneity it cannot cover the case where the heterogeneity scales linearly with the number of arms. More recent work \cite{GGY23, GGY23b} shows that such priority policies can be \textbf{exponentially close to optimal} under three key conditions centered around the optimal fixed point: unichain and aperiodicity of the underlying Markov process, the UGAP condition for the dynamical system and non-degeneracy of the Linear Program for the fixed point solution.
\subsection{Breaking the UGAP assumption.}
 On the other hand, the UGAP assumption was first eliminated in \cite{HXCW23}, where the authors constructed an unconstrained virtual system driven to the fixed point. To obtain an asymptotically optimal policy, they introduced a \emph{synchronization assumption}, ensuring the real system aligns with the virtual one. This method of \textbf{steering systems to optimal fixed points} has inspired subsequent works such as \cite{HXCW24, yan2024}. For instance, \cite{HXCW24} proposed a \emph{two-set policy} to actively synchronize non-aligned arms, recovering exponential optimality by leveraging the fixed-point structure as in \cite{GGY23b}. In contrast, \cite{yan2024} designed an \emph{align-and-steer} policy that guides the mean field control system to the fixed point under a \emph{reachability} assumption.

\subsection{Homogeneous vs Heterogeneous Setting}
In contrast to these rich developments in the homogeneous case, results for \textbf{heterogeneous RMABs} are relatively scarce. For the infinite-horizon case, \cite{Avrachenkov2024} constructed a fluid-policy framework which resembles \cite{yan2024} align and steer policy, however, its application for a heterogeneous problem where each arm lives in its own state space is unclear since such a system does not have a notion of empirical measure that can be aligned with the optimal fixed point policy. Close in spirit to our results is the recent work of ~\cite{zhang2025projLyapunov} which proves an order \(\frac{1}{\sqrt{N}}\) result using a Lyapunov projection method. The core idea of their work is to create an expanding set of arms whose actions align with the optimal fixed point actions (used previously in the homogeneous setting in \cite{HXCW24}), aligning the system with the invariant measure. In terms of theoretical results, our work primarily differs from their work by the assumptions made for ensuring performance. While our assumptions require a form of ergodicity in the \(-0-\) action, their assumptions require a mixing time result for the \emph{optimal single arm policy}. To the best of our knowledge these assumptions cannot directly be compared to each other. Another point of difference is, algorithmically speaking we wish to prove that a well known control policy Model Predictive Control (MPC) will perform well in the heterogeneous setting whereas their work proves results for a new Lyapunov based priority policy on the arms.  
\subsection{Model Predictive Control and LP-Update Policy.}
 The model predictive control method is a widely used and widely studied advanced control technique. The controller uses a model to predict the system behavior for a finite time in the future, denoted by $\THor$ for the rest of this paper, to take a single action. The system is then allowed to evolve by a single step and the process is reiterated. It has found applications in a range of different fields including control of autonomous ground vehicles ~\cite{Yu2021ModelPC}, renewable energy and management of electric microgrids ~\cite{2022SulmanMPCmicrogrid} and  water resource management ~\cite{CASTELLETTI2023442}. For a more detailed look at model predictive control, we direct the readers attention to the survey ~\cite{schwenzer_review_2021} and the references therein. In the homogeneous setting, the model predictive control strategy is equivalent to solving a finite-horizon linear program iteratively, called the LP-update policy. \cite{gast2024MPCoptimal} showed that the LP-update policy achieves asymptotic optimality in the  infinite horizon setting. Our work builds on those ideas to generalize their results to the infinite horizon heterogeneous RMAB problem.
 
Several authors have proposed variations of the LP-update policy, we highlight here the work by~\cite{brown2020index} for the finite-horizon heterogeneous RMAB, who posed the infinite-horizon extension as an open problem. It should be noted here that their optimality gap grows exponentially over time, but since the horizon is limited this result does not feature in the order wise bound. Other works that have built upon variations of this LP-update algorithm, include ~\cite{hu2017asymptotically, zayas2019asymptotically, ghosh2022indexability, zhang2021restless, GGY23} and references therein. 
\subsection{Points of Comparison.}
Closest to our work is \cite{Brown&Zhang25} which looks to use LP-update policies for a finite horizon problem and the infinite horizon discounted reward problem. Our work differs from theirs in several respects. First, in terms of the model predictive control idea, their finite-horizon results do not directly translate to the infinite-horizon case because infinite-horizon problems lack terminal conditions. In particular, although a finite-horizon policy with horizon $T$ may recompute a policy after one step using horizon $T - 1$, this behavior is not comparable to the infinite-horizon setting, where the computational horizon remains finite over time. Second, their sufficiency conditions for order $T\sqrt{N}$ bounds (specifically the conditions for Theorem 4.1) imply our sufficient conditions on Ergodicity of the $-0-$ action making our conditions strictly less restrictive. Finally, it should be noted that their infinite horizon discounted result [Theorem 4.5] relies on truncation which cannot be applied to the average reward problem. We elaborate a little more on this final point in Section \ref{SEC::PROOF}. However, as we demonstrate in Section \ref{SUBSEC::MPC2}, one can still leverage these ideas to create an infinite horizon solution to the average reward problem by repeatedly segmenting the time horizon into pieces of equal length. The non-trivial question under this variation of the algorithm is precisely how large do those pieces need to be? We answer this question in Theorem \ref{THM::ALMOSTOPT2}.

\subsection{Main Contributions}

Our main contributions can be summarized as follows:
\begin{enumerate}
    \item We establish an \( O\left(\frac{1}{\sqrt{N}} \right) \) bound on the difference between the optimal value of the average rewards in the Restless Multi-Armed Bandit (RMAB) setting and the average reward accumulated by the LP-update policy under an easily verifiable ergodicity condition. In doing so, we resolve a key open problem posed by \cite{brown2020index}. In addition, our LP-update strategy has an excellent empirical performance and is easy to implement.
    \item Our proof technique is modular in nature, meaning that our ergodicity conditions or concentration techniques can be independently replaced or extended. Our proof can be broken into two main components, 1) a convergence result for the approximate dynamical system and 2) a concentration result which we call the \emph{Jensen gap result}. Each of these problems can be solved independently, making the techniques we use easy to generalize to a wide range of constrained Markov Decision processes. Our dynamical system uses a novel dissipative framework, previously only used by \cite{gast2024MPCoptimal} in the homogeneous setting. This is the key concept that allows us to show that a finite horizon linear program solution can be used to solve an infinite horizon problem. An important example, the weakly coupled heterogeneous Markov decision process setting, where we recover the same results without much modification to our algorithm.
    \item 
    In Section \ref{SUBSEC::MPC2}, we show a variation of the finite horizon LP-update policy that segments the infinite horizon into pieces of finite length in terms of a tolerance error. This resolves any concern one may have about the asymptotic nature of our time horizon $\tau$. Our analysis is tied closely with the mechanisms we established in the conventional LP-update strategy.
    \item We show that our policy reduces to the well-known LP-index policy when we set the horizon $\tau=1$. Increasing the horizon to a higher value allows us to avoid the necessity of the notoriously hard to verify UGAP assumption.   
    \item Numerically our algorithm outperforms the state-of-the-art algorithm by \cite{zhang2025projLyapunov} in all cases even when the time horizon is set to a fixed value of $\THor = 5$. These cases include bandit instances developed specifically to ensure LP-priority policies perform poorly. Indicating that slightly increasing $\THor$ allows the LP-update to avoid the pit-falls surrounding LP-index policy when the UGAP assumption fails.
\end{enumerate}

\textbf{Road Map for the Rest of the Paper:\\}
In Section \ref{SEC::SYS}, we present a concrete system model for the Restless Multi-Armed Bandits (RMABs). Next, we present the LP-update algorithm in Section \ref{SEC::IPL} with different time horizons $\THor$. In Section \ref{SEC::MAINRESULTS} we introduce our main result, Theorem \ref{THM::ALMOSTOPT} for the LP-update algorithm with a fixed value of $\THor$, we then refine the result under additional conditions and finally present its generalization to Weakly coupled MDPs under equivalent ergodicity assumptions. Further, in the same section we present Theorem \ref{THM::ALMOSTOPT2} when our algorithm divides the time into finite time segments using ideas from Theorem \ref{THM::ALMOSTOPT}. Section \ref{SEC::PROOF} introduces the main ideas necessary for the proof including the concept of \emph{dissipativity} and \emph{Jensen's gap}. We end the paper with Section \ref{SEC:NUM} begins by comparing our policy explicitly with the LP-priority policy followed by numerical simulations where we compare our algorithm to state-of-the-art results. The appendix contains the details of our proofs, additional details for our algorithms and numerical experiments.

\section{System Model}\label{SEC::SYS}

We consider an infinite-horizon discrete-time heterogeneous Restless multi-armed bandit problem parameterized by the tuple $\langle \sspace, \{\Pmat{n}{0}\}_{n = 1}^{n = N}, \{\Pmat{n}{1}\}_{n = 1}^{n = N}, \Rvec{n}; \budget \rangle$. A decision maker controls $N$ heterogeneous arms and must choose a subset of at most $\budget N$ arms to pull at each time, where $0 < \budget \leq 1$ is a known constant. If arm $n$ in state $s$ is pulled (the action is denoted by $1$), it follows a transition kernel $\Pmat{n}{1}$ obtains a reward $\rew{n}{s}{1}$. On the other hand, if it is left alone (the corresponding action is denoted by $0$), the arm transitions according to $\Pmat{n}{0}$, and obtains a reward $\rew{n}{s}{0}$. The reward vector $\Rvec{n}$ is the concatenation of all such rewards. We will assume that each arm $n$ lives in a finite state space $\sspaceN{n}$ and that all rewards lie between $0$ and $1$ for each arm. So long as all rewards are finite, this assumption can be made without loss of generality by shifting and scaling rewards appropriately. We denote by $\sspace$ the joint state space i.e, the cartesian product of the individual spaces of each arm, $\Pi_{n = 1}^{N}\sspaceN{n}$. Given the state of the arms at time $t$, $[\bsn{1} \dots \bsn{N}] =:\bsN \in \sspace$ and the decision maker's action $[\ban{1} \dots \ban{N}] =: \baN \in \{0,1\}^N$, a reward $\frac{1}{N}\sum_{n} \rew{n}{\bsn{n}}{\ban{n}}$ is earned and the system of arms will transition jointly to the next state $\bsNN \in\sspace$ as follows:
\begin{align}\label{EQ:MKEVOL}
            &\ProbP(\bSN(t + 1) = \bsNN|\bSN(t) = \bsN,\bAN(t) = \baN, \dots \bSN(0), \bAN(0)) \nonumber\\
            &= \ProbP(\bSN(t + 1) = \bsNN|\bSN(t) = \bsN,\bAN(t) = \baN)
            = \Pi_{n = 1}^{N} P^{(n)}_{\bsnN{n}|\bsn{n}, \ban{n}}.
\end{align}
Where $P^{n}_{\bsnN{n}|\bsn{n}, \ban{n}}$ denotes the $(\bsn{n}, \bsnN{n})$ entry of the corresponding transition matrix $\Pmat{n} {\ban{n}}$ indicating the probability that an arm in state $\bsn{n}$ transitions to state $\bsnN{n}$. \emph{It is important to note that the arms are only coupled through the budget constraint $\budget$: We say that the arms are weakly coupled.}

Let $\pi$ be a stationary policy mapping each state to a probability of choosing joint actions across the $N$ arms, i.e, $\pi : \sspace \to \Delta (\{0, 1\}^{N})$ subject to the budget constraint, $\budget$ at each time i.e, if $\baN := \{a_1, a_2 \dots a_N\}$ is the realization of the joint actions of all the arms then, $\sum_{n}^{N}a_{n} \leq \budget N.$
Let $\pspaceN{N}$ denote the space of all such policies. Given a policy $\pi \in \pspaceN{N}$, and an initial set of arm states $\bSN(0) = \bs$, we define the average gain of policy $\pi$ as
\begin{equation}\label{EQ::VPOLT}
    \VpolavgN{\pi}{N}{\bs} = \lim_{T \to \infty}\frac{1}{T}\E_{\pi}\left[ \sum_{t = 0}^{T - 1} \frac{1}{N} \sum_{n = 1}^N \rew{n}{A_n(t)}{S_n(t)}\Bigg|\bSN(0) = \bs \right].
\end{equation}
Here $(S^N_n(t), A^N_n(t))$ denotes the state-action pair of the $n^\text{th}$ arm at time $t$ and $\rew{n}{A_n(t)}{S_n(t)}$ denotes the $S_n(t)^{\text{th}}$ entry of the $\Rvec{n}(\cdot, A_n(t))$ vector. From \cite{puterman2014markov}, one can show that for any stationary policy $\pi$, this limit is well defined. The heterogeneous infinite-horizon average reward problem or the heterogeneous restless multi-armed bandit problem dubbed RMAB is to find the optimal policy $\pi \in \pspaceN{N}$ that maximizes this limit. We call this optiomal value $\VoptN{N}$:
\begin{align}\label{EQ::VOPTT}
    \VoptN{N} &:= \max_{\pi \in \pspaceN{N}} \VpolavgN{\pi}{N}{\bs}.
\end{align}
The optimal policy exists and the limit is well defined. It is not difficult to show that under mild conditions, the value of $\VoptN{N}$ is independent of the initial state. Hence, we refrain from adding any explicit dependence on the initial state in our notation.

It is shown in \cite{PT99} showed that such problems are PSPACE hard. \emph{The RMAB problem then boils down to finding a computation policy whose optimality gap is small when the number of arms is large}, where the optimality gap of a policy $\pi$ is: $\VoptN{N} - \VpolavgNnos{\pi}{N}$. In this paper we will show that a simple Linear Program (LP) update based algorithm with updates that compute a \textbf{finite time horizon LP} is \emph{almost asymptotically optimal} for this problem. 

\paragraph*{Notation} The number of arms, $N$ is fixed; hence, for notational convenience, we will drop the explicit dependence on $N$ in all future sections. In particular, when referring to joint states, actions or value functions we will drop the superscript $N$, bold variations of parameters without a subscript or superscript indicate a concatenation of the corresponding parameters over all $N$ arms. 
We use subscript $n$ for scalar variables, for example, an arm $n$ is in state $s_n$ with state space $\sspaceN{n}$. The same is true for action $a_n$. We use super-script $n$ to indicate the $n^{\text{th}}$ vector in a tuple of $N$ vectors. As an example a joint probability distribution of the arms, $\bx$ is a collection of $N$ vectors $\{\bx^{(1)}, \dots \bx^{(N)}\}$ where $\bx^{(n)}_s$ is the probability of arm $n$ occupying state $s$. We distinguish between the dot product $a \cdot b$ and the use of the Dirac bracket notation $\langle \cdot, \cdot \rangle$ as follows, for any two vectors of the same dimension, $a = \{a_1 \dots a_N\}$ and $b = \{b_1 \dots b_N\}$, $a \cdot b := \sum_{i = 1}^{N} a_ib_i$ whereas $\langle a, b \rangle := \frac{1}{N}\sum_{i = 1}^{N} a_ib_i $. 

\section{The LP-update Policy}\label{SEC::IPL}

We will begin by stating a related fixed-point problem that stems from \emph{Whittle's relaxation}, \cite{Wh88,avrachenkov2022whittleIB}. Our decision variables correspond to state-action probabilities on the arms denoted by $\bYnf{n}{s}{a}$. Concretely, $\bYnf{n}{s}{a}$ is the probability of an arm $n$ in state $s$ playing an action $a$, it follows that the sum of these elements of any arm $n$ must sum to $1$. We denote by $\bYNn{n}$ the corresponding state-action vector. As stated in the notation section, we will denote the concatenation of all $N$ vectors by $\bYN$, it follows that $\bYN$ is an element of $\Pi_{n = 1}^{N}\Delta(\Real{\sspaceN{n} \times 2})$. We will call this set of elements $\mcl{Y}$. Now, consider the following fixed point problem :  
\begin{align}
    \gstar := \max_{\bYN \in \mcl{Y}} \frac{1}{N}\sum_{n = 1}^{N} & \Rvec{n}\cdot \bYinf{(n)} \label{EQ::WHITF1}\\
    \text{such that,}& \nonumber \\
    \bYZinf{(n)} + \bYUinf{(n)} =&  \bYZinf{(n)} \cdot \Pmat{n}{0} + \bYUinf{(n)} \cdot \Pmat{n}{1}  \hspace{0.2 in } \forall n \label{EQ::WHITF2}\\
    \sum_{n = 1}^{N} \| \bYUinf{n}\|_1 \leq& N. \alpha \hspace{0.2 in } \label{EQ::WHITF3}
\end{align}
Where, the $l^1$ norm in the last constraint is short for $\| \bYUinf{(n)}\|_1 := \sum_{s \in \sspaceN{n}} \bYnf{n}{s}{1}$. 
This problem is called a \emph{relaxation} of \eqref{EQ::VOPTT} since it can be shown that it is equivalent to a problem where the action constraints are satisfied on average over time i.e., $ \lim_{T \to \infty}\frac{1}{T} \sum_{t = 0}^{T - 1} \sum_{n = 1}^{N}\| \bYUinf{(n)}(t)\|_1 \leq N \alpha$ rather than at each time instance. Note that our reward is linear and the space of solutions is non-empty, hence, a unique maximizer must exist for this problem along with an equivalent dual problem. 

The following definitions will be very helpful in the sections to come. 
\begin{definition}\label{DEF:GSTAR}\label{DEF:WHITMUL}
    \begin{itemize}
        \item We denote by $\gstar$ the solution of \eqref{EQ::WHITF1}-\eqref{EQ::WHITF3} which we shall call \textbf{gain}
        \item $\bYN^{\star}$ is a corresponding fixed point solution to \eqref{EQ::WHITF1}-\eqref{EQ::WHITF3}.
        \item $\laginf := \{\laginf^{(1)}, \laginf^{(2)} \dots \laginf^{(N)}\}$ is the concatenated Lagrange multipliers\footnote{In the following, we assume that $\min_{n, s_n}\mu^{(n)}_{s_n}=0$. This can be done without loss of generality, since adding a constant to all values of $\mu$ gives an equivalent valid multiplier vector.} corresponding to the fixed point constraints \eqref{EQ::WHITF2}. Here each $\laginf^{(n)}$ is an $\sspaceN{n}$ dimensional vector. 
    \end{itemize}
\end{definition}
Note that $\gstar$ is unique by construction but $\bYN^{\star}$ is does not have to be unique and $\laginf$ is never unique. In the remainder of the paper, we fix $\bYN^{\star}$ and $\laginf$ to any particular solution.  

The following lemma follows from arguments of \cite{yan2024}, note that since the problem \eqref{EQ::WHITF1}-\eqref{EQ::WHITF3} can be viewed as a relaxation on the action constraints, its solution upper bounds the value of \eqref{EQ::VOPTT}.
\begin{lemma}\label{LEM:VGBOUND}
    The gain is an upper bound on the maximum value that can be obtained: $\VoptNoN \leq \gstar$.

\end{lemma}

\subsection{A finite horizon Linear Program}\label{SUBSEC:FLP}

Given an initial state $\bs := \{s_{1}, s_{2} \dots s_{N}\}$, we denote by $\bxs{\bs} := \{\bx^{(1)}, \bx^{(2)} \dots \bx^{(N)}\} \in \Real{|\sspace|}$ the one hot encoded vector,  
\begin{equation*}
    \bx^{(1)}_s :=
    \begin{cases}
        1 &\text{if } s = s_1\\
        0 &\text{otherwise}
    \end{cases}
\end{equation*}
which represent the empirical distribution of the joint states over the arms. Next, we may extend this definition to a more general $\bx := \{\bx^{(1)}, \bx^{(2)} \dots \bx^{(N)}\}$, where each $\bx^{(i)}$ represents the state distribution of an arm $i$.
Finally, one may define the state-action variable $\bYn{n}{s}{a}{t}$ which represents the probability of arm $n$ being in state $s$ and choosing action $a$ at time $t$. We call the corresponding state-action vector $\bYVn{(n)}{t}$ with the concatenation of the $N$ vectors denoted by $\bYNt{t}$. Finally, we use $\bYVZn{n}{t}$ and $\bYVUn{n}{t}$ to denote the vectors in $\Real{|\sspaceN{n}|}$ corresponding to action $0$ and $1$ respectively for arm $n$. We now define a finite-horizon linear program for fixed horizon $\THor$ as follows:     
\begin{align}
    \VLPROB{\THor}{\bx} :=& \max_{\bYNt{t}: t \in \{0, 1 \dots \THor\}}\left[ \sum_{t = 0}^{\THor - 1} \frac{1}{N} \sum_{n = 1}^N \Rvec{n} \cdot\bYVn{(n)}{t} \right] + \langle \laginf, (\bYVZN{\THor} + \bYVUN{\THor}) \rangle\label{EQ::FLUIDLP1}\\
    &\text{such that,} \nonumber \\
    &\bYVZn{n}{0} + \bYVUn{n}{0} = \bx^{(n)} \hspace{0.2 in } \forall n  \label{EQ::FLUIDLP2} \\
    &\bYVZn{n}{t + 1} + \bYVUn{n}{t + 1} =  \bYVZn{n}{t} \cdot  \Pmat{n}{0} + \bYVUn{n}{t} \cdot \Pmat{n}{1} \hspace{0.2 in } \forall n, t \label{EQ::FLUIDLP3}\\
    &\sum_{n = 1}^{N} \| \bYVUn{n}{t}\|_1 \leq N \alpha \hspace{0.2 in } \forall t \label{EQ::FLUIDLP4}
\end{align}
Here, $\laginf$ is the multiplier that was specified in definition \ref{DEF:WHITMUL}. Constraint \eqref{EQ::FLUIDLP3} corresponds to the averaged Markov constraint. Further, we are relaxing the hard constraint on the actions of the original RMAB problem by only ensuring that they are solved in expectation \eqref{EQ::FLUIDLP4} at each time instance.  

Given any initial distribution $\bx$, we will refer to the set of all feasible state action vectors (i.e., satisfying condition \eqref{EQ::FLUIDLP2} - \eqref{EQ::FLUIDLP4}) by $\mcY{\bx}$. Now, one can restate the $\THor$ horizon problem using the following dynamic program:
\begin{align}\label{EQ::DYN1}
    \VLPROB{\THor}{\bx} =
    \begin{cases}
        &\max_{\bYNt{0} \in \mcY{\bxs{\bs}}}\frac{1}{N} \sum_{n = 1}^N  \Rvec{n} \cdot \bYVn{n}{0}  + \VLPROB{\THor - 1}{\bYVUN{0} \cdot \PmatNN{1} + \bYVZN{0} \cdot \PmatNN{0}}\\
        & \langle \laginf, \bxs{\bs} \rangle \hspace{0.3 in} \text{when $\THor = 0$}
    \end{cases}
\end{align}

\begin{definition}\label{DEF:FLUID}
   For any fixed initial state $\bx$ and computational horizon $\THor$, let $[\bYFLT{\bx (t), t'}]_{t' = t}^{t + \THor}$ be the solution to the dynamic program, \eqref{EQ::DYN1}. In keeping with the naming convention from \cite{DBert2016}, we will call $[\bYFLT{\bx (t), t'}]_{t' = t}^{t + \THor}$ the  \emph{aggregate flow} and $[\bYFLTn{\bx(t), t'}{n}]_{t' = t}^{t + \THor}$ the \emph{flow per arm}.  
\end{definition}

$\VLPROB{\THor}{\bx}$ is Lipschitz and concave. We elaborate further on the properties of the dynamic program with proof in \ref{APP::FLPROP}. We will refer to the first state-action vector $\bYFLT{t}$ as the \emph{fluid solution} to the finite horizon dynamic program with horizon $\THor$ to \eqref{EQ::DYN1}.

\subsection{The LP-update policy}\label{SUBSEC:ALG}

We will now use the fluid policy described in the previous subsection to inform our infinite-horizon policy. However, we run into a few problems in using the fluid policy to solve our problem, \textbf{1}) the policy is defined only for a finite horizon, \textbf{2}) how does one go about converting the continuous state-action vector defined for a \emph{deterministic fluid problem} into a binary decision variable for a \emph{stochastic problem }? We use a \emph{model predictive policy} with horizon $\THor$ to address the first problem and a \emph{randomized rounding procedure} to solve the second problem. 

Our algorithm at time $t$ can be summarized by the relation:
\begin{equation}\label{EQ::SUMREL}
    \bs(t) \xrightarrow[\text{encoding}]{\text{One-hot}} \bxs{\bs(t)} \xrightarrow[\eqref{EQ::FLUIDLP1}-\eqref{EQ::FLUIDLP4}]{\text{$\THor$-horizon LP}} \bYFLT{t} \xrightarrow[\text{rounding}]{\text{Randomized }} \bANN(t) \xrightarrow[]{\eqref{EQ:MKEVOL}} \bs(t + 1).
\end{equation}
More precisely, at time $t$, given a state $\bs(t)$, we use the one hot coded state $\bxs{\bs(t)}$ as the initial state to solve a $\THor$ horizon Linear program \eqref{EQ::FLUIDLP1}-\eqref{EQ::FLUIDLP4} to obtain a fluid solution $\bYFLT{t}$. We set the choice of horizon $\THor$ in line $4$ according to a function $\hat{\THor}(\bx, t)$ that could in general depend on the state $\bx$ of the system and time $t$. 
In our first variant, we choose $\hat{\THor}(\bx, t)$ to be a \textbf{fixed constant} as is traditional in model predictive control
.\footnote{We show in equation \eqref{EQN::INFCOST3} that one can choose a constant that is at most $\epsilon$ away in value from the optimal solution.} A second variant is shown in Section \ref{SUBSEC::MPC2} where the horizon changes as a function of time horizon to provide different guarantees for the algorithm. 

Denoting by $\bYFLT{t}$ the first vector of the solution of length $\THor$ to \eqref{EQ::DYN1}, we use a randomized rounding procedure generates a random binary decision vector for pulling each arm with marginal probability proportional to $\bYFLUn{(n)}{\bsn{n}}{t}$. 
Given the state and the joint actions, one may now observe the next set of states at time $t + 1$ since the arms will now evolve independently according to \eqref{EQ:MKEVOL}. We then return to step $1$ and solve a $\THor$ horizon LP and the process continues iteratively.


The following lemma, proven in \ref{APP:RR}, characterizes the $l^1$ norm distance between the action vector chosen by the rounding procedure (that is fully described in Algorithm~\ref{algo::rounding} in ~\ref{APP:RR} and is inspired by \cite{ioannidis2016adaptive}). 
With a slight abuse of notation, we let $\bANN := \{\bANN^{(1)}, \bANN^{(2)} \dots \bANN^{(N)}\}$ be the random joint action vector that uses one hot encoding to indicate which arms are being pulled.
\begin{lemma}\label{LEM:RR1}
    The $l^1$ norm distance between $\bYFLU{t, \bs(t)}$ and the $\bANN$ produced by Algorithm~\ref{algo::rounding} is bounded in expectation, given $\bsN(t)$:
    \begin{equation}
        \frac{1}{N}\sum_{n = 1}^{N}\|\E [\bANN^{(n)}(t)] - \bYFLUn{(n)}{\cdot}{t}\|_1  \leq \frac{\left( \alpha N - \lfloor \alpha N \rfloor \right)}{N} \leq \frac{1}{N}
    \end{equation}
\end{lemma}

We are now ready to introduce our LP-update algorithm, the explanation will follow.

\begin{algorithm}
\begin{algorithmic}[1]
    \caption{LP-update policy}
	\label{algo::MPC}
            \State \textbf{Input}:(Initial state $\bs(0)$, model parameters $\PmatN, \RvecN$,Total time $T$, Procedure $\hat{\THor}(\cdot)$)
            \State \textbf{Solve} the fixed point equation \eqref{EQ::WHITF1}-\eqref{EQ::WHITF3} to find $\mu$.
            \For{$t=0$ to $T-1$}
            \State Set: $\THor \leftarrow \hat{\THor}(\bx, t)$ 
            \State $\bYFLT{t} \leftarrow \text{Fluid solution} \hspace{0.05 in} \eqref{EQ::FLUIDLP1}\!\!-\!\!\eqref{EQ::FLUIDLP4}$ with initial condition $\bxs{\bs(t)}$ and time horizon $\tau$.
            \State $\ba (t)$ $\leftarrow$ Randomized Rounding($\bYFLT{t}$) (described in Algorithm~\ref{algo::rounding}). 
            \State Observe:  $R(\bs(t), \ba(t))$.
            \State Observe: $\bs(t + 1)$ according to transition \eqref{EQ:MKEVOL}.
            \State Set: Total-reward $\leftarrow$ Total-reward + $R(\bs(t), \ba(t))$
            \EndFor
            \State $\VLPN = \frac{\text{Total-reward}}{T}$
    \end{algorithmic}
\end{algorithm}


 As we shall show in the next section, our algorithm requires a mild ergodicity assumption to provide guarantees, we make these assumptions precise in \eqref{eq:ergodic_coef}. 

Having defined the algorithm, we denote the average \textbf{value function under our LP-update policy by $\VLPN$.} 

\begin{remark}
    It is important to note that $\bYFLT{t'}$ is a function of the initial state $\bs (t)$. Here, we have used $\bYFLT{t'}$ as a short-hand for $\bYFLT{\bxs{\bs(t)}, t'}$, suppressing the dependence on the initial state for brevity. The LP-update policy adopts only the first value in the collection of $\THor$ state-action vectors derived from the finite horizon process, justifying our notation. 
\end{remark}

\section{Main Theorem\\}\label{SEC::MAINRESULTS}

Our main result will show that the LP update policy is asymptotically optimal as the number of arms $N$ goes to infinity, under an easily verifiable mixing assumption on the transition matrices. To express this condition, for a fixed integer and a sequence of action $\ba=(a_1\dots a_k)\in\{0,1\}^k$, we denote by $P^{(n)}_{j|i, \ba}$ the $(i,j)^{\text{th}}$ entry of the matrix $\prod_{t=1}^k \Pmat{n}{a_t}$. We then denote by $\rho_k$ the the synchronization constant defined below: 
\begin{align}
    \label{eq:ergodic_coef}
      \syncconst{k} \triangleq \min_{n, s, s'\in\sspaceN{n}, \ba\in\{0,1\}^k} \sum_{s^* \in \sspaceN{n}}\min\{P^{(n)}_{s^*|s,(a_1, a_2, \dots a_k)}, P^{(n) }_{s^*|s', (0, 0, \dots 0)}\}
\end{align}
In the equation above, the minimum is taken over all arms, all possible initial states $s,s'$ and all possible sequence of actions over $k$ steps. The quantity $\syncconst{k}$ can be viewed as the probability (under the best coupling) that two arms starting in states $s$ and $s'$ reach the same state after $k$ iterations, if the sequence $a_1\dots a_k$ is used for the first arm while the second arm only uses the action $0$. The assumption that we use for our result is that $\syncconst{k}>0$ for some integer $k$.
\begin{assumption}\label{AS::SA}
    There exists a finite $k$ such that $\syncconst{k} > 0$.
\end{assumption}

\begin{remark}
    We pause to state that while our assumption may look complicated at first, a sufficient condition for the assumption to be satisfied is if the transition kernels of all arms are ergodic under action $0$. This follows since if $\Pmat{n}{0}$ is ergodic, then there exists $k>0$ such that $\left[\Pmat{n}{0}\right]^k > 0$ for all entries. This provides a simple sufficient condition that implies Assumption~\ref{AS::SA}.
\end{remark}
We now come to the main theorem. We show that the LP-update policy is asymptotically optimal as the number of arms tends to infinity. 
\begin{theorem}\label{THM::ALMOSTOPT}
    Under Assumption \ref{AS::SA} with constant $\syncconst{k}$ for some fixed integer $k$, fix a positive constant $\epsilon > 0$; then there exists a constant $\THor(\epsilon)$ such that setting the planning horizon to $\hat{\THor}(\bx, t) = \THor(\epsilon)$ for the LP-update policy we have the following performance guarantee:
    \begin{align}
        \gstar \geq \VLPN \geq~& \gstar - \epsilon - \frac{\THor(\epsilon)}{2} \Bigg(\Bigg(1 + \frac{k}{\syncconst{k}} \Bigg)\frac{k}{\syncconst{k}}  + \|\laginf\|_{\infty} + 1\Bigg)\sqrt{\frac{\pi}{N}} \nonumber\\
    &- \left[1 + 2\|\mu\|_{\infty} + \left(1 + \frac{k}{\syncconst{k}} \right)\frac{k}{\syncconst{k}}\right] \frac{\alpha N - \lfloor \alpha N \rfloor}{N}.
    \end{align}
\end{theorem} 

 There are three error terms that determine the upper bound between the difference in performance of the LP-update policy and the optimal policy. 
 The first error term appears since we are using a finite time horizon LP to solve an infinite horizon problem. We bound this term using the idea of \textbf{dissipativity}. The last error term occurs due to the difference between choosing a \textbf{fluid policy} and its implementation using \textbf{randomized rounding}. Our LP-policy solves a relaxed problem where the constraints are satisfied in expectation \textit{at each time step}, translating this solution to a viable joint action gives rise to this error term. The remaining error term is derived from a \textbf{Jensen gap} result (in keeping with the naming convention of ~\cite{Brown&Zhang25}). Due to the heterogeneity of our system we are unable to directly leverage concentration results such as the ones used \cite[Lemma 1]{GGY23b}. We need to show that the stochastic value of the next state does not stray too far from its deterministic counterpart. Our difficulty lies in the fact that the value function is concave. Hence, loosely speaking this first term stems from an `inverse' Jensen inequality.

\subsection{Refining the Main Result.\\}\label{SUBSEC::REF}

The main result allows us to distill the difference between the optimal value that can be obtained and the value that is obtained as an outcome of our algorithm into error components. We will dedicate this subsection to refining the result in Theorem \ref{THM::ALMOSTOPT} by examining the dependence of the error components on the underlying parameters and the role of our time horizon under additional assumptions.  In light of Remark \ref{REM::TURNPIKE}, we conjecture the following about the nature of the necessary computation time $\THor$.

\begin{conjecture}\label{CONJ:1}
    Given a collection of \(N\) heterogeneous arms that satisfy the Assumption \ref{AS::SA}, the function $\THor(\epsilon): \Real{} \to \mathbb{N}$ is independent of \(N\) and depends only on the maximum size of the state space, $\max_{n}\{|\sspaceN{n}|\}$, the budget $\budget$ and the ergodicity constant $(k$, $\syncconst{k})$. 
\end{conjecture}

The first result in refining Theorem \ref{THM::ALMOSTOPT} will be to bound the Lagrange multiplier $\laginf$. This result is proved explicitly in \cite{gast2024MPCoptimal} [Corollary 15], hence, we will state it here without a proof. 

\begin{lemma}\label{LEM:LAGBOUND}
    The lagrange multiplier is bounded in the infinite norm as follows,
    \begin{equation}
        \|\laginf\|_{\infty} \leq \frac{k}{\syncconst{k}}\left[1 + \frac{\alpha k}{\syncconst{k}} \right]
    \end{equation}
\end{lemma}

The following corollary takes into account the bound on the Lagrange multiplier Lemma \ref{LEM:LAGBOUND}, the relationship between $\VoptNoN$ and $\gstar$ in Lemma \ref{LEM:VGBOUND} and Conjecture \ref{CONJ:1} to refine Theorem \ref{THM::ALMOSTOPT}:

\begin{corollary}\label{COR::REFRES}
    Fix $\epsilon > 0$, under Assumption \ref{AS::SA} with constant $(k, \syncconst{k})$ and Conjecture \ref{CONJ:1}, there exists a fixed planning horizon $\THor(\epsilon)$ that is independent of the state such that the LP-update algorithm 
    has the following guarantee of performance:
    \begin{equation}
        0 \leq \VoptNoN - \VLPN \leq \epsilon + \mcl{O}\left(\frac{\THor(\varepsilon)}{\sqrt{N}}\right).
    \end{equation}
    Where $\mcl{O}$ hides a $C(\budget, k, \syncconst{k}) $ factor that can be computed explicitly using Lemma \ref{LEM:LAGBOUND}. 
\end{corollary}


Before we end this section we make the following important remark:
\begin{remark}\label{REM::TURNPIKE}
    One major problem that we face in our algorithm are guarantees regarding the time horizon as a function of the error $\epsilon$. We note that in general the time horizon can scale with the dimension of $\sspace$ which in itself depends on $N$ in the fully heterogeneous case. To put the reader's mind at ease, we remark that in all our simulations a time horizon of $\hat{\THor}(\bx, t) = 5$ for the LP-update policy suffices to solve most problems even in adversarial settings. 
    If the gap between $\VLPN$ and the optimal value decreases geometrically with horizon Conjecture $1$ allows us to state that $\tau(\epsilon)=O(\log(\epsilon))$. This phenomenon is called the \textbf{exponential turnpike property} and is beyond the scope of the current paper and left for future work.
\end{remark}

\subsection{Generalization to weakly coupled MDPs.\\}\label{SUBSEC::WMDP}

In this section, we show that we can generalize the main result to the case of weakly coupled MDPs. Consider a variant of our system model, where we have $N$ Markov processes that are weakly coupled through their actions,
$\langle \sspace, \mcl{A}, \{\Pmat{n}{a}\}_{n = 1, a \in \mcl{A}}^{n = N}, \Rvec{n}; \budgetK \rangle$.  As before, the decision maker must control these $N$ heterogeneous processes which operate on separate state spaces $\sspaceN{n}$ and a common finite action space $\mcl{A}$. The reward of arm $n$ in state $s$ under action $a$ is given by $\Rvec{n}(s, a)$ and the arm transitions according to the transition kernel $\Pmat{n}{a}$. Given a joint state $\bSN$ and action $\bAN$ the arms evolve according to \eqref{EQ:MKEVOL}.  

The main generalization of the restless multi-armed bandit problem to the weakly coupled Markov decision process (WMDPs) comes in the form of the generalized action constraint set, 
\begin{align*}
    \sum_{n = 1}^{N} \constrWMDP{(n, e)}{S_n(t)}{A_n(t)} \leq N\budgetk{e} \hspace{0.2 in} \text{for all $e$ and $t$}
\end{align*}
where $\budgetK := \{\budgetk{e}\}_{e = 1}^{E}$ is a vector of dimension $\Real{E}$ that specifies the budget constraints for our actions. Where we use $\constrWMDP{(n, e)}{\cdot}{\cdot}$ to denote the $e^{\text{th}}$ cost function. We assume here that each of our constraints are non-negative and bounded above by a constant $b$ independent of $N$, i.e. $0 \leq \constrWMDP{(n, e)}{\cdot}{\cdot} \leq b$.

Let $\pspaceN{N}$ denote the space of all policies that satisfy this constraint. Given a policy $\pi \in \pspaceN{N}$, and an initial set of arm states $\bSN(0) = \bs$, we can once again define the average reward of the policy $\pi$ by:
\begin{equation*}
    \VpolavgN{\pi}{N}{\bs} = \lim_{T \to \infty}\frac{1}{T}\E_{\pi}\left[ \sum_{t = 0}^{T - 1} \frac{1}{N} \sum_{n = 1}^N \rew{n}{A_n(t)}{S_n(t)}\Bigg|\bSN(0) = \bs \right].
\end{equation*}
The heterogeneous infinite-horizon average reward problem for the weakly coupled Markov decision process is to find the optimal policy $\pi \in \pspaceN{N}$ or:
\begin{align}\label{EQ::VOPTWMDP}
    \VoptWMDP{N} &:= \max_{\pi \in \pspaceN{N}} \VpolavgN{\pi}{N}{\bs}.
\end{align}
As before, given an initial state $\bx$, the corresponding finite-horizon linear program is given by:
\begin{align*}
    \VLWMDP{\THor}{\bx} :=& \max_{\bYNt{t}: t \in \{0, 1 \dots \THor\}}\left[ \sum_{t = 0}^{\THor - 1} \frac{1}{N} \sum_{n = 1}^N \Rvec{n} \cdot\bYVn{(n)}{t} \right] + \langle \laginf, \sum_{a \in \mcl{A}} \bYWMDPn{\THor}{} \rangle\\
    &\text{such that,} \nonumber \\
    &\sum_{a \in \mcl{A}} \bYWMDPn{0}{(n)} = \bx_{(n)} \hspace{0.2 in } \forall n  \\
    &\sum_{a \in \mcl{A}} \bYWMDPn{t + 1}{(n)} = \sum_{{a \in \mcl{A}}}\Pmat{n}{a} \cdot \bYWMDPn{t}{(n)}   \hspace{0.2 in } \forall n, t \\
    &\sum_{n = 1}^{N} \sum_{(s,a)}\bYWMDPn{t, s}{n}\constrWMDP{(n), (e)}{s}{a} \leq N \budgetk{e} \hspace{0.2 in} \text{for all $e$ and $t$}
\end{align*}

One difficulty that arises when using the LP-update rule for the weakly coupled problem rather than the markovian bandit problem is the use of randomized rounding to draw $A_n(t) = a$ from the fluid solution $\bYFLT{t}$ in order to construct a feasible joint action. This difficulty can easily be circumvented by drawing an action with probability proportional to $\bYFLAn{t}{a}{n}$, ties being broken arbitrarily. 

Next, when considering the weakly coupled problem, one must re-examine the assumptions under which feasible solutions exist. A natural extension of the Ergodicity Assumption \ref{AS::SA} is stated in terms of the next two assumptions. 

\begin{assumption}[Feasibility]\label{AS:SA2}
There exists an action $a^{\star} \in \mcl{A}$ that is \emph{always feasible} in the following sense: Given a joint state $\bsN := \{s_{1}, s_{2} \dots s_{N}\}$ let $\bAN := \{a_{1}, a_{2} \dots a_{i} \dots a_{N}\}$ be a joint action that satisfies the constraints, i.e, 
\[
\sum_{n = 1}^{N} \constrWMDP{(n),(e)}{s_{n}}{a_{n}} \leq N \budgetk{e}
\]
for all $e$. Then replacing an action $a_{i}$ by $a^{\star}$ for any $i$ will continue to satisfy the constraints,
\[
\sum_{n \neq i}\constrWMDP{(n), (e)}{s_{n}}{a_{n}} +  \constrWMDP{(i), (e)}{s_{i}}{a^{\star}}\leq N \budgetk{e}
\]
\end{assumption}
Note, this is a natural generalization of the $0-$ action in the RMAB case where choosing the \(0-\) action for more arms will never violate the budget constraint. Here we have simply replaced the $0$ action with the $a^{\star}$ action and the \(\mathbf{1}_{a_{i} = 1}\) cost with the more general cost constraint, $\constrWMDP{(i), (e)}{s_{i}}{a_{i}}$.

We will extend our idea of Ergodicity to the WMDP using the aforementioned action $a^{\star}$. 
\begin{assumption}[Ergodicity]\label{AS:SA3}
     Consider the joint action $\bAN(t) := \{a^{\star}, a^{\star} \dots \}$ where $a^{\star} \in \mcl{A}$ is described in Assumption \ref{AS:SA2}. This joint action satisfies the following generalization of inequality \eqref{eq:ergodic_coef},
    \begin{align}
    \label{eq:ergodic_WMDP}
      \syncconst{k} \triangleq \min_{n, s, s'\in\sspaceN{n}, \{a_1, a_2 \dots a_k\}\in\mcl{A}^k} \sum_{s_* \in \sspaceN{n}}\min\{P^{(n)}_{s_*|s, a_1, a_2, \dots a_k}, P^{(n)}_{s_*|s',a^{\star}, a^{\star}, \dots a^{\star}}\} > 0   
\end{align}
for some fixed integer $k$. 
\end{assumption}

Noting that the leading term for our error comes from the Jensen Gap term, we obtain the following corollary by following the same steps for the proof as Theorem \ref{THM::ALMOSTOPT}:

\begin{corollary}\label{COR::GENMDP}
    Under Assumption \ref{AS:SA2}, Assumption \ref{AS:SA3} and Conjecture \ref{CONJ:1} there exists a planning horizon $\THor(\epsilon)$ such that the LP-update algorithm has the following guarantee of performance:
    \begin{equation}
        \VoptWMDP{N} - \VLPinfWMDP{\THor} \leq \epsilon + \mcl{O}\left(\frac{\THor(\epsilon)}{\sqrt{N}}\right).
    \end{equation}
    Where $\mcl{O}$ hides a constant factor $C(\budget, k, \syncconst{k})$ which depends only on the $\budget$, $k$ and $\syncconst{k}$.
\end{corollary}


The proof of this corollary is parallel to the proof for Theorem \ref{THM::ALMOSTOPT} and can be extrapolated similar to the proof of Corollary \ref{COR::REFRES}. We will differ a sketch of the proof to \ref{APP:GENMDP} since it requires the framework presented in Section \ref{SEC::PROOF}. 

\subsection{Modified Algorithm with Time Horizon Guarantees}\label{SUBSEC::MPC2}

As mentioned in Remark \ref{REM::TURNPIKE}, the relationship between $\THor(\epsilon)$ and $\epsilon$ is unclear based on the method of our proof. This poses a problem for a designer since a long horizon greatly increases the amount of time the algorithm takes for computation. This is because unlike a priority policy, the LP-update must be computed at each step. We now present a modification to the LP-update algorithm that allows us to bridge this gap. Let $B:=\left[\|\laginf\|_{\infty} +  \left(1 + \frac{k}{\syncconst{k}} \right) \frac{k}{\syncconst{k}}\right] + \|\laginf\|_{\infty}$, the main idea of this modified algorithm is to cut our time into segments of length $B/\epsilon$. The value of $B$ is a bound on an equivalent cost minimization problem for the fluid problem, more on this topic can be found in our Section on dissipativity, \ref{SUBSEC::DISSIPATIVITY}. 

 Concretely, we treat the problem as a finite horizon problem where we can use recomputation of the fluid solution to bound our error. More specifically, at time $t$ in $\{1, 2 \dots B/\epsilon\}$ with joint state $\bs(t)$, we compute a fluid solution with initial state $\bs(t)$ and length $B/\epsilon - t \mod B/\epsilon$. As before, we pick the first solution produced in this process and apply it to move a single step forward. It is not hard to check that our main error in this case once again stems from ``Jensen's gap". We continue this process till the end of the finite time segment. When our time horizon hits $0$ we reset the time horizon once more to $B/\epsilon$ and begin again the process of computing a fluid solution for a finite horizon problem.We denote the output of our finite segment algorithm by the value function $\VLPF$.

Our algorithm has the following guarantees that accounts for a relationship between $\epsilon$ and the time horizon $\THor$.

\begin{theorem}\label{THM::ALMOSTOPT2}
    Under Assumption \ref{AS::SA} with constant $\syncconst{k}$ for some fixed integer $k$, fix a positive constant $\epsilon > 0$; then 
    one can choose the time horizon $\hat{\THor}(\bx, t) := B/\epsilon - t \mod B/\epsilon$ in line 4 of the LP-update policy (Algorithm \ref{algo::MPC}) to obtain the following guarantee:
    \begin{align}
        \gstar \geq \VLPF \geq& \gstar - \epsilon - \frac{B}{\epsilon} \Bigg(\Bigg(1 + \frac{k}{\syncconst{k}} \Bigg)\frac{k}{\syncconst{k}}  + \|\laginf\|_{\infty} + 1\Bigg)\sqrt{\frac{\pi}{N}} \nonumber\\
    &- \left[1 + 2\|\mu\|_{\infty} + \left(1 + \frac{k}{\syncconst{k}} \right)\frac{k}{\syncconst{k}}\right] \frac{\alpha N - \lfloor \alpha N \rfloor}{N}.
    \end{align}
    Where the constant $B$ is equal to $\left[\|\laginf\|_{\infty} +  \left(1 + \frac{k}{\syncconst{k}} \right) \frac{k}{\syncconst{k}}\right] + \|\laginf\|_{\infty}.$
\end{theorem} 

\begin{remark}
    We remark that several works have looked at finite horizon variants of the LP-update policy with recomputation, recently \cite{Brown&Zhang25} and \cite{GGY23b}. The crucial difference is knowing the exact value of the horizon $\hat{\THor}(\bx, t)$ to choose so that this strategy can be directly used to break an infinite horizon problem into finite horizon problems of fixed length. We refer the reader to Remark \ref{REM::DISS} for more details. It worth noting that both variants of our algorithm are independent of the state $\bx$ however, our first variant could scale arbitrarily with the dimension $N$ which necessitates the conjecture. We pay a price for avoiding this conjecture in the second variant by noting that the best scaling we can obtain is $\mathcal{O}(\frac{1}{N^{1/4}})$. 
\end{remark}






\section{Technical Framework for the Main Theorem}\label{SEC::PROOF}
Several works have viewed the problem of using the LP fluid policy to solve a finite horizon heterogeneous restless multi-armed bandit problem. One of the main difficulties associated with extending the finite horizon policy to its infinite horizon average reward counterpart is that the difference in the value of solutions can grow exponentially over time, see for example the value of the Lipschitz constant in \cite{GGY23b} as a function of time. As a result, the optimal solution induced by the finite horizon problem may have nothing to do with the infinite horizon counterpart. Another approach that is quite common in analytical results on Model Predictive Control is first solving a problem using the discounted model and then allowing the discount factor to converge towards $1$. The reader is directed towards [Theorem 4.5, \cite{Brown&Zhang25}] for an illustration of this method. Unfortunately, we cannot apply this method to obtain our results. Typically, the idea in discounted reward problems is to use truncation of the time horizon, unfortunately as seen in [Theorem 4.5, \cite{Brown&Zhang25}], for discount factor $\beta$, the horizon for truncation is order $\log_{1/\beta} \sqrt{N}$ which very quickly goes to infinity as the discount factor goes to $1$. We are therefore forced to develop new tools to tackle the infinite horizon average reward problem. A more acute version of this problem is seen in \cite{brown2020index} without recomputation in the fully heterogeneous scenario.

\emph{Dissipativity} is a key feature of our framework which provides us with an equivalent surrogate loss function. We introduce this framework in Section \ref{SUBSEC::DISSIPATIVITY}. As a result, the finite horizon problem cannot be too far away from the infinite horizon limit in value and hence, one is free to use a model predictive control approach to solve the average reward problem. As a bonus, this framework readily allows us to analyze the second variant of Algorithm \ref{algo::MPC} where we cut the horizon into finite time intervals to give us explicit bounds on optimality. We relegate our proof of this result to \ref{APP::THM2}. 

A second key challenge is to show that under heterogeneity we are still able to obtain an approximation result which allows us to use a deterministic fluid system to solve a stochastic problem. We address this challenge using the \emph{Jensen's gap lemma}, introduced in \ref{SUBSEC:JENSENGAP}. 

We begin our framework by discussing an important continuity result for the bias function. 

\subsection{Continuity of value and bias functions}\label{SUBSEC::BIASCONT}

This section covers two key continuity properties that we will need throughout our proof. These properties will prove important for bounding differences between value at a state and its fluid counter-part.
\begin{restatable}[Continuity of the Fluid Value Function]{lemma}{lemValC}
\label{LEM::VALT}
    Let $\bx := \bxc{i}$ and $\Tilde{\bx} := \bxcn{i}$ be two joint state distributions that only differ in the $i^{\text{th}}$ component. Then under Assumption~\ref{AS::SA}, with synchronization constant $\syncconst{k} > 0$ for some positive integer $k$:
    \begin{equation}
        \VLPROB{t}{\bx} - \VLPROB{t}{\Tilde{\bx}} \leq \frac{1}{N}\left[\|\mu\|_{\infty} + \left(1 + \frac{k}{\syncconst{k}} \right)\sum_{l = 1}^{t}\left(1 - \frac{\syncconst{k}}{k} \right)^{l - 1}\right]\|\bx^{(i)} - \Tilde{\bx}^{(i)}\|_1
    \end{equation}
\end{restatable}
We differ the proof of this lemma and of other properties of the fluid problem to \ref{APP::FLPROP}.



Almost as a direct consequence of Lemma \ref{LEM::VALT}, we show that we can define a \textbf{bias function} $\hstar{\cdot} : \Pi_{n = 1}^{N}\Delta(\Real{|\sspaceN{n}|}) \to \Real{1}$ that is  Lipschitz-continuous and convex (see the proof 
in \ref{APP:CONTWFL}).

\begin{restatable}[Continuity of the Bias Function]{proposition}{Contbias}\label{PROP::BIAS_CONT}
    For any initial condition $\bx$, let $\hstar{\bx}$ be:
    \begin{equation}\label{EQ::BIASFN}
        \hstar{\bx} := \lim_{T \to \infty}T \gstar - \VLPROB{T}{\bx}.
    \end{equation}
    Under Assumption \ref{AS::SA}, the bias function $\hstar{\cdot}$ is well defined, convex and Lipschitz-continuous in the following sense: let $\bx := \bxc{i}$ and $\Tilde{\bx} := \bxcn{i}$ be two state distribution, we have
    \begin{align}\label{EQ::HCONT}
        \hstar{\bx} - \hstar{\Tilde{\bx}} \leq \frac1N \left[\|\laginf\|_{\infty} + \left(1 + \frac{k}{\syncconst{k}} \right)\frac{k}{\syncconst{k}}\right] \|\bx^{(i)} - \Tilde{\bx}^{(i)} \|_1
    \end{align}
\end{restatable} 

\subsection{Dissipativity}\label{SUBSEC::DISSIPATIVITY}

In the previous section we showed that $\hstar{\cdot}$ is a well defined and continuous function. This allows us to define the notion of dissipativity, which in turn will allow us to define our surrogate loss function for the \emph{fluid problem}, $\CostT{\bx}{\infty}$ in terms of $\laginf$ and $\hstar{\cdot}$. We have already described the most relevant properties of $\laginf$ in section \ref{SEC::IPL}. We will give an equivalent set of definitions where we delineate between the state and control variable instead of treating the problem as an LP for the state-action variable as is more traditional in the control literature. Given a state $\bx$, one may define a one-to-one relationship between $\bYN$ and control $\bu$ simply by setting $\bu = \bYN_{\cdot, 1}$ and noting the constraint \eqref{EQ::FLUIDLP2}, which gives us $\bx = \bYN_{\cdot, 0} + \bYN_{\cdot, 1}$. Consequently, the transition kernel given $\bx$ and $\bu$ is defined by $$\bx_{\text{next}} = \bx \cdot \PmatNN{0}  + \bu \cdot (\PmatNN{1} - \PmatNN{0})  =: \Phi(\bx, \bu).$$ 
We now use the following general definitions from \cite{DTGLS14}, for a discrete time \emph{deterministic} minimization problem $\CostT{\bx}{T} := \min_{\bu(t)}\sum_{t = 0}^{T - 1}\costt{\bx(t)}{\bu(t)}$,starting at state $\bx(0) = \bx$ with state $\bx(t)$ and control $\bu(t)$. The state transitions according to transition kernel $\Phi$ and we let the stage cost at time $t$ be given by $\costt{\bx(t)}{\bu(t)}$. 

The storage function and rotated cost are defined as follows:
\begin{definition}
    A \textbf{rotated cost} denoted by $\rotcost{c}{\bx,\bu}$ for the dynamical system $\Phi(\cdot, \cdot)$ with stage cost $\costt{\cdot}{\cdot}$ with \textbf{storage function} $\stor{\bx}$ is defined by:
    \begin{equation}\label{EQ::ROTCOSTDEF}
        \rotcost{c}{\bx,\bu} := \costt{\bx}{\bu} + \stor{\bx} - \stor{\Phi(\bx, \bu)}
    \end{equation}
\end{definition}
 Let $\bx_e, \bu_e$ be a fixed point for the dynamical system, i.e,
\[
\bx_e := \Phi(\bx_e, \bu_e)
\]
We say that a problem is \textbf{dissipative with respect to the fixed point $({\bx_e},{\bu_e})$} if the rotated cost satisfies, for any $(\bx, \bu)$:
\begin{equation}\label{EQ::DISSDEF}
    \rotcost{c}{\bx,\bu} \geq \rotcost{c}{\bx_e,\bu_e}.
\end{equation}
Without loss of generality we set the rotated cost to $0$, $\rotcost{c}{\bx_e, \bu_e} = 0$. 

Recall, $\bYN^{\star}$ is the fixed point solution of \eqref{EQ::WHITF1} - \eqref{EQ::WHITF3} from Definition \ref{DEF:GSTAR}. We set $(\bx^{\star}, \bu^{\star})$ to be the corresponding fixed point in terms of the state and control variable. 

With the definitions out of the way we can now state the following result for our system:
\begin{lemma}\label{LEMM::BIAS}
    By setting the stage cost to 
    \[
    \costt{\bx}{\bu} := \gstar - \left(\langle \Rvec{0}, \bx \rangle + \langle \Rvec{1} - \Rvec{0}, \bu \rangle \right)
    \]
    and storage function,
    \[
    \stor{\bx} := \langle \laginf, \bx \rangle
    \]
    the rotated cost becomes
    \begin{align}
        \rotcost{c}{\bx, \bu} =& \gstar - \left(\langle \Rvec{0} \cdot \bx \rangle + \langle \Rvec{1} - \Rvec{0}, \bu \rangle \right) + \langle \laginf, \bx - \Phi(\bx, \bu) \rangle. \nonumber
    \end{align}
    The rotated cost function so defined is dissipative with respect to the fixed point of \eqref{EQ::WHITF1} - \eqref{EQ::WHITF3}, $(\bx^{\star}, \bu^{\star})$. 
\end{lemma}
We define the following important \textbf{surrogate loss function} for our LP-update policy.

\begin{definition}\label{DEF::SURCST}
   We define the following cost minimization problem under constraints \eqref{EQ::FLUIDLP2} -\eqref{EQ::FLUIDLP4} with fixed time horizon $\THor$ and initial state $\bx(0) = \bx$,
    \begin{equation}\label{EQN::COSTMINDEF}
        \CostT{\bx}{\THor} := \min_{\bYN s.t. \eqref{EQ::FLUIDLP2} -\eqref{EQ::FLUIDLP4}}\sum_{t = 0}^{\THor - 1}\rotcost{c}{\bx, \bu}.
    \end{equation}
    $\CostT{\bx}{\THor}$ is the \textbf{surrogate loss function of length $\THor$. }
\end{definition}
An immediate consequence of dissipativity is the following lemma.

\begin{lemma}\label{LEM::INFCOST}
    There exists a fixed time horizon $\THor(\epsilon)$ independent of the state $\bx$ such that,
    \begin{equation}\label{EQN::INFCOST3}
        |\CostT{\bx}{\THor(\varepsilon) + 1} - \CostT{\bx}{\THor(\varepsilon)}| \leq \epsilon.
    \end{equation}
\end{lemma}
The full proof of this theorem can be found in  \ref{APP:DISS}. We present below a sketch of the proof.
\begin{proof}
    The main idea of the proof is to note that the surrogate loss function, \eqref{EQN::COSTMINDEF} is equivalent to the reward maximization problem \eqref{EQ::FLUIDLP1}-\eqref{EQ::FLUIDLP4}. Further, as a consequence of dissipativity, the cost is non-decreasing in $\THor$ and the limit as $\THor \to \infty$ is given by:
    \begin{equation}\label{EQN::INFCOST1}
        \CostT{\bx}{\infty} := \langle \laginf, \bx \rangle + \hstar{\bx}.
    \end{equation}
    The rest of the proof is a straightforward consequence of monotonicity.
\end{proof}

This rotated cost will therefore allow us to demonstrate that a finite horizon solution to the LP is sufficient to determine a near optimal value for the infinite horizon reward maximization problem.

\begin{remark}\label{REM::DISS}
 The bound on the surrogate loss function, equation \eqref{EQN::INFCOST1} is exactly the value of $B$ chosen in Theorem \ref{THM::ALMOSTOPT2} in order to find an explicit finite horizon for the second variant of our algorithm. Dissipativity proves vital for the analysis of both variants of our LP-update strategy, the fixed time horizon as well as a changing time horizon.     
\end{remark}

\subsection{Jensen Gap lemma}\label{SUBSEC:JENSENGAP}

 Let $\bx := \{\bx^{(1)}, \bx^{(2)} \dots \bx^{(N)}\}$ be an initial state defined on $\Pi_{n = 1}^{N}\Delta(\Real{\sspaceN{n}})$. Next, let $X := \{X^{(1)}, X^{(2)} \dots X^{(N)}\}$ be a state so that each component is $X^{(i)}$ is sampled independently from the distribution $\bx^{(i)}$. In \ref{APP::FLPROP}, Lemma~\ref{APP::LEMCONCAVE}, we show that  $\VLPROB{t}{\bx}$ is concave in $\bx$. It follows from Jensen's lemma that $\left(\VLPROB{t}{\bx} - \E_{X \sim \bx} \left[ \VLPROB{t}{X}\right]\right) \geq 0$. We call this non-negative difference between a concave function and its expected value the \textbf{Jensen's gap}. We will show a `converse' result that indicates that Jensen's gap is, in fact, quite small for our problem indicating that our value is ``roughly" linear.  
 This lemma will prove essential for our concentration bound in the presence of heterogeneity. 
\begin{restatable}[Jensen Gap lemma]{lemma}{LemJen}
\label{LEM::JENGAP}
 The Jensen gap between $\VLPROB{t}{\bx}$ and the expected value of $\VLPROB{t}{X}$ where the state of each arm $i$ is drawn independently according to $\bx^{(i)}$ for any fixed $t$ is bounded by 
 \begin{equation}
     \VLPROB{t}{\bx} - \E_{X \sim \bx} \left[\VLPROB{t}{X}\right] \leq  \frac{t}{2}\left( 1 + \|\mu\|_{\infty} + \left(1 + \frac{k}{\syncconst{k}} \right)\frac{k}{\syncconst{k}} \right)\sqrt{\frac{\pi}{N}}.
 \end{equation} 
\end{restatable}
Our proof of this lemma proceeds by induction. There are two related technical challenges associated with the proof, first, the solution to $\VLPROB{t}{\bx}$ is a solution to a \emph{fluid problem}. Since, $X$ is sampled according to $\bx$, we need to generate a feasible joint action over the arms using our solution to $\bx$, in state $X$. Second, we must ensure that the value of $\VLPROB{t}{X}$ so produced are not too different from the values produced at $\bx$ for one more step. We then use steps one and two combined with the inductive hypothesis to complete the proof. It is important to note that if the feasible action is chosen incorrectly, the inductive hypothesis could potentially lead to an exponential dependence on the horizon $t$. This lemma is our most technically challenging result.
Our proof can be found in \ref{APP:JENGAP}.
\subsection{Proof of Theorem \ref{THM::ALMOSTOPT}}\label{SUBSEC::STEPSOFPROOF}
 Throughout the proof, the fluid problem computes the $\THor$ horizon LP problem and picks the first vector from the fluid problem for each arm $n$, denoted by $\bYFLTn{t}{n}$. The corresponding joint action using $\bYFLTn{t}{n}$ is denoted by $\bANN(t)$. Recall, the distribution of the next state under $\bYFLTn{t}{n}$ is $\bYFLTn{t}{n} \cdot \Pmat{n}{}$.
  We will use $\bx(t)$ to denote $\bxs{\bs(t)}$, the one hot coded empirical distribution given joint state $\bs(t)$. It now follows that, 
 \begin{align*}
    \gstar - &\VLPN =  \lim_{T \to \infty} \E \left[\gstar  - \frac{1}{T} \sum_{t = 0}^{T - 1} \frac{1}{N} \sum_{n = 1}^{N}\rew{n}{A_n(t)}{S_n(t)}\right]\\
     &= \lim_{T \to \infty} \E \left[\gstar - \frac{1}{T} \sum_{t = 0}^{T - 1} \frac{1}{N} \sum_{n = 1}^{N}\Rvec{n} \cdot \bYFLTn{t}{n} \right]\\
     &\hspace{0.4 in}+ \frac{1}{T} \E \Bigg[\sum_{t = 0}^{T - 1} \frac{1}{N} \sum_{n = 1}^{N} \left( \Rvec{n} \cdot \bYFLTn{t}{n} - \rew{n}{A_n(t)}{S_n(t)}\right)\Bigg]\\
     &= \lim_{T \to \infty} \E \left[\gstar - \frac{1}{T} \sum_{t = 0}^{T - 1} \frac{1}{N} \sum_{n = 1}^{N} \Rvec{n} \cdot \bYFLTn{t}{n}  + \langle \laginf, \bxs{t} -  \bYFLT{t} \cdot \PmatN \rangle\right]\\
    &\quad+ \frac{1}{T} \E\sum_{t = 0}^{T - 1} \frac{1}{N} \sum_{n = 1}^{N} \left( \Rvec{n} \cdot \bYFLTn{t}{n} - \rew{n}{A_n(t)}{S_n(t)}\right) - \frac{1}{T} \E\sum_{t = 0}^{T - 1} \langle \laginf, \bxs{t} -  \bYFLT{t} \cdot \PmatN \rangle\\
     &\stackrel{(a)}{\leq} \lim_{T \to \infty} \E \frac{1}{T}  \sum_{t = 0}^{T - 1} \rotcost{c}{\bxs{t}, \ufl{t}} + \frac{1}{T} \E\sum_{t = 0}^{T - 1} \frac{1}{N} \sum_{n = 1}^{N} \left( \Rvec{n} \cdot \bYFLTn{t}{n} - \rew{n}{A_n(t)}{S_n(t)}\right)\\
     &\hspace{0.4 in}+ \|\laginf\|_{\infty} \left(\frac{\alpha N - \lfloor \alpha N \rfloor}{N} \right) + \frac{\langle \mu, \bx  - \bYFLT{T - 1} \cdot \PmatN\rangle}{T} \\
     &\stackrel{(b)}{\leq} \lim_{T \to \infty} \epsilon + \E \frac{1}{T}  \sum_{t = 0}^{T - 1} \left[\CostT{\bxs{t + 1}}{\THor} - \CostT{\bYFLT{t} \cdot \PmatN}{\THor}\right] + \left(1 + \|\laginf\|_{\infty} \right)\frac{\alpha N - \lfloor \alpha N \rfloor}{N} + O\left(\frac{1}{T}\right)
 \end{align*}
 Step $(a)$ follows from the definition of the rotated cost and our choice of constructing Binary actions using \emph{randomized rounding}. We have also used an equivalent notation of $\bxs{t}, \ufl{t}$ in place of $\bYFLT{t}$. Step $(b)$ follows from rewriting the reward as an inner product between the action vector and the reward vector, the dynamic programming principle, and dissipativity. 
 
 Substituting the value of $\CostT{\cdot}{\THor}$ in terms of $\VLPROB{\THor}{\cdot}$, this is equal to
 \begin{align*}   
     &\stackrel{(c)}{=} \lim_{T \to \infty} \epsilon + \frac{1}{T}\E \sum_{t = 0}^{T - 1} \left[\VLPROB{\THor}{\bYFLT{t} \cdot \PmatN} - \VLPROB{\THor}{\bxs{t + 1}} \right]
     + \left(1 + \|\laginf\|_{\infty}\right)\frac{\alpha N - \lfloor \alpha N \rfloor}{N} + O\left(\frac{1}{T} \right)\\
 \end{align*}
 By using Jensen Gap (lemma \ref{LEM::JENGAP}), and the boundedness of $\VLPROB{\THor}{\cdot}$ since the rewards are bounded between $0$ and $1$, we obtain:
 \begin{align*}
     \stackrel{(d)}{\leq} &\lim_{T \to \infty} \epsilon + \THor \Bigg(\Bigg(1 + \frac{k}{\syncconst{k}} \Bigg)\frac{k}{\syncconst{k}}  + \|\laginf\|_{\infty} + 1\Bigg)\sqrt{\frac{\pi}{N}} \\
    &+ \left[1 + 2\|\mu\|_{\infty} + \left(1 + \frac{k}{\syncconst{k}} \right)\sum_{l = 1}^{\THor}\left(1 - \frac{\syncconst{k}}{k} \right)^{l - 1}\right] \frac{\alpha N - \lfloor \alpha N \rfloor}{N} + O\left(\frac{1}{T} \right).
 \end{align*}
 We obtain the result by using that $\sum_{l = 1}^{\THor}\left(1 - \frac{\syncconst{k}}{k} \right)^{l - 1}\le k/\rho_k$. 
The details for steps (a) - (d) can be found in \ref{APP:SEC:FP}.


\begin{remark}
    We remark that our proof is highly modular. Any set of assumptions that ensures the continuity of the bias function $\hstar{\cdot}$ can be slotted in to replace Assumption \ref{AS::SA}, leaving the proof largely unchanged. Similarly, a different style of concentration can be used to replace steps $(b)$ and $(d)$ depending on the nature of heterogeneity and coupling of the arms. 
\end{remark}

\section{Comparison With Other Policies}\label{SEC:NUM}
\subsection{Connection Between the LP-Priority Policy and the LP-update Policy.\\}
In this section, we take a detailed look at another class of policies that is commonly used for restless multi-armed bandit problems, the LP-priority index \citep{GGY23b, verloop2016asymptotically}. 
In this section, we fundamentally draw a  connection between the LP-update policy and the LP-priority policy. 


We return once more to the relaxed linear program \eqref{EQ::WHITF1}-\eqref{EQ::WHITF3}. 

\begin{align*}
    \gstar := \max_{\bYN \in \mcl{Y}} \frac{1}{N}\sum_{n = 1}^{N} & \Rvec{n} \cdot \bYinf{(n)} \\
    \text{such that,}& \nonumber \\
    \bYZinf{(n)} + \bYUinf{(n)} =&  \bYZinf{(n)} \cdot \Pmat{n}{0} + \bYUinf{(n)} \cdot \Pmat{n}{1}  \hspace{0.2 in } \forall n \\
    \sum_{n = 1}^{N} \| \bYUinf{n}\|_1 \leq& N. \alpha \hspace{0.2 in }
\end{align*}

Let $\blamT_{rel}$ be the Lagrange multiplier corresponding to the budget constraint. As shown by \cite{puterman2014markov}, the optimal policy to solve the LP above corresponds to finding a pair $(\gstar, \mu_{rel})$ of a constant and a vector, respectively, that satisfy the following equation:
\begin{equation}\label{EQ::LPPRIO}
    \gstar + \mu_{rel, s_i}^{(i)} = \max_{a \in \{0, 1\}}\{r(s_i, a) - a \blamT_{rel} + \sum_{s' \in \sspaceN{i}}P^{(i)}_{s'|s_i, a}\mu_{rel, s'}^{(i)}\}
\end{equation}

We will denote by $Q_{rel}^{(i)}(s_i, a) := r^{(i)}(s_i, a) - a \blamT_{rel} + \sum_{s' \in \sspaceN{i}}P^{(i)}_{s'|s_i, a}\mu_{rel, s'}^{(i)}$. When the state is $\bsN := \{s_1 \dots s_N\}$, then the LP priority index for arm $n$ in state $s_n$ is defined by: 
\[
\indxlp{n} := Q_{rel}^{(j)}(s_n, 1) - Q_{rel}^{(j)}(s_n, 0).
\]


We now present a striking result that links the LP-priority policy to the LP-update policy directly.

\begin{proposition}\label{PROP:WHITLUPDATE}
    The LP-priority policy corresponds to the solution of the LP-update problem when the computation horizon is set to $1$. 
\end{proposition}

The proof of this result can be found in \ref{APP::PROOFWHITLUPDATE}. We now pause to emphasize the importance of this proposition in the remark below.
\begin{remark}
    LP-priority policies are very similar to the Whittle's index policy, in fact at the fixed point, the priority order induced by both policies are identical. Both these policies are based around Whittle's relaxation and hence, require UGAP as a prerequisite for convergence. However, the LP priority policy does not require indexability. In this context one notes that simply adding a few more steps of computation (as our simulations indicate, a time horizon $\THor$ of $4$ or $5$ is sufficient for this purpose) allows us to replace the hard to verify condition of UGAP with the easy to verify ergodicity condition. In other words, \emph{adding a few more steps of computation allows us to compensate for the shortcomings of Whittle's relaxation.} One of the shortcomings of our policy is the computational complexity associated with computing a few steps of the fluid solution for each time step, the LP-priority policy computes a fixed priority order that can be computed in advance. This makes us hopeful of reducing our computational burden by using the result from Proposition \ref{PROP:WHITLUPDATE}. This is an important matter of future work.
\end{remark}

\subsection{Numerical Experiments}
In this section we detail the numerical experiments we carried out to verify our Theoretical result. We begin by listing our policies, then we describe our experimental set up and finally elaborate on our results.  
\subsubsection*{Set of Policies}
We compare the LP-update policy with the LP-priority policy described in the previous section. Unless specified otherwise, for the purposes of this simulation, we set the LP-update computation horizon parameter, $\THor$ to $4$. Since most common implementations of the LP-priority policy use an equality constraint for the budget instead of the inequality constraint, we will use the same constraint in our own LP-update policy. Further, we will use the water-filling implementation described in the previous section rather than randomized rounding to make the comparison more fair. An interested reader can read \cite[Remark 3.1]{verloop2016asymptotically} or \cite{GGY23b} for an equivalence between the problem with the equality and inequality constraints. 
We also compare the LP-update policy to the arm id reassignment policy from \cite{zhang2025projLyapunov} which assigns a priority order among bandit arms based on the fixed-point policy of each arm. Note, the LP-priority policy can be thought of as fixing the order across the \emph{union of states}, on the other hand, under id reassignment, the priority order is assigned to the \emph{arms} themselves. We were unable to find a code for \cite{zhang2025projLyapunov}, we have therefore written our own implementation. This policy is known to produce an error of $O(1/\sqrt{N})$ when the optimal single-armed policy has an upper bound on the mixing time. 

\subsubsection*{Experiment}
\textbf{Setup}
We simulate a set of $N$ arms over a horizon $T = 1000$ time steps. We then difference scenarios. In some scenarios, arms are \emph{drawn randomly}, and in other cases, we manually construct scenarios in which the LP-priority is not asymptotically optimal. The construction are detailed below. Additional experimental details, including the exact matrices used, are provided in \ref{APP:EXP}.

When we say that our arms are \emph{drawn randomly}, we first select a state space $\{\sspaceN{n}\}$ uniformly at random between $1$ to $10$. Next, we assign to each arm a transition kernel $\Pmat{n}{a}$ of size $\sspaceN{n} \times 2 \times \sspaceN{n}$ whose entries are chosen according to an exponential random variable with parameter $1$ and normalized in order to construct a stochastic matrix. Finally, we construct a reward vector for each arm $\Rvec{n}$ whose entries are chosen according to an exponential random variable with parameter $1$. Each arm is set with a random initial state but in all our comparisons we ensure that all policies begin with the same initial conditions but will evolve over the time horizon depending on the actions chosen by the policies.

In order to illustrate what happens when LP-priority policy is not asymptotically optimal, we produce three examples: two homogeneous and one heterogeneous. Our first counter-example is taken from \cite{HXCW23} with budget parameter $\budget = 0.5$, we call it Counter-example-Hong. Our second is the counter-example from \cite{GGY23} with the budget parameter $\budget = 0.4$, we call this example Counter-example-Yan. The third counter-example is heterogeneous, we set half the arms to Counter-example-Hong and the other half to Counter-example-Yan, and set $\budget = 0.4$.

\textbf{Metric of Comparison: } 
Our metric of comparison is the normalized average reward over time. Concretely we measure at each time $t \geq 1$, the following average reward, normalized by the upper bound $\gstar$ on the average reward:
\[
\frac{ \sum_{t = 1}^{T} \frac{1}{N} \sum_{n = 1}^{N} \Rvec{n}(S_n(t), A_n(t))}{T \times \gstar}.
\]
Note, this metric can exceed $1$ for small values of $T$, however, as $T$ tends to infinity, this number is bounded above by $1$.

\begin{figure}[ht]
    \centering
    \begin{tabular}{@{}c@{}c@{}@{}c}
        \includegraphics[width=0.33\linewidth]{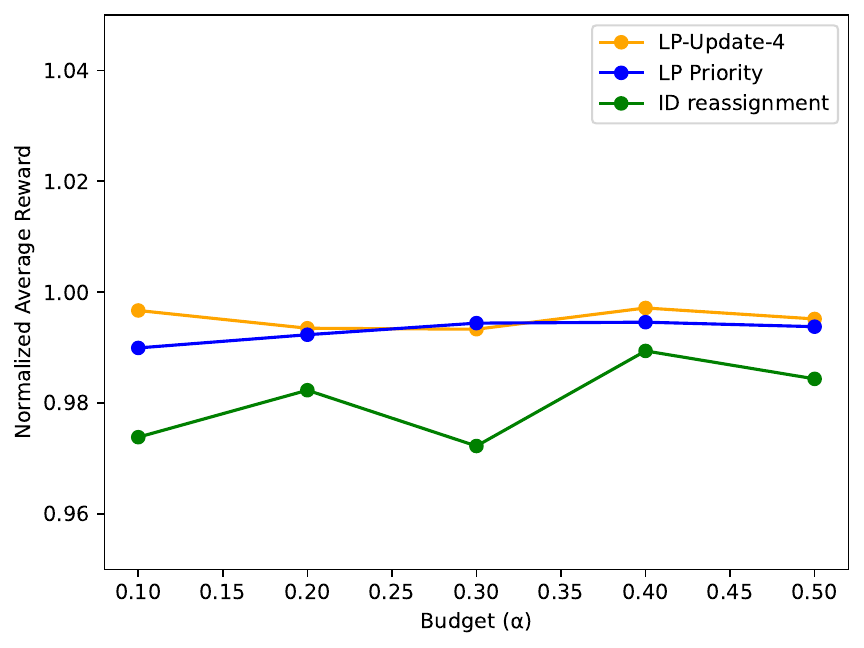}
        &\includegraphics[width=0.33\linewidth]{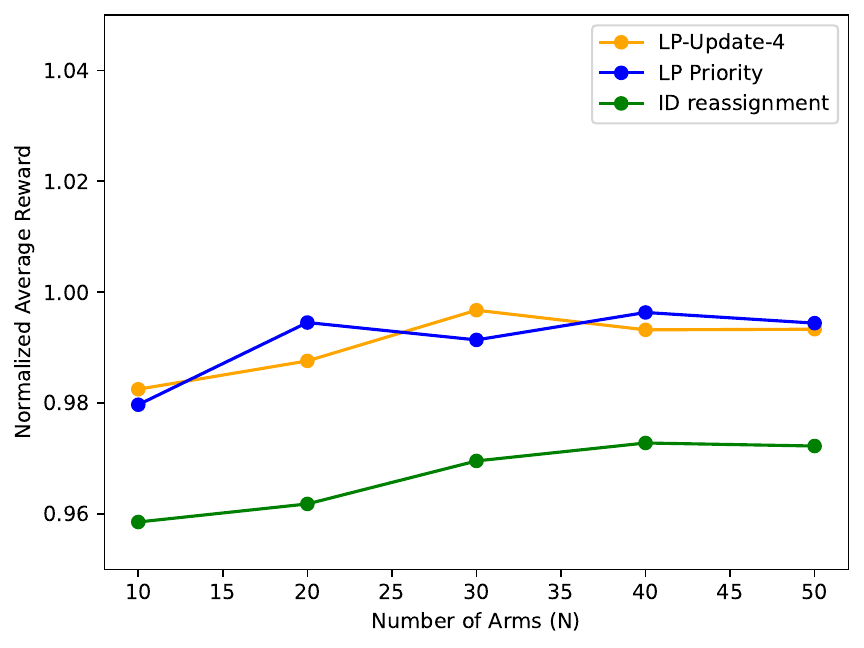}
        &\includegraphics[width=0.33\linewidth]{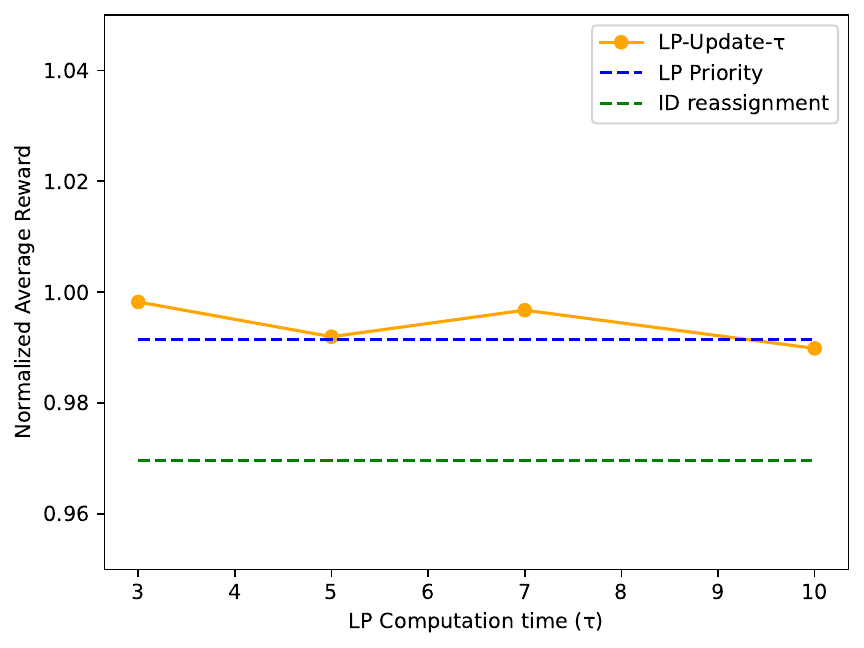}
        \\
        \multicolumn{1}{c}{(a) Influence of $\budget$.} &\multicolumn{1}{c}{(b) Influence of  $N$. }&\multicolumn{1}{c}{(c) Influence of  $\THor$.} 
    \end{tabular}
    \caption{Comparison of the normalized average reward for the three policies on random examples. We study the influence of various parameters ($\budget$, $N$ and $\tau$)}\label{FIG:NARM}
\end{figure}

\textbf{Results: }
In Figure~\ref{FIG:NARM}, Panel~(a) studies the influence of the budget $\budget$ for a collection of $N = 50$ arms. All the policies perform fairly consistently for a randomly drawn configuration over different values of the budget. Next, in Panel~(b) we show the dependence of our reward as a function of the number of arms $N$; varying the values of $N$ between $10$ to $50$ in increments of $10$, setting the budget to $\alpha = 0.3$. In Panel~(c), we study the influence of the computation horizon $\THor$ for a collection of $N = 50$ arms. Surprisingly, our policy is largely invariant to the computation time horizon $\THor$ in these randomly drawn examples. While it is clear that our LP-update policy outperforms the Id-reassignment policy, it seems to perform almost exactly as well as the LP-priority policy. This difference does not disappear even as the number of arms grows, although the overall performance of all policies improves with the number of arms, as can be seen in Figure~\ref{FIG:NARM}(b). We believe that one of the main reasons why the LP-update policy and the LP-priority policies perform similarly for the randomly drawn arms is because (although this is not verifiable) most arm configurations satisfy the UGAP condition and hence, the LP-priority policy is able to steer the system towards the optimal invariant state-action measure just as the LP-update policy does.
 

\begin{figure}[ht]
    \centering
    \begin{tabular}{@{}c@{}c@{}@{}c}
        \includegraphics[width=0.34\linewidth]{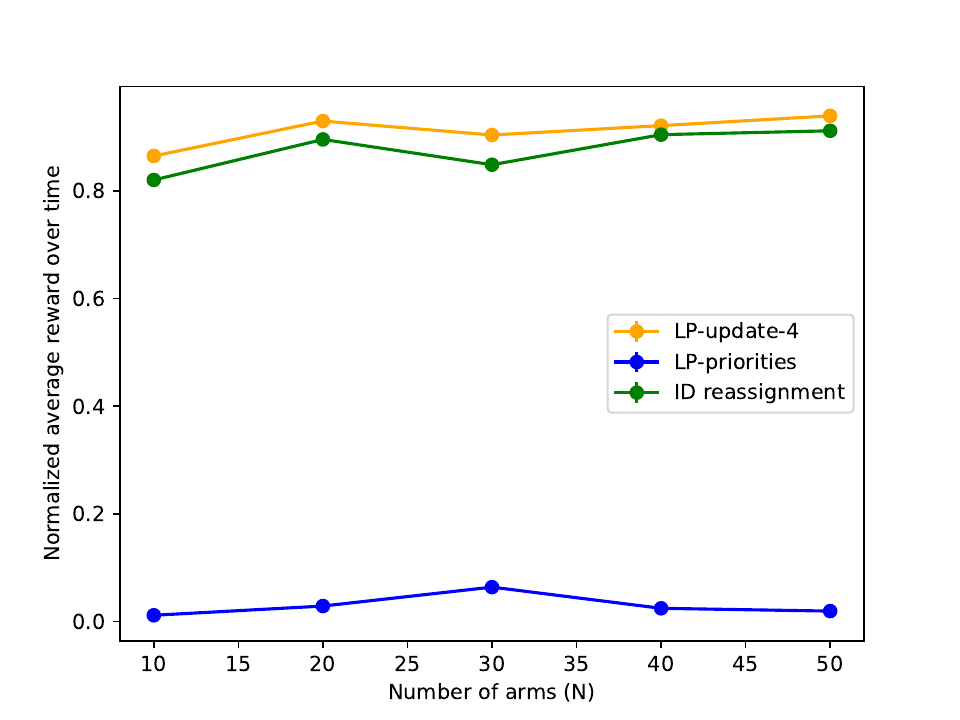}
        &\includegraphics[width=0.34\linewidth]{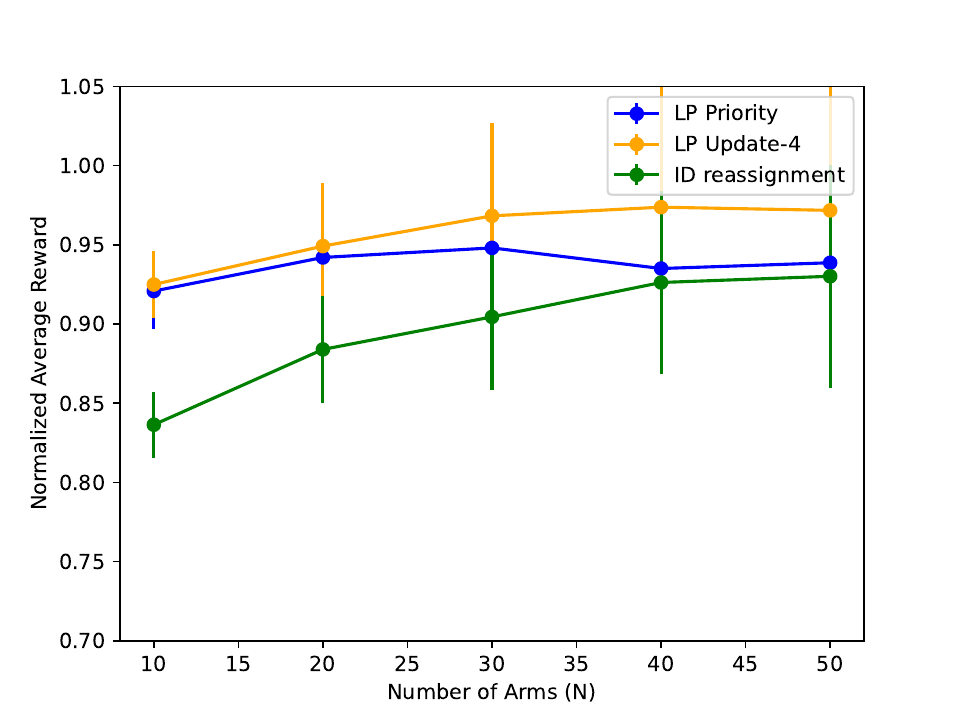}
        &\includegraphics[width=0.31\linewidth]{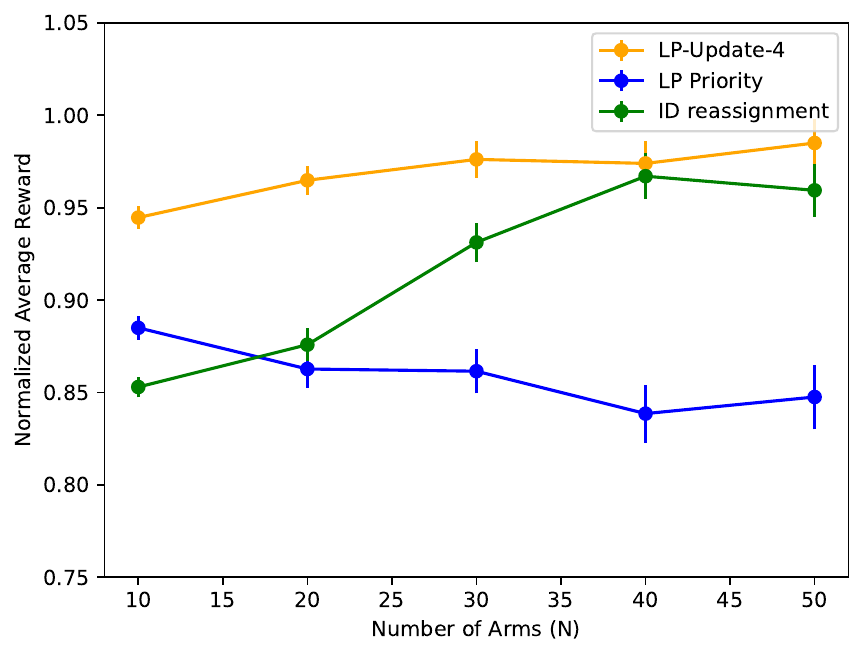}
        \\
        \multicolumn{1}{c}{(a) Counter-example-Hong} &\multicolumn{1}{c}{(b) Counter-example-Yan.} &\multicolumn{1}{c}{(c) Mixed Counter Example. }
    \end{tabular}
    \caption{Comparison of the normalized average reward of the three policies for the three counter examples to the LP-priority Policy.}\label{FIG:OTHER}
\end{figure}

The simulation results in Figure~\ref{FIG:NARM} seem to suggest that the LP-priority policy is sufficient to obtain excellent performance for most systems. What is more, when compared to our own LP-update policy, the LP-priority policy requires only one pre-processing step to decide the priority order. This begs the question, why should we follow either the LP-update policy or the ID-reassignment policy at all? Figure~\ref{FIG:OTHER} tells us another story. Figure~\ref{FIG:OTHER} examines the performance of the three policies as a function of the number of arms for the three counter examples. Under the counter examples the LP-priority policy performs poorly. The LP-update policy performs consistently well. In Panel~(b), ID-reassignment appears to be performing poorly, but the performance improves as the number of arms grows larger indicating that it will eventually beat the LP-priority policy. Figure (c) explicitly shows how both the LP-update-$4$ policy and the ID-reassignment which were constructed to accommodate a violation of the UGAP condition perform well as the number of arms grow larger.

\begin{figure}[ht]
    \centering
    \begin{tabular}{@{}c@{}c@{}@{}c}
        \includegraphics[width=0.34\linewidth]{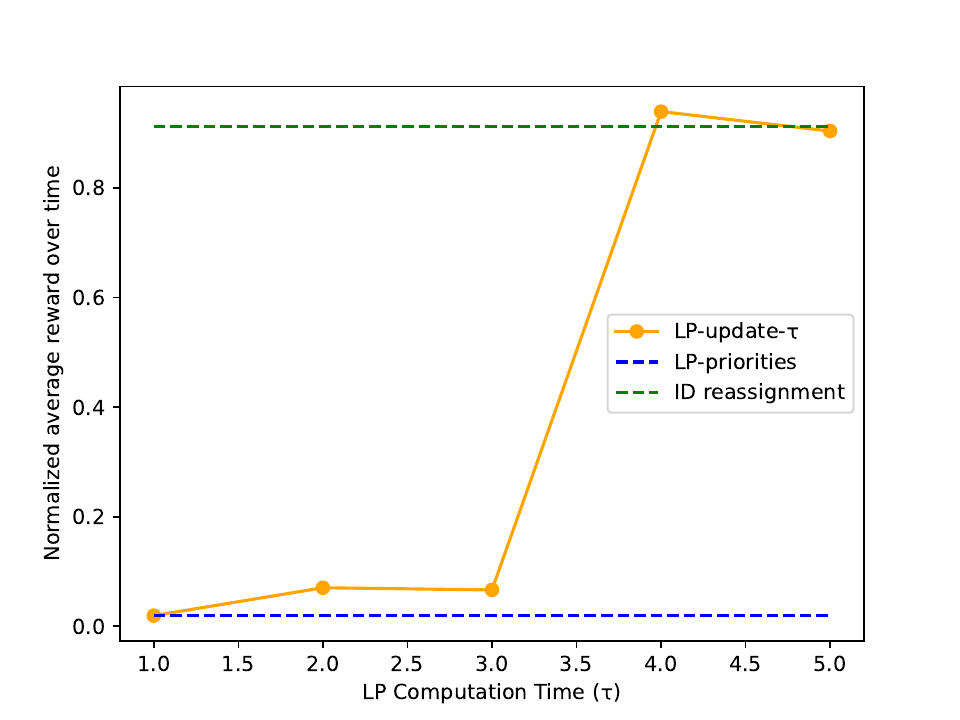}
        &\includegraphics[width=0.32\linewidth]{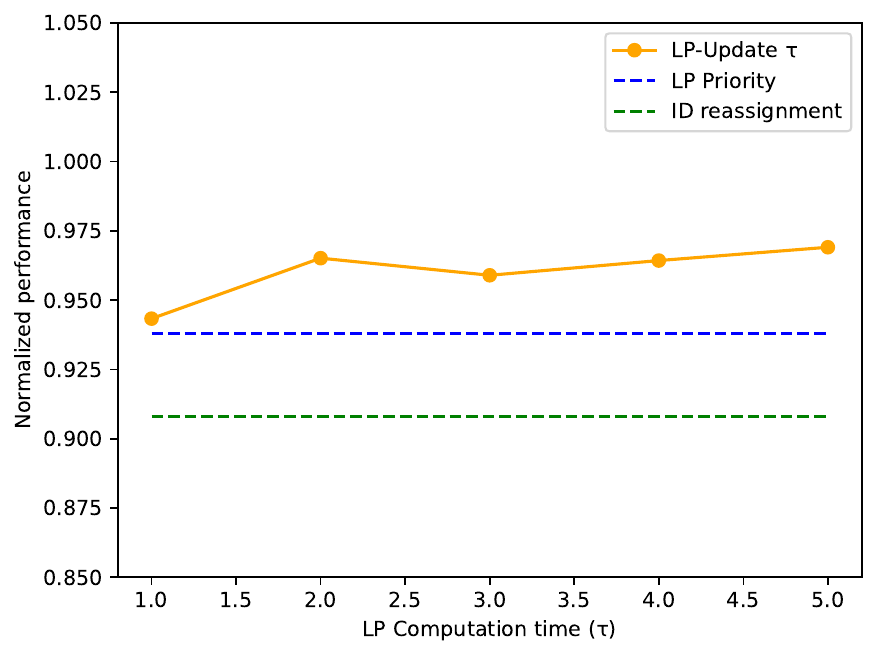}
        &\includegraphics[width=0.32\linewidth]{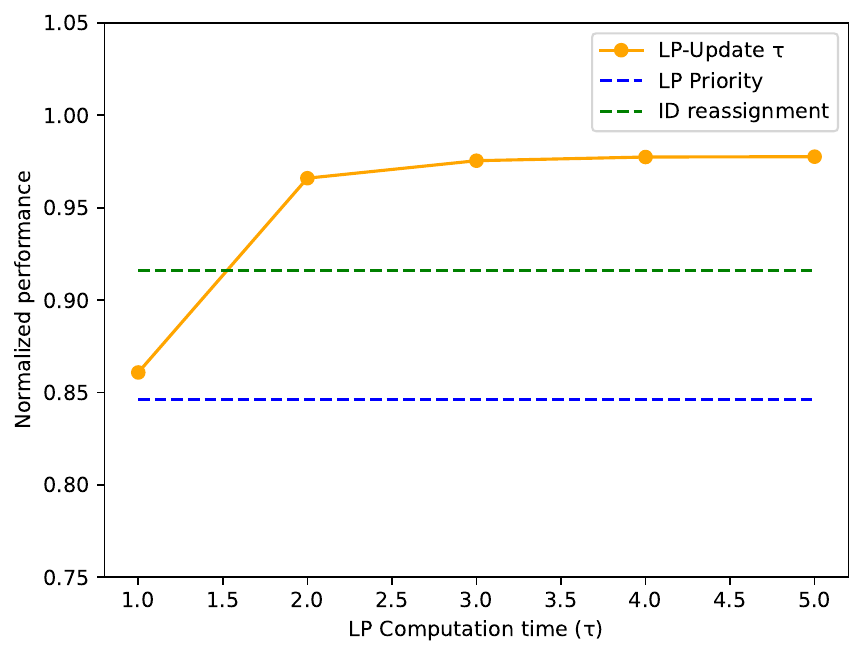}
        \\
        \multicolumn{1}{c}{(a) Counter-Example-Hong.} &\multicolumn{1}{c}{(b) Counter-Example-Yan. }&\multicolumn{1}{c}{(c) Mixed Counter Example.} 
    \end{tabular}
    \caption{Comparison of the average normalized reward as a function of $\THor$, the time horizon used in the LP, for the three counter examples. Recall that the LP-priority and the ID reassignement policies do not depend on $\tau$ (which is why their performance is a straight line as a function of $\tau$).}\label{FIG:TimeComp}
\end{figure}

Our last figure, Figure \ref{FIG:TimeComp} studies the performance for $N = 30$ arms over different values of $\THor$ specifically for the three counter examples. 
The figures indicate a threshold behavior with respect to the time horizon indicating that unless we set $\THor$ to at least $4$ we do not observe good performance from the LP-update algorithm. This indicates the need to more closely examine the role of $\THor$ and \emph{turnpike properties} related to Conjecture \ref{CONJ:1} for the model predictive control algorithm.  These topics are beyond the scope of the current paper but an interesting topic for future work. 

\section{Conclusion}\label{sec:Conclusion}
Through this work we show the asymptotic optimality of the LP-index update policy for the notoriously difficult heterogeneous multi-armed bandit problem under an easily verifiable assumption. What is more, our framework is modular enough so that a different assumption that proves the continuity of our bias function can be slotted in yielding much the same result. Furthermore, we observe the close relationship between our own notion of index which depends on the time horizon and the Whittle's index policy. 
Crucially, using the framework of dissipativity, we are able to show that a finite-horizon computation of the index subverts the need to relax constraints (Whittle's relaxation) while giving us a near-optimal policy for the infinite-horizon average reward problem. Along the way there are several open questions that we may pose given some of the shortcomings of our analysis. The most direct question requires a resolution of our Conjecture \ref{CONJ:1}, which we strongly believe can come through \emph{turnpike properties}. A very related question now follows: can one optimize $\THor$ for the LP-update? Our simulations clearly indicate that for all intents and purposes one can fix a very small time horizon, say $\THor < 5$ in order to solve this problem. 
 Such questions are highly relevant, not only in solving weakly coupled MDPs but also in the more complex learning setting. We conclude by stating that Model Predictive Control methods are an underexplored avenue towards practical solutions to constrained MDP problems as we have demonstrated through this work.
 \bibliographystyle{elsarticle-num} 
  \bibliography{neurips_2025}

@article{zhang2021restless,
  title={Restless bandits with many arms: Beating the central limit theorem},
  author={Zhang, Xiangyu and Frazier, Peter I},
  journal={arXiv:2107.11911},
  year={2021}
}

@article{brown2020index,
  title={Index policies and performance bounds for dynamic selection problems},
  author={Brown, David B and Smith, James E},
  journal={Management Science},
  volume={66},
  number={7},
  pages={3029--3050},
  year={2020},
  publisher={INFORMS}
}

@article{zayas2019asymptotically,
  title={An asymptotically optimal heuristic for general nonstationary finite-horizon restless multi-armed, multi-action bandits},
  author={Zayas-Caban, Gabriel and Jasin, Stefanus and Wang, Guihua},
  journal={Advances in Applied Probability},
  volume={51},
  number={3},
  pages={745--772},
  year={2019}
}

@article{hu2017asymptotically,
  title={An asymptotically optimal index policy for finite-horizon restless bandits},
  author={Hu, Weici and Frazier, Peter},
  journal={arXiv:1707.00205},
  year={2017}
}

@article{ioannidis2016adaptive,
  title={Adaptive caching networks with optimality guarantees},
  author={Ioannidis, Stratis and Yeh, Edmund},
  journal={ACM SIGMETRICS Performance Evaluation Review},
  volume={44},
  number={1},
  pages={113--124},
  year={2016},
  publisher={ACM New York, NY, USA}
}

@article{verloop2016asymptotically,
  title={ASYMPTOTICALLY OPTIMAL PRIORITY POLICIES FOR INDEXABLE AND NONINDEXABLE RESTLESS BANDITS},
  author={Verloop, IM},
  journal={The Annals of Applied Probability},
  pages={1947--1995},
  year={2016},
  publisher={JSTOR}
}

@article{ghosh2022indexability,
  title={Indexability is not enough for whittle: Improved, near-optimal algorithms for restless bandits},
  author={Ghosh, Abheek and Nagaraj, Dheeraj and Jain, Manish and Tambe, Milind},
  journal={arXiv:2211.00112},
  year={2022}
}

@phdthesis{chen:tel-04068056,
  TITLE = {{Close-to-opimal policies for Markovian bandits}},
  AUTHOR = {Yan, Chen},
  NUMBER = {2022GRALM046},
  SCHOOL = {{Universit{\'e} Grenoble Alpes [2020-2023]}},
  YEAR = {2022},
  MONTH = Dec,
  TYPE = {Theses},
  HAL_ID = {tel-04068056},
  HAL_VERSION = {v2},
}

@article{PT99,
 ISSN = {0364765X, 15265471},
 author = {Papadimitriou, Christos H. and Tsitsiklis, John N.},
 journal = {Mathematics of Operations Research},
 number = {2},
 pages = {293--305},
 publisher = {INFORMS},
 title = {The Complexity of Optimal Queuing Network Control},
 urldate = {2024-09-16},
 volume = {24},
 year = {1999}
}

@article{Wh88,
 ISSN = {00219002},
 author = {P. Whittle},
 journal = {Journal of Applied Probability},
 pages = {287--298},
 publisher = {Applied Probability Trust},
 title = {Restless Bandits: Activity Allocation in a Changing World},
 urldate = {2024-09-16},
 volume = {25},
 year = {1988}
}

@inproceedings{
HXCW24,
title={Achieving Exponential Asymptotic Optimality in Average-Reward Restless Bandits without Global Attractor Assumption},
author={Yige Hong and Qiaomin Xie and Yudong Chen and Weina Wang},
booktitle={NeurIPS 2025 Workshop MLxOR: Mathematical Foundations and Operational Integration of Machine Learning for Uncertainty-Aware Decision-Making},
year={2025},
url={https://openreview.net/forum?id=igPYJRM3Yi}
}

@inproceedings{HXCW23,
title={Restless Bandits with Average Reward: Breaking the Uniform Global Attractor Assumption},
author={Hong, Yige and Xie, Qiaomin and Chen, Yudong and Wang, Weina},
booktitle={Thirty-seventh Conference on Neural Information Processing Systems},
year={2023},
}

@Article{GGY23,
  author={Gast, Nicolas and Gaujal, Bruno and Yan, Chen},
  title={{Exponential asymptotic optimality of Whittle index policy}},
  journal={Queueing Systems: Theory and Applications},
  year=2023,
  volume={104},
  number={1},
  pages={107-150},
  month={June},
  doi={10.1007/s11134-023-09875-}
}

@article{GGY23b,
author = {Gast, Nicolas and Gaujal, Bruno and Yan, Chen},
title = {Linear Program-Based Policies for Restless Bandits: Necessary and Sufficient Conditions for (Exponentially Fast) Asymptotic Optimality},
journal = {Mathematics of Operations Research},
year = {2023},
doi = {10.1287/moor.2022.0101}
}

@INPROCEEDINGS{yan2024,
  author={Yan, Chen},
  booktitle={2024 IEEE 63rd Conference on Decision and Control (CDC)}, 
  title={An Optimal-Control Approach to Infinite-Horizon Restless Bandits: Achieving Asymptotic Optimality with Minimal Assumptions}, 
  year={2024},
  volume={},
  number={},
  pages={6665-6672},
  doi={10.1109/CDC56724.2024.10886365}}

@article{WW90,
 ISSN = {00219002},
 author = {Weber, Richard  and Weiss, Gideon},
 journal = {Journal of Applied Probability},
 number = {3},
 pages = {637--648},
 publisher = {Applied Probability Trust},
 title = {On an Index Policy for Restless Bandits},
 urldate = {2024-09-16},
 volume = {27},
 year = {1990}
}

@article{DTGLS14,
author = {Damm, Tobias and Gr\"{u}ne, Lars and Stieler, Marleen and Worthmann, Karl},
title = {An Exponential Turnpike Theorem for Dissipative Discrete Time Optimal Control Problems},
journal = {SIAM Journal on Control and Optimization},
volume = {52},
number = {3},
pages = {1935-1957},
year = {2014},
doi = {10.1137/120888934},
eprint = {    
        https://doi.org/10.1137/120888934
}
}

@article{avrachenkov2022whittleIB,
  title={Whittle index based Q-learning for restless bandits with average reward},
  author={Avrachenkov, Konstantin E and Borkar, Vivek S},
  journal={Automatica},
  volume={139},
  pages={110186},
  year={2022},
  publisher={Elsevier}
}

@book{puterman2014markov,
  title={Markov Decision Processes: Discrete Stochastic Dynamic Programming},
  author={Puterman, M.L.},
  isbn={9781118625873},
  series={Wiley Series in Probability and Statistics},
  year={2014},
  publisher={Wiley}
}

@article{veatch1996scheduling,
  title={Scheduling a make-to-stock queue: Index policies and hedging points},
  author={Veatch, Michael H and Wein, Lawrence M},
  journal={Operations Research},
  volume={44},
  number={4},
  pages={634--647},
  year={1996},
  publisher={INFORMS}
}

@article{nino2002dynamic,
  title={Dynamic allocation indices for restless projects and queueing admission control: a polyhedral approach},
  author={Nino-Mora, Jos{\'e}},
  journal={Mathematical programming},
  volume={93},
  number={3},
  pages={361--413},
  year={2002},
  publisher={Springer}
}

@article{dance2019optimal,
  title={Optimal policies for observing time series and related restless bandit problems},
  author={Dance, Christopher R and Silander, Tomi},
  journal={Journal of Machine Learning Research},
  volume={20},
  number={35},
  pages={1--93},
  year={2019}
}

@article{MearaWebcrawl01,
  title={A topic-specific Web robot model based on restless bandits},
  author={O'Meara, Tadhg and Patel, Ahmed},
  journal={IEEE Internet Computing},
  volume={5},
  number={2},
  pages={27--35},
  year={2001},
  publisher={IEEE}
}

@Article{NinoMora23,
AUTHOR = {Niño-Mora, José},
TITLE = {Markovian Restless Bandits and Index Policies: A Review},
JOURNAL = {Mathematics},
VOLUME = {11},
YEAR = {2023},
NUMBER = {7},
ARTICLE-NUMBER = {1639},
ISSN = {2227-7390},
DOI = {10.3390/math11071639}
}

@book{altman1999constrained,
  title={Constrained Markov Decision Processes},
  author={Altman, E.},
  isbn={9780849303821},
  lccn={99210415},
  series={Stochastic Modeling Series},
  url={https://books.google.fr/books?id=3X9S1NM2iOgC},
  year={1999},
  publisher={Taylor \& Francis}
}

@Article{DBert2016,
  author={Bertsimas, Dimitris and Mišić, Velibor V. },
  title={{Decomposable Markov Decision Processes: A Fluid Optimization Approach}},
  journal={Operations Research},
  year=2016,
  volume={64},
  number={6},
  pages={1537-1555},
  month={December},
  doi={10.1287/opre.2016.1531}
}

@inbook{McDiarmid_1989,
place={Cambridge}, 
series={London Mathematical Society Lecture Note Series},
title={On the method of bounded differences}, booktitle={Surveys in Combinatorics, 1989: Invited Papers at the Twelfth British Combinatorial Conference},
publisher={Cambridge University Press},
author={McDiarmid, Colin},
year={1989}, pages={148–188},
collection={London Mathematical Society Lecture Note Series}}

@InProceedings{gast2024MPCoptimal,
  title = 	 {Model predictive control is almost optimal for restless bandits},
  author =       {Gast, Nicolas and Narasimha, Dheeraj},
  booktitle = 	 {Proceedings of Thirty Eighth Conference on Learning Theory},
  pages = 	 {2326--2361},
  year = 	 {2025},
  editor = 	 {Haghtalab, Nika and Moitra, Ankur},
  volume = 	 {291},
  series = 	 {Proceedings of Machine Learning Research},
  month = 	 {30 Jun--04 Jul},
  publisher =    {PMLR},
  url = 	 {https://proceedings.mlr.press/v291/gast25a.html}
}

@ARTICLE{Borkar2018,
  author={Borkar, Vivek S. and Kasbekar, Gaurav S. and Pattathil, Sarath and Shetty, Priyesh Y.},
  journal={IEEE Transactions on Control of Network Systems}, 
  title={Opportunistic Scheduling as Restless Bandits}, 
  year={2018},
  volume={5},
  number={4},
  pages={1952-1961},
  doi={10.1109/TCNS.2017.2774046}}

@article{zhang2025projLyapunov,
  title={Projection-based Lyapunov method for fully heterogeneous weakly-coupled MDPs},
  author={Zhang, Xiangcheng and Hong, Yige and Wang, Weina},
  journal={Advances in Neural Information Processing Systems},
  volume={38},
  pages={91054--91076},
  year={2026}
}

@article{Brown&Zhang25,
author = {Brown, David B. and Zhang, Jingwei},
title = {Fluid Policies, Reoptimization, and Performance Guarantees in Dynamic Resource Allocation},
year = {2025},
issue_date = {March-April 2025},
publisher = {INFORMS},
address = {Linthicum, MD, USA},
volume = {73},
number = {2},
issn = {0030-364X},
doi = {10.1287/opre.2022.0601},
journal = {Oper. Res.},
month = mar,
pages = {1029–1045},
numpages = {17}
}

@article{Yu2021ModelPC,
  title={Model predictive control for autonomous ground vehicles: a review},
  author={Shuyou Yu and Matthias Hirche and Yanjun Huang and Hong Chen and Frank Allg{\"o}wer},
  journal={Autonomous Intelligent Systems},
  year={2021},
  volume={1},
}

@ARTICLE{2022SulmanMPCmicrogrid,
  author={Shahzad, Sulman and Abbasi, Muhammad Abbas and Chaudhry, Muhammad Akmal and Hussain, Muhammad Majid},
  journal={IEEE Access}, 
  title={Model Predictive Control Strategies in Microgrids: A Concise Revisit}, 
  year={2022},
  volume={10},
  number={},
  pages={122211-122225},
  doi={10.1109/ACCESS.2022.3223298}}

@article{schwenzer_review_2021,
	title = {Review on model predictive control: an engineering perspective},
	volume = {117},
	issn = {1433-3015},
	doi = {10.1007/s00170-021-07682-3},
	number = {5},
	journal = {The International Journal of Advanced Manufacturing Technology},
	author = {Schwenzer, Max and Ay, Muzaffer and Bergs, Thomas and Abel, Dirk},
	month = nov,
	year = {2021},
	pages = {1327--1349},
}

@article{CASTELLETTI2023442,
title = {Model Predictive Control of water resources systems: A review and research agenda},
journal = {Annual Reviews in Control},
volume = {55},
pages = {442-465},
year = {2023},
issn = {1367-5788},
doi = {https://doi.org/10.1016/j.arcontrol.2023.03.013},
author = {Andrea Castelletti and Andrea Ficchì and Andrea Cominola and Pablo Segovia and Matteo Giuliani and Wenyan Wu and Sergio Lucia and Carlos Ocampo-Martinez and Bart {De Schutter} and José María Maestre}
}

@article{Avrachenkov2024,
title={Asymptotically optimal policies for weakly coupled Markov decision processes}, volume={63},
DOI={10.1017/jpr.2026.10079},
number={3},
journal={Journal of Applied Probability}, 
author={Goldsztajn, Diego and Avrachenkov, Konstantin},
year={2026},
pages={1213–1243}}

\appendix


\section{Properties of the Finite Horizon Dynamic Program\\}\label{APP::FLPROP}

In this section we cover key properties of the dynamic program \eqref{EQ::DYN1}. We will begin by proving strong duality to characterize the solution and move on to proving continuity of the solution in terms of the initial conditions. 

Consider the Lagrangian relaxation on the action constraint. For a fixed $\THor > 0$, one may begin by defining for any $0 \leq \blamT :=\{ \blam{0}, \blam{1} \dots \blam{\THor - 1} \}$, a dynamic program as follows:
\begin{align}\label{EQ:LAGRANGE1}
    \LrelT{\bx}{\THor}{\blamT} := \max_{\bYNt{0} \text{satisfies \eqref{EQ::FLUIDLP2}-\eqref{EQ::FLUIDLP3}}} &\Bigg\{ \frac{1}{N} \sum_{n = 1}^N  \Rvec{n} \cdot \bYVn{n}{0}  + \blam{\THor}(\alpha  - \frac{1}{N}\sum_{n = 1}^{N}\|\bYVUn{n}{0}\|_1 ) \nonumber\\
     \hspace{0.5 in} + \hspace{0.2 in}& \LrelT{\bYNt{0} \cdot \PmatN}{\THor - 1}{\blamT} \Bigg\} 
\end{align}
When $\THor = 0$ we set $\LrelT{\bx}{\THor}{\blamT} = \langle \laginf, \bxs{\bs} \rangle$.

By treating the $\blam{t}$ parameter as the cost of choosing the pull action at time $t$, one relaxes the requirement that at most $\alpha$ fraction of the arms can be pulled. We begin with a decomposition result given a fixed collection of multipliers $\blamT$ and relaxing the action constraint. Further, we will fix the terminal condition $\mu$ in what follows.  
\begin{proposition}\label{PROP:DUAL1}
   The lagrangian defined by \eqref{EQ:LAGRANGE1} satisfies strong duality, is piecewise linear and convex in $\blamT$ and can be decomposed in the following sense :
   \begin{equation}\label{EQ:DECOMP1}
       \LrelT{\bx}{\THor}{\blamT} := \alpha \sum_{t = 0}^{\THor - 1} \blam{t} + \frac{1}{N}\sum_{n = 1}^{N} \VLPROBn{\THor}{\bx^{(n)}}{n}{\blamT}
   \end{equation}
   where $\VLPROBn{\THor}{\bx^{(n)}}{n}{\blamT}$ is defined recursively by,
   \begin{align}
       \VLPROBn{\THor}{\bx^{(n)}}{n}{\blamT} := \max_{\bYVn{n}{t} \text{s.t. \eqref{EQ::FLUIDLP3} holds}} &\Bigg\{ \RCvec{n}(0) \cdot \bYVn{n}{0} + \VLPROBn{\THor - 1}{\bYVn{n}{0} \cdot \Pmat{n}{}}{n}{\blamT}\Bigg\}
   \end{align}
   and $\VLPROBn{0}{\bx^{(n)}}{n}{\blamT} = \laginf^{(n)} \cdot \bx^{(n)}$. Where $\RCvec{n}$ indicates the cost adjusted reward vectors i.e. $\RCvec{n, 0}(t) := \Rvec{n, 0}$ and $\RCvec{n, 1}(t) := \Rvec{n, 1} - \blam{t}\Onne$. Next, 
\end{proposition}
\begin{proof}{Proof}
   Weak duality follows since any feasible policy for \eqref{EQ::FLUIDLP1}- \eqref{EQ::FLUIDLP4} is valid for the relaxed problem \eqref{EQ:LAGRANGE1} and when $\blamstarT > 0$, $\sum_{t = 0}^{\THor - 1} \blam{t}\left(\alpha - \frac{1}{N}\sum_{n = 1}^{N} \|\bYVUn{t}{n}\|_1\right) \geq 0$. Hence, one has:
   \begin{equation}\label{EQ::WKDUAL}
       \VLPROB{t}{\bx} \leq \LrelT{\bx}{t}{\blamT}
   \end{equation}
    The decomposition result follows simply by definition. We will call $\VLPROBn{t}{\cdot}{n}{\blamT}$ the item specific dynamic program in keeping with the notation of \cite{brown2020index}. 
    Note that the all-zero action is always feasible, concretely if one sets the action vector $\bYFLU{t} = 0$ for all $t$ is always a feasible state-action vector for our problem. This means that Slater's condition is always satisfied and hence, for any $N$ our solution space is non-empty. Further, due to the compactness of our solution space, a maximizer always exists and hence, the problem is well defined.
    Next, let $\bYVn{n}{[0 : \THor - 1]}$ be a fixed feasible set of state action vectors. Now, given $\bYVn{n}{[0 : \THor - 1]}$ it is not hard to verify that if $\VLPROBn{[0 : \THor - 1]}{\bYVn{n}{[0 : \THor - 1]}}{n}{\blamstarT}$ denotes the value function induced by this policy under fixed $\bYVn{n}{[0 : \THor]}$, the value function is affine (and decreasing) in $\blamT$. Since, the value $\LrelT{\bx}{t}{\blamT}$ is a point-wise maximum of a finite collection of affine functions, $\LrelT{\bx}{t}{\blamT}$ is a convex, piece-wise linear function of $\blamT$ and strong duality is satisfied.  
\end{proof}

A direct consequence of the Strong duality above is the following result which we formalize in the form of a corollary.

\begin{corollary}\label{COR::EQPRDU}
    There exists a $\blamstarT := \{\blamstar{0}, \blamstar{1} \dots \blamstar{\THor - 1}\}$ such that problem \eqref{EQ:LAGRANGE1} and \eqref{EQ::FLUIDLP1}-\eqref{EQ::FLUIDLP4} are equivalent. Further, $\blamstarT$ satisfies the following for each time $t \in \{0, 1, \dots \THor - 1\}$,
    \begin{equation}\label{EQ:CORFLPOL1}
        \blamstar{t}\left(\alpha - \frac{1}{N}\sum_{n = 1}^{N}\|\bYVn{n}{t}\|_1 \right) = 0.
    \end{equation}
\end{corollary}

\begin{definition}\label{DEF:LAGMUL}
    When we wish to distinguish the $\THor$ sequence of multipliers that depend on a state we will denote it by $\blamstar{\bx} := \{\blamstar{0}, \blamstar{1} \dots \blamstar{\THor - 1}\}$ corresponding to the multipliers from Corollary \ref{COR::EQPRDU} with initial condition $\bx$. Hence, for fixed $\THor$, given a state $\bx(t)$ at time $t$, we can define the Lagrange multipliers, $\blamstarT(\bx(t)) := \{\blamstar{t}, \blamstar{t + 1} \dots \blamstar{t + \THor - 1}\}$ that depend on the initial state $\bx(t)$.  
\end{definition}

It follows that when $\blamstar{t} > 0$, $\frac{1}{N}\sum_{n = 1}^{N}\|\bYFLUn{t}{n}\|_1 = \alpha$ and when $\blamstar{t} = 0$ we have $\frac{1}{N}\sum_{n = 1}^{N}\|\bYFLUn{t}{n}\|_1 \leq \alpha$.

\begin{lemma}\label{APP::LEMCONCAVE}
    The value of the dynamic program, 
    $\VLPROB{t}{\bx}$ is concave in the initial state $\bx$.
\end{lemma}
\begin{proof}{Proof:}
    Let $\VLPROB{\THor}{\bx_1}$ and $\VLPROB{\THor}{\bx_2}$ be two values associated with the initial states $\bx_1$ and $\bx_2$ respectively and solutions $\bYVN{[0 : \THor]}{1}$ and  $\bYVN{[0 : \THor]}{2}$ respectively. 
    Given any $\eta \in [0, 1]$, consider a state $\bx_{\eta} := \eta \bx_1 + (1 - \eta) \bx_2$, it is not hard to check that $\bYVN{[0 : \THor]}{\eta} := \eta \bYVN{[0 : \THor]}{1} + (1 - \eta)\bYVN{[0 : \THor]}{2}$ is a feasible solution to the dynamic program. It now follows that
    \begin{align*}
        &\VLPROB{\THor}{\bx_{\eta}} \geq \VLPROB{\THor, \bYVN{[0 : \THor]}{\eta}}{\bx_{\eta}}\\
        =& \eta \VLPROB{\THor,\bYVN{[0 : \THor]}{1} }{\bx_1} + (1 - \eta) \VLPROB{\THor,\bYVN{[0 : \THor]}{2} }{\bx_2}\\
        =& \eta \VLPROB{\THor}{\bx_1} + (1 - \eta) \VLPROB{\THor}{\bx_2}
    \end{align*}
    which completes the claim. 
\end{proof}

The following lemma (Lemma 14) from \cite{gast2024MPCoptimal} gives us the bound for maximum value of $\blam{t}$. We will state it without proof. 
\begin{lemma}\label{LEM::lambound}
   The maximum value of the Lagrange multiplier is bounded above by $\frac{k}{\syncconst{k}}$. 
\end{lemma}

\subsection*{Proof of Lemma \ref{LEM::VALT}}
We are now ready to prove a key result, Lemma \ref{LEM::VALT}, which states that the finite horizon linear program is Lipschitz in the state. We draw the reader's attention to the proof technique used in the lemma as it will be referred to in our proof of continuity in the sections to come.
\lemValC*

Before we begin our proof, we will add some necessary definitions. 
Let the component $i$ be fixed for the purposes of this proof. Consider a state $\bx^{(i)} \in \Delta(\sspaceN{i})$ that is rational (all the components of $\bx^{(i)}$ are rational numbers). For sufficiently large $M$, one may create $M$ identical copies of the arm with joint state $\bs \in \sspaceN{i}^{M}$ (that are unique up to permutation) such that the empirical measure of the $M$ length vector is equal to $\bx^{(i)}$. Concretely, the probability that an arm is in state $\state$ can be represented by: 
\begin{equation*}
    \bhx^{(i)}_{\state}(\bs) := \frac{1}{M}\sum_{m = 1}^{M}\Onne_{s_{m} = \state } = \bx^{(i)}_{\state}. 
\end{equation*}
Since $i$ is fixed, we will assume that the $\sspaceN{i}$ is also fixed. As an example, if $\sspaceN{i} = \{0, 1\}$ and $\bx^{(i)} = [0.2, 0.8]$, then we can have $M = 5$ identical arms with joint state $[0, 1, 1, 1, 1]$ so that $\bhx^{(i)}_{0}(\bs) = \bx^{(i)}_{0} = 0.2$.

We let $\Delta^{M} \subset \Delta (\sspaceN{i})$ denote the corresponding space of all such empirical measures. Note that each $\bx^{(i)} \in \Delta^{M}$ has a $M$ length vector that is unique up to a permutation of its  elements. For the sake of clarity we will set the synchronization constant i.e, $\syncconst{k} = \syncconst{1} = \rho,k = 1$. Extending the proof to a more general $k$ is straight-forward, but only increases the already cluttered notation. 

\begin{proof}
    Our proof begins with a result for distributions that are within $2/M$ and extends the result to a more general distribution using continuity.
    
    Let $\bhx^{(i)} := \bhx^{(i)}(\bs),\underline{\bhx}^{(i)} := \bhx^{(i)}(\Tilde{\bs}) \in \Delta^{M}$ be two distributions such that $\|\bx^{(i)} - \underline{\bhx}^{(i)}\|_1 \leq 2/M$. This implies that there are two states $\state_1, \state_2$ such that, $\bhx^{(i)}_ {\state_1} = \underline{\bhx}^{(i)}_{\state_1} + 1/M$ and $\bhx^{(i)}_ {\state_2} = \underline{\bhx}^{(i)}_{\state_2} - 1/M$ and we permute the $M$-component vector so that for any $m \geq 2$, $s_m = s_m'$. We consider the difference in value functions 
    when the $i^{\text{th}}$ component, $\bhx^{(i)}$ is replaced by $\underline{\bhx}^{(i)}$.  
    \begin{align}\label{EQ::APCONT1}
        \VLPROB{t}{\bx} - \VLPROB{t}{\Tilde{\bx}} &=\LrelT{\bx}{t}{\blamstarT(\bx)} - \LrelT{\Tilde{\bx}}{t}{\blamstarT(\Tilde{\bx})}
        \leq \LrelT{\bx}{t}{\blamstarT(\Tilde{\bx})} - \LrelT{\Tilde{\bx}}{t}{\blamstarT(\Tilde{\bx})} \nonumber\\
        &= \frac{1}{N}\sum_{n \neq i}\left(\VLPROBn{t}{\bx^{(n)}}{n}{\blamstarT(\Tilde{\bx})} - \VLPROBn{t}{\bx^{(n)}}{n}{\blamstarT(\Tilde{\bx})} \right) \nonumber\\
        &\hspace{0.3 in} + \frac{1}{N} \left( \VLPROBn{t}{\bhx^{(i)}}{i}{\blamstarT(\Tilde{\bx})} - \VLPROBn{t}{\underline{\bhx}^{(i)}}{i}{\blamstarT(\Tilde{\bx})}\right) \nonumber\\
        &= \frac{1}{N} \left( \VLPROBn{t}{\bhx^{(i)}}{i}{\blamstarT(\Tilde{\bx})} - \VLPROBn{t}{\underline{\bhx}^{(i)}}{i}{\blamstarT(\Tilde{\bx})}\right)
    \end{align}
    The first equality is derived from Corollary \ref{COR::EQPRDU} and the second follows from equation \eqref{EQ::WKDUAL}. Fix, $t$, $\laginf$ and $\blamstarT(\Tilde{\bx})$ and define a Q-function at state $\bs$ with an $M$ component action vector $\ba = \{a_1 \dots a_M\}$ as follows:
    \begin{equation}
        Q_i(t, (\bs, \ba)) := \frac{1}{M}\sum_{m = 1}^{M} \Tilde{r}(s_m, a_m) + \sum_{\Tilde{\bs} \in \sspaceN{i}^{M}} \VLPROBn{t - 1}{\underline{\bhx}^{(i)}}{i}{\blamstarT(\Tilde{\bx})}\Pi_{m = 1}^{M} P^{(i)}_{\Tilde{s}_m|s_m, a_m}
    \end{equation}
    Our proof proceeds by induction. The statement is clearly true for $t = 0$, assume the inductive hypothesis for some $t > 0$. We will show that the hypothesis holds for $t + 1$. Next, note that there exists $\ba^{*}$ such that $Q_i({t + 1},(\bs, \ba^{*})) := \VLPROBn{t + 1}{\bhx^{(i)}}{i}{\blamstarT(\Tilde{\bx})}$. We allow $\Tilde{\ba}$ to be equal to $\ba^{*}$ for all $m > 1$ and $0$ for the first action. It now follows that,
    \begin{align*}
        &\VLPROBn{t + 1}{\bhx^{(i)}}{i}{\blamstarT(\Tilde{\bx})} - \VLPROBn{t + 1}{\underline{\bhx}^{(i)}}{i}{\blamstarT(\Tilde{\bx})} \leq Q_i(t + 1,(\bs, \ba^{*})) - Q_i({t + 1},(\bs, \Tilde{\ba}))\\
        &=\frac{1}{M} \left( \Tilde{r}(s_1, a_1) -  r(s'_1, 0)\right) + \sum_{\Tilde{\bs} \in \sspaceN{i}^{M}} \Pi_{m = 2}^{M} P^{(i)}_{\Tilde{s}_m|s_m, a_m} \VLPROBn{t}{\underline{\bhx}^{(i)}}{i}{\blamstarT(\Tilde{\bx})}\left(P^{(i)}_{\Tilde{s}_1|s_1,a_1} -  P^{(i)}_{\Tilde{s}_1 |s'_1, 0}\right)\\
        &= \frac{1}{M} \left( \Tilde{r}(s_1, a_1) -  r(s'_1, 0)\right)\\
        &\hspace{0.3 in} + \sum_{\{\Tilde{s}_2 \dots \Tilde{s}_M\} \in \sspaceN{i}^{M - 1}} \Pi_{m = 2}^{M} P^{(i)}_{\Tilde{s}_m|s_m, a_m}\sum_{k \in \sspaceN{i}}\VLPROBn{t}{\underline{\bhx}^{(i)}}{i}{\blamstarT(\Tilde{\bx})}\left(P^{(i)}_{\Tilde{s}_1|s_1,a_1} - P^{(i)}_{\Tilde{s}_1 |s'_1, 0}\right)
    \end{align*}
    By rewriting $\Tilde{\bs}$ in terms of its $M$ component vector $\Tilde{\bs} = \{l, \Tilde{s}_2 \dots \Tilde{s}_M\}$ and using the product form of the transition kernel under fixed joint actions we may now derive the following relations: 
    \begin{align*}
       &\sum_{k} \VLPROBn{t}{\underline{\bhx}^{(i)}}{i}{\blamstarT(\Tilde{\bx})}\left(P^{i}_{l| s_1, a_1} -  P^{i}_{l| s_1', 0}\right) \\
       &= \sum_{k} \VLPROBn{t}{\underline{\bhx}^{(i)}}{i}{\blamstarT(\Tilde{\bx})}\left(P^{i}_{l|s_1, a_1} - \rho(k) \right) - \sum_{k} \VLPROBn{t}{\underline{\bhx}^{(i)}}{i}{\blamstarT(\Tilde{\bx})}\left(P^{i}_{l|s'_1, 0} - \rho(k) \right)
    \end{align*}
    Where $\rho_i(l) := \min \{P^{i}_{l|s_1, a_1}, P^{i}_{l|s'_1, 0}\}$. It follows that,
    \begin{align*}
      &\sum_{k} \VLPROBn{t}{\underline{\bhx}^{(i)}}{i}{\blamstarT(\Tilde{\bx})}\left(P^{i}_{l|s_1, a_1} - \rho_i(k) \right) - \sum_{k} \VLPROBn{t}{\underline{\bhx}^{(i)}}{i}{\blamstarT(\Tilde{\bx})}\left(P^{i}_{l|s'_1, 0} - \rho_i(k) \right)\\
      &\leq \max_{k} \VLPROBn{t}{\bhx^{(i)}(\bs_{k})}{i}{\blamstarT(\Tilde{\bx})}(1 - \sum_k \rho_i(k)) - \min_{k'} \VLPROBn{t}{\bhx^{(i)}(\bs_{k'})}{i}{\blamstarT(\Tilde{\bx})}(1 - \sum_k \rho_i(k)), 
    \end{align*}
    where $\bs_{k} := \{k, \Tilde{s}_2 \dots \Tilde{s}_M\}$. 
    By setting $\rho = \sum_l \rho_i(l)$ and applying the induction hypothesis we have:
    \begin{align*}
        &\sum_{l} \VLPROBn{t}{\underline{\bhx}^{(i)}}{i}{\blamstarT(\Tilde{\bx})}\left(P^{i}_{l|s_1, a_1} - \rho_i(l) \right) - \sum_{l} \VLPROBn{t}{\underline{\bhx}^{(i)}}{i}{\blamstarT(\Tilde{\bx})}\left(P^{i}_{l|s'_1, 0} - \rho_i(l) \right)\\
        &\leq \frac{1}{M} \left((1 - \rho)\left[\|\laginf\|_{\infty} + \left(1 + \frac{1}{\syncconst{}} \right)\sum_{l = 1}^{t}\left(1 - \rho \right)^{l - 1}\right] \right)
    \end{align*}
    Combining this inequality with Lemma \ref{LEM::lambound}, it follows that:
    \begin{align*}
        &\VLPROBn{t + 1}{\bx^{(i)}}{i}{\blamstarT(\Tilde{\bx})} - \VLPROBn{t + 1}{\Tilde{\bx}^{(i)}}{i}{\blamstarT(\Tilde{\bx})}\\
        &\leq\frac{1}{M} \left[\left( \Tilde{r}(s_1, a_1) -  r(s'_1, 0)\right) + \left((1 - \syncconst{})\left[\|\laginf\|_{\infty} + \left(1 + \frac{1}{\syncconst{}} \right)\sum_{l = 1}^{t}\left(1 - \syncconst{} \right)^{l - 1}\right] \right)\right]\\
        &\leq \frac{1}{M} \left[1 + \blam{t + 1} + \left((1 - \rho)\left[\|\laginf\|_{\infty} + \left(1 + \frac{1}{\syncconst{}} \right)\sum_{l = 1}^{t}\left(1 - \syncconst{} \right)^{l - 1}\right] \right)\right]\\
        &\leq \frac{1}{M}\left[\|\laginf\|_{\infty} + \left(1 + \frac{1}{\syncconst{}} \right)\sum_{l = 1}^{t + 1}\left(1 - \syncconst{} \right)^{l - 1} \right].
    \end{align*}
    The first inequality follows from the induction hypothesis and the second follows from Lemma \ref{LEM::lambound}. By noting that any two points on $\sspaceN{i}^{M}$ can be represented as a sequence of points which differ by one-coordinate we can now show that $$\VLPROBn{t + 1}{\bx^{(i)}}{i}{\blamstarT(\Tilde{\bx})} - \VLPROBn{t + 1}{\Tilde{\bx}^{(i)}}{i}{\blamstarT(\Tilde{\bx})} \leq \left[\|\laginf\|_{\infty} + \left(1 + \frac{1}{\syncconst{}} \right)\sum_{l = 1}^{t + 1}\left(1 - \rho \right)^{l - 1}\right]\|\bx^{(i)} - \Tilde{\bx}^{(i)}\|_1.$$ Our result follows from taking the limit as $M$ tends to infinity. Plugging this result into \eqref{EQ::APCONT1} completes our proof.
\end{proof}

\section{Continuity of the Bias Function}\label{APP:CONTWFL}


We begin this section with a definition of a discounted variant of the reward maximization problem. For any $0 < \beta < 1$, let $\bYN(t)$ be, as before, the joint state action vector. The reward and transition kernel for our discounted problem will remain the same. Let $\bx(0)$ denote the initial state of the joint system, we can state a discounted infinite horizon LP as:  
\begin{align}
    V_{\beta}(\bx(0)) =& \max_{\bYN(t) = [\bYVZN{t}, \bYVUN{t}]}\sum_{t = 0}^{\infty} \beta^{t} \langle \RvecN, \bYN(t) \rangle \label{EQ:APB20}\\
    \text{such that}& \nonumber \\
    \bYVZn{n}{0} + \bYVUn{n}{0} =& \bxn{0}{(n)} \hspace{0.2 in } \forall n  \label{EQ:APB21} \\
    \bYVZn{n}{t + 1} + \bYVUn{n}{t + 1} =& \bYVZn{n}{t} \cdot \Pmat{n}{0} + \bYVUn{n}{t} \cdot \Pmat{n}{1} \hspace{0.2 in } \forall n, t \label{EQ:APB22}\\
    \sum_{n = 1}^{N}\|\bYVUn{n}{t}\|_1 & \leq \budget N \hspace{0.2 in }\forall  t.\label{EQ:APB23}
\end{align}

We denote the set of valid action vector that satisfy \eqref{EQ:APB21} - \eqref{EQ:APB23} by $\mcl{Y}(\bx)$. The existence of an optimal policy follows using standard arguments in optimal control. The following lemma follows directly from [Lemma 4, 5 from \cite{gast2024MPCoptimal}], these arguments mimic the proof of Lemma \ref{LEM::VALT} but in the discounted problem and hence, we will state them without proof :

\begin{lemma}\label{LEM::B21}
    Fix a positive integer $N > 0$ and discount factor $0 \leq \beta < 1$. Under Assumption \ref{AS::SA}, 
    \begin{equation}\label{EQ::CONTDISC}
        |V_{\beta}(\bx) - V_{\beta}(\Tilde{\bx})| \leq \frac{1}{N}\sum_{n = 1}^{N}\|\bx^{(n)} - \Tilde{\bx}^{(n)}\|_1
    \end{equation}
    Further, there exists a constant $g_{\beta}$ and bias function $h_{\beta}(\cdot)$ from the state space to the real numbers such that:
    \begin{equation}
        V_{\beta}(\bx(0)) = (1 - \beta)^{-1} g_{\beta} + h_{\beta}(\bx(0)) 
    \end{equation}
    where $g_{\beta}$ and $h_{\beta}(\cdot)$ satisfy the following fixed point equation:
    \begin{equation}
        g_{\beta} + h_{\beta}(\bx) = \max_{\bYN \in \mcl{Y}(\bx)}\langle \RvecN, \bYN \rangle + \beta h_{\beta}(\bYN \cdot \PmatN)
    \end{equation}
\end{lemma}

A direct consequence of Lemma \ref{LEM::VALT} and Lemma \ref{LEM::B21} above is Proposition \ref{PROP::BIAS_CONT} which guarantees the continuity of the bias function \eqref{EQ::BIASFN}.


\Contbias*

\begin{proof}
    By using standard arguments from \cite{puterman2014markov} one can choose an appropriate sequence $\beta_j$ converging to $1$ as $j$ tends to infinity so that $(1 - \beta_i) V_{\beta_i}(\bx) = \lim_{T \uparrow \infty}\frac{1}{T} \VLPROB{T}{\bx}$. Therefore, by Lemma \ref{LEM::B21} we have,
    \begin{align*}
        |(1 - \beta) V_{\beta}(\bx) - (1 - \beta) V_{\beta}(\Tilde{\bx})| \leq (1 - \beta)\frac{1}{N}\sum_{n = 1}^{N}\|\bx^{(n)} - \Tilde{\bx}^{(n)}\|_1
    \end{align*}
    It now follows that for the sequence $\beta_j$ one has:
    \begin{equation*}
        (1 - \beta_j)\|V_{\beta_j}(\bx) - V_{\beta_j}(\Tilde{\bx})\| \to 0.
    \end{equation*}
    This in turn implies that there exists a constant equal to $\lim_{T \uparrow \infty} \frac{1}{T} \VLPROB{t}{\bx}$, further note that by choosing $\bYN^{\star}$ as the fixed point to the problem \eqref{EQ::WHITF1} - \eqref{EQ::WHITF3} with initial condition $\bYUstarN + \bYZstarN$ we have: 
    \begin{equation*}
        \gstar \leq \lim_{T \uparrow \infty} \frac{1}{T} \VLPROB{T}{\bYUstarN + \bYZstarN} \leq \gstar.
    \end{equation*}
    The lower bound follows since we have fixed our policy at $\bYN^{\star}$ which means we will never leave the initial state $\bx^{\star}$ and the upper bound follows since $\gstar$ is the optimal value corresponding to the relaxed problem.
    It now follows that the limit, $g_{\beta_j} \to \gstar$.
    
    Next, note that
    \begin{align*}
        |\hstar{\bx} - \hstar{\Tilde{\bx}}| =& \lim_{T \to \infty} |T \gstar - \VLPROB{T}{\bx} - T \gstar + \VLPROB{T}{\Tilde{\bx}}|\\
        =& \lim_{T \to \infty} | \VLPROB{T}{\Tilde{\bx}} - \VLPROB{T}{\bx}|
    \end{align*}
    Finally, by Lemma \ref{LEM::VALT} one has,
    \begin{align}
        |\hstar{\bx} - \hstar{\Tilde{\bx}}| \leq& \lim_{T \to \infty} \left[\|\laginf\|_{\infty} + \left(1 + \frac{k}{\syncconst{k}} \right)\sum_{l = 1}^{T}\left(1 - \frac{\syncconst{k}}{k} \right)^{l - 1}\right]\frac{1}{N}\sum_{n = 1}^{N}\|\bx^{(n)} - \Tilde{\bx}^{(n)}\|_1 \nonumber\\
        =& \left[\|\laginf\|_{\infty} + \left(1 + \frac{k}{\syncconst{k}} \right)\frac{k}{\syncconst{k}}\right]\frac{1}{N}\sum_{n = 1}^{N}\|\bx^{(n)} - \Tilde{\bx}^{(n)}\|_1
    \end{align}
    We have thus shown that the bias function is Lipschitz and hence bounded on the simplex of distributions. Convexity of $\hstar{\cdot}$ is a direct consequence of the concavity of the dynamic program.
    The fixed point equation now follows by taking the limit. Hence, we can conclude $\eqref{EQ::HCONT}$. This completes the proof.
\end{proof}

\section{Proof of Dissipativity}\label{APP:DISS}

\subsection{Proof of Lemma \ref{LEMM::BIAS}}
\begin{proof}
Recall the fixed point problem, \eqref{EQ::WHITF1} - \eqref{EQ::WHITF3},
\begin{align*}
   \gstar =& \max_{\bYN} \frac{1}{N}\sum_{n = 1}^{N} \Rvec{n} \cdot \bYinf{n} \\
    &\text{such that,} \\
    &\bYZinf{n} + \bYUinf{n} = \bYUinf{n} \cdot \Pmat{n}{1} + \bYZinf{n} \cdot \Pmat{n}{0} \hspace{0.1 in} \forall \hspace{0.1 in} n \\
    &\sum_{n = 1}^{N} \| \bYUinf{n}\|_1 \leq N \alpha. 
\end{align*}

Recall by Definition \ref{DEF:GSTAR} we denote the solution to the problem above by $\bYstarN$ and the value by $\gstar$. Further, there exist Lagrange multiplier $\laginf$ such that,
\begin{align*}
    \gstar =& \frac{1}{N} \sum_{n = 1}^{N}\Rvec{n} \cdot \bYstar{n} - \langle \laginf, (\bYZstarN + \bYUstarN) - \bYstarN\cdot\PmatN   \rangle =  \frac{1}{N} \sum_{n = 1}^{N}\Rvec{n} \cdot \bYstar{n}\\
    \geq& \frac{1}{N} \sum_{n = 1}^{N}\Rvec{n} \cdot \bYinf{n} - \langle \laginf, (\bYZN + \bYUN) - \bYN \cdot \PmatN \rangle \hspace{0.2 in} \text{for any } \bYN
\end{align*}
The proof now follows by noting the equivalence between $(\bx, \bu)$, $\bYN$ hence,
\begin{equation}
    0=\gstar -  \frac{1}{N} \sum_{n = 1}^{N}\Rvec{n} \cdot \bYstar{n} \leq \gstar - \frac{1}{N} \sum_{n = 1}^{N}\Rvec{n} \cdot \bYinf{n} + \langle \laginf, (\bYZN + \bYUN) - \bYN \cdot \PmatN \rangle.
\end{equation}
By setting the storage cost to $\psi(\bx) = \langle \laginf, \bx \rangle $ we complete the proof.
\end{proof}

\subsection{Proof of Lemma \ref{LEM::INFCOST}}
\begin{proof}
Recall, our surrogate loss function requires a minimization of the accumulated non-negative rotated cost function \eqref{EQN::COSTMINDEF},
\begin{equation*}
    \CostT{\bx}{\THor} := \min_{\bYN s.t. \eqref{EQ::FLUIDLP2} -\eqref{EQ::FLUIDLP4}}\sum_{t = 0}^{\THor - 1}\rotcost{c}{\bx(t)}
\end{equation*}
Now consider,
    \begin{align*}
        \CostT{\bx}{\THor} =& \min_{\bYN s.t. \eqref{EQ::FLUIDLP2} -\eqref{EQ::FLUIDLP4}}\sum_{t = 0}^{\THor - 1}\rotcost{c}{\bx(t)}\\
        =& \tau \gstar  + \langle \laginf, \bx \rangle - \max_{\bYNt{t}: t \in \{0, 1 \dots \THor\}}\left[ \sum_{t = 0}^{\THor - 1}\left( \frac{1}{N} \sum_{n = 1}^N \Rvec{n} \cdot\bYVn{t}{n} \right) + \langle \laginf, (\bYVZN{\THor} + \bYVUN{\THor}) \rangle \right]
    \end{align*}
    The next result follows by the telescoping sum of the storage costs. Substituting the definition of $\VLPROB{\THor}{\bx}$ from \eqref{EQ::FLUIDLP1}-\eqref{EQ::FLUIDLP4} gives us:
    \begin{align*}
        \CostT{\bx}{\THor} = \langle \laginf, \bx \rangle + \tau \gstar - \VLPROB{\THor}{\bx}
    \end{align*}
    Now taking the limit as $\THor$ tends to infinity and using equation \ref{EQ::BIASFN} gives us:
    \begin{equation}
        \CostT{\bx}{\infty} = \langle \laginf, \bx \rangle + \hstar{\bx}.
    \end{equation}
    Finally, note that by Proposition (\ref{PROP::BIAS_CONT}), $\hstar{\cdot}$ is Lipschitz continuous, hence bounded by $$\left[\|\laginf\|_{\infty} +  \left(1 + \frac{k}{\syncconst{k}} \right) \frac{k}{\syncconst{k}}\right] =: \|\hstar{\cdot}\|_{\infty}$$ on the simplex of distributions. It follows that $\CostT{\cdot}{\infty}$ is bounded by $\|\hstar{\cdot}\|_{\infty} + \|\laginf\|_{\infty}$. By Lemma \ref{LEMM::BIAS} the rotated cost is non-negative, hence, $\CostT{\cdot}{t}$ is non-decreasing in $t$.

    Fix any initial state $\bx$. Now note that $\CostT{\bx}{t}$ is a monotone increasing sequence in $t$ for the compact domain of $\bx$ that converges point-wise to $\CostT{\bx}{\infty}$. Further, under Assumption \ref{AS::SA}, $\CostT{\bx}{t}$ is a continuous (Lipschitz) function  for any fixed $t$ by Lemma \ref{LEM::INFCOST}. The uniform convergence result now follows by Dini's theorem. Hence, there exists a $\THor(\epsilon)$ such that for any $\THor \geq \THor(\epsilon)$, $|\CostT{\bx}{\THor + 1} - |\CostT{\bx}{\THor}| < \epsilon $. Completing the proof. 
\end{proof}

\section{Proof of Jensen's Gap lemma}\label{APP:JENGAP}

The proof technique in this section is inspired by the result of the same name in \cite{Brown&Zhang25}. Much like in Subsection \ref{SUBSEC::DISSIPATIVITY}, for the purposes of this section we will explicitly distinguish between the states and actions for our fluid policies. Recall, given any joint distribution over the state $\bx$ and the corresponding state-action vector $\bYN$, one can recover the action by setting $\bu = \bYUN$ and the state to $\bx = \bYUN + \bYZN$. Hence, we have an equivalence between $(\bx, \bu)$ and $\bYN$. We can also describe the constraint set for our actions to be the set of actions such that $\bu \leq \bx$ in a component-wise sense and a bound on the $l^1$ norm of the action $\sum_{n = 1}^{N}\|\bu^{(n)}\|_1 \leq \alpha \sum_{n = 1}^{N}\|\bx^{(n)}\|_1 = \alpha N$. We will use $\mcU{\bx}$ to denote this set. Where, as has been the case in the rest of the paper, $\bu^{(n)}$ will indicate the action(sometimes referred to as \emph{control} in the literature) associated with the $n^{\text{th}}$ arm. Note, this is a $\sspaceN{n}$ dimensional vector rather than a binary decision since we are solving a fluid problem. We now denote the transition kernel for the next state by $\Phi(\bx, \bu)$. Concretely, given a state $\bx$ and action $\bu$, recall that $\Phi$ takes the form:
\[
\Phi(\bx, \bu) := \bx \cdot\PmatNN{0}  + \bu \cdot (\PmatNN{1} - \PmatNN{0})
\]
Recall the dynamic program \eqref{EQ::DYN1}, we rewrite it explicitly in terms of the state and control below:
\begin{align}\label{EQ::APDYN}
    \VLPROB{t}{\bx} =
    \begin{cases}
        &\max_{\bu \in \mcU{\bx}} \langle \Rvec{0}, \bx \rangle + \langle \Rvec{1} - \Rvec{0}, \bu  \rangle + \VLPROB{t - 1}{\Phi(\bx, \bu)}\\
        &\langle \laginf, \bx \rangle \hspace{0.3 in} \text{when $t = 0$}
    \end{cases}
\end{align}
Let time $t$ and state $\bx$ be fixed, we will denote by $\ufl{\bx}$ the \emph{fluid solution} to the dynamic program.
We recall the following properties of the dynamic program and its solution:
\begin{enumerate}
    \item Lemma \ref{APP::LEMCONCAVE} states that $\VLPROB{t}{\bx}$ is concave in $\bx$.
    \item By Lemma \ref{LEM::VALT}, $\VLPROB{t}{\bx}$ is Lipschitz continuous in $\bx$ with Lipschitz constant bounded above by $\left( \|\laginf\|_{\infty} + \left(1 + \frac{k}{\syncconst{k}}\right)\frac{k}{\syncconst{k}}\right)$.
\end{enumerate}
 It will be helpful to denote by $\Qval{t}{\bx}{\bu}$ the finite time horizon $Q$-function :
\begin{equation}\label{EQ::QVAL}
    \Qval{t}{\bx}{\bu} := \langle \Rvec{0}, \bx \rangle + \langle  \Rvec{1} - \Rvec{0}, \bu \rangle + \VLPROB{t - 1}{\Phi(\bx, \bu)}
\end{equation}
It is important to note that the $Q$-function we have defined here does not necessarily restrict $\bu$ to be a feasible policy. Therefore, typical properties associated with a Q-function do not follow \emph{unless we restrict ourselves to policies that are feasible.}

We are now ready to show the converse result on Jensen's lemma, it shows that for any fixed $t$, the Jensen's gap as defined in Section \ref{SUBSEC:JENSENGAP}, is order $\left( t \sqrt{\frac{1}{N}}\right)$. 

Let us first make some notation precise before we proceed with the lemma. Let, \\$\bx:= \{\bx^{(1)}, \bx^{(2)}, \dots \bx^{(N)}\}$ denote any initial state in $\Pi_{n = 1}^{N} \Delta(\Real{\sspaceN{n}})$. 
Given a component $\bx^{(i)}$, the number $S_i$ is sampled from this distribution with probability $\bx^{(i)}_{S_i}$. We then use $X^{(i)}$ to denote the one-hot encoded random vector:
 \begin{equation}
    X^{(i)}_{\hat{s}} =
     \begin{cases}
            1 \hspace{0.2 in}& \hspace{0.2 in} \hat{s} = S_i \\
            0 \hspace{0.3 in}&\text{otherwise.}
        \end{cases}
 \end{equation}
We will do so for each of the $N$ components independently at random and call the random vector $X = \{X^{(1)}, X^{(2)} \dots X^{(N)}\}$, 
\emph{a random vector generated by sampling each component of $\bx$}, denoted by $X \sim \bx$.

The following intermediate result regarding the continuity of the $Q$ function, $\Qval{t}{\bx}{\bu}$ is necessary for the main result we are looking to prove. The proof of this result follows directly from Lemma \ref{LEM::VALT} and so we avoid stating it here. 
\begin{lemma}\label{LEM::QVALLIP}
   The Q function, $\Qval{t}{\bx}{\bu}$ is Lipschitz in the policy $\bu$. Concretely, let $\bu$ and $\Tilde{\bu}$ be any two actions, then we have that
   \begin{equation}
       \left|\Qval{t}{\bx}{\bu} - \Qval{t}{\bx}{\Tilde{\bu}} \right| \leq \frac{1}{N}\left( 1 + \|\mu\|_{\infty} + \left(1 + \frac{k}{\syncconst{k}} \right)\sum_{l = 1}^{t}\left(1 - \frac{\syncconst{k}}{k} \right)^{l - 1} \right) \|\bu - \Tilde{\bu}\|_1.
   \end{equation}
   Where we have fixed the state $\bx$ for both functions. 
\end{lemma}


\LemJen*

\begin{proof}
    Our proof will proceed through induction. We begin with the verification step for $t = 0$,
    \begin{align*}
        &\VLPROB{0}{\bx} - \E_{X \sim \bx} \left[ \VLPROB{0}{X}\right] \\
        =& \langle \laginf, \bx \rangle - \E_{X \sim \bx} \langle \laginf, X \rangle = 0
    \end{align*}
    The inductive hypothesis assumes that for any $\bx$ distributed over $\sspace$ and $t > 0$, we have:
    \begin{equation}
        \VLPROB{t}{\bx} - \E_{X \sim \bx} \left[ \VLPROB{t}{X}\right] \leq \left[\left(1 + \frac{k}{\syncconst{k}} \right)\frac{k}{\syncconst{k}}  + \|\laginf\|_{\infty} + 1\right]\left( \frac{t}{2} \sqrt{\frac{\pi}{N}} \right)
    \end{equation}
    At time $t + 1$, let the state be $\bx := \{\bx^{(1)}, \bx^{(2)} \dots \bx^{(N)}\}$, we now draw the random variable $X$ according to $\bx$, $X \sim \bx$. Next, let the fluid policy give action $\ufl{\bx} := \{u^{(1)}, u^{(2)} \dots u^{(N)}\}$ for the $t$ horizon dynamic program. We may now generate a policy $\bu(X)$, component-wise as follows:
    \begin{equation}\label{EQ::DEFPOL1}
        \bu^{(i)}(X) :=
        \begin{cases}
            \frac{u^{(i)}_{s_i}}{\bx^{(i)}_{s_i}} \hspace{0.2 in} &\text{if  $i^{\text{th}}$ arm is in state $s_i$} \\
            0 \hspace{0.2 in}&\text{otherwise}
        \end{cases}
    \end{equation}
    with $\bu_{s}^{(i)} = 0$ if $x_{s} = 0$ for consistency. Now note that this new policy may not satisfy the action constraint i.e, $\frac{1}{N}\|\bu(X)\|_1$ need not be less than $\alpha$. In order to ensure that the action prescribed by the fluid solution corresponds to a valid policy we construct a policy $\hat{\bu}(X)$ so that,
    \begin{equation}\label{EQ::DEFPOL2}
        \hat{\bu}^{(i)}(X) :=
        \begin{cases}
            \bu^{(i)}(X) \hspace{0.2 in}&\text{if $\sum_{l = 1}^{i} \bu^{(l)}(X) < N \alpha$} \\
            0 \hspace{0.3 in}&\text{otherwise}
        \end{cases}
    \end{equation}
    Clearly, $\hat{\bu}(X)$ is always a valid policy. Hence, the fluid solution will always dominate the Q-function when choosing action $\hat{\bu}(X)$,
    \[
    \VLPROB{t + 1}{X} \geq  \Qval{t + 1}{X}{\hat{\bu}(X)}.
    \]
    If $\frac{1}{N} \E\|\bu(X) - \hat{\bu}(X)\|_1 \leq \sqrt\frac{\pi}{N}$ (proof differed till Lemma \ref{LEM::extract}), then from Lemma \ref{LEM::QVALLIP} it follows that,
    \begin{align*}
        \E &\left[\VLPROB{t + 1}{X}\right] \geq \E \left[\Qval{t + 1}{X}{\hat{\bu}(X)}\right]\\
        \geq& \E \left[\Qval{t + 1}{X}{\bu(X)}\right] - \left( 1 + \|\mu\|_{\infty} + \left(1 + \frac{k}{\syncconst{k}} \right)\frac{k}{\syncconst{k}} \right) \frac{1}{2}\sqrt\frac{\pi}{N}\\
        = & \E \left[\langle \Rvec{0}, X \rangle + \langle  \Rvec{1} - \Rvec{0}, \bu(X) \rangle + \VLPROB{t}{\Phi(X, \bu(X))}\right]\\
        &\hspace{0.4 in}- \left( 1 + \|\mu\|_{\infty} + \left(1 + \frac{k}{\syncconst{k}} \right)\frac{k}{\syncconst{k}} \right)\frac{1}{2} \sqrt\frac{\pi}{N} .
    \end{align*}
 
    Finally, using Jensen's inequality we have: 
    \begin{align*}
        &\E_{X \sim \bx} \left[\langle \Rvec{0}, X \rangle + \langle  \Rvec{1} - \Rvec{0}, \bu(X) \rangle + \VLPROB{t}{\Phi(X, \bu(X))}\right]\\
        \geq & \E_{X \sim \bx} \E_{\Tilde{X} \sim \Phi(X, \bu(X))} \left[\langle \Rvec{0}, X \rangle + \langle  \Rvec{1} - \Rvec{0}, \bu(X) \rangle + \VLPROB{t}{\Tilde{X}}\right].
    \end{align*}
    Combining these inequalities with the inductive hypothesis leads us to the following:
    \begin{align*}
        &\VLPROB{t + 1}{\bx} - \E_{X \sim \bx} \left[ \VLPROB{t + 1}{X}\right]\\
        \leq & \left( 1 + \|\mu\|_{\infty} + \left(1 + \frac{k}{\syncconst{k}} \right)\frac{k}{\syncconst{k}} \right)\frac{1}{2} \sqrt\frac{\pi}{N} + \VLPROB{t}{\Phi(\bx, \bu)} - \E_{X \sim \bx} \E_{\Tilde{X} \sim \Phi(X, \bu(X))} \left[\VLPROB{t}{\Tilde{X}}\right]\\
        = & \left( 1 + \|\mu\|_{\infty} + \left(1 + \frac{k}{\syncconst{k}} \right)\frac{k}{\syncconst{k}} \right)\frac{1}{2} \sqrt\frac{\pi}{N} + \VLPROB{t}{\Phi(\bx, \bu)} - \E_{\Tilde{X} \sim \Phi(\bx, \bu(\bx))} \left[\VLPROB{t}{\Tilde{X}}\right]\\
        \leq & \frac{(t + 1)}{2}\left( 1 + \|\mu\|_{\infty} + \left(1 + \frac{k}{\syncconst{k}} \right)\frac{k}{\syncconst{k}} \right)\sqrt{\frac{\pi}{N}}.
    \end{align*}
    \qed
\end{proof}

We now return to the differed proof regarding the expected actions in the lemma below.
\begin{lemma}\label{LEM::extract}
    Let $\bu(X)$ and $\hat{\bu}(X)$ be drawn according to \eqref{EQ::DEFPOL1} and \eqref{EQ::DEFPOL2} respectively, then the expected difference between the two actions is bounded as follows:
    \begin{equation}
        \frac{1}{N} \E\|\bu(X) - \hat{\bu}(X)\|_1 \leq \frac{1}{2}\sqrt\frac{\pi}{N}
    \end{equation}
\end{lemma}
\begin{proof}
    Note, for any two $X := \{X^{(1)}, X^{(2)} \dots X^{(i)}, \dots X^{(N)}\}$, $\Tilde{X} := \{X^{(1)}, X^{(2)} \dots \Tilde{X}^{(i)}, \dots X^{(N)}\}$ which deviate in only one component, we have:
    \begin{align*}
        &\frac{1}{N} \left(\|\bu(X)\|_1 -  \|\bu(X')\|_1\right)\\
        =& \frac{1}{N}\left(\frac{u^{(i)}(S_i)}{\bx^{(i)}(S_i)} - \frac{u^{(i)}(\Tilde{S}_i)}{\bx(\Tilde{S}_i)} \right) \leq \frac{1}{N}.
    \end{align*}
    The last inequality follows from noting that $u^{(i)}(S_i) \leq x(S_i)$ for any fluid policy. Hence, the ratio is strictly non-negative and less than $1$.
    
    We split the problem into two cases.
    
    \textbf{Case 1: } When $X$ is such that $\|\bu(X)\|_1 \leq \alpha N$, $\bu(X) = \hat{\bu}(X)$.  Hence, $\|\bu(X) - \hat{\bu}(X)\|_1 = 0$.
    
    \textbf{Case 2: } When $\|\bu(X)\|_1 \geq \alpha N \geq \E_X \|\bu(X)\|_1$, by Mcdiarmid's inequality \citep{McDiarmid_1989}, we may conclude that :
    \begin{equation}
        \Pr \left[ \frac{1}{N} \left(\|\bu(X)\|_1 - \E_X \|\bu(X)\|_1 \right) > \varepsilon \right] \leq \exp{(-N \varepsilon^2)}.
    \end{equation}
    Note, under Case 2, $\|\bu(X)\|_1 \geq \alpha N = \|\hat{\bu}(X)\|_1 \geq \E \|\bu(X)\|_1.$
    \begin{align*}
        \Pr \left[ \frac{1}{N} \left(\|\bu(X)\|_1 - \|\hat{\bu}(X)\|_1 \right) \geq \varepsilon\right] \leq \Pr \left[ \frac{1}{N} \left(\|\bu(X)\|_1 - \E_{X}\|\bu(X)\|_1 \right) \geq \varepsilon\right].
    \end{align*}
    We now perform a standard trick to conclude the proof:
    \begin{align*}
        \E \left[ \frac{1}{N} \left(\|\bu(X)\|_1 - \|\hat{\bu}(X)\|_1 \right)\right] =& \int_{0}^{N} \Pr \left[ \frac{1}{N} \left(\|\bu(X)\|_1 - \|\hat{\bu}(X)\|_1 \right) \geq \varepsilon\right] d\varepsilon\\
        \leq& \int_{0}^{\infty} \exp{(-N \varepsilon^2)}d \varepsilon
        = \frac{1}{2}\sqrt\frac{\pi}{N}.
    \end{align*}
    Combining the two cases gives us the result. 
\end{proof}

\section{Final Steps of Proof of Theorem \ref{THM::ALMOSTOPT}}\label{APP:SEC:FP}

\subsection{Deriving Inequality (a)}

We note from the tower property and a little algebraic manipulation, we obtain the following:
\begin{align*}
    &\frac{1}{T} \E\sum_{t = 0}^{T - 1} \langle \laginf, \bxs{t} \rangle - \frac{1}{T} \E\sum_{t = 0}^{T - 1} \langle \laginf, \bYFLT{t} \cdot \PmatN \rangle \\
    =& \frac{1}{T}\langle \laginf, \bxs{0} \rangle + \frac{1}{T} \E\sum_{t = 0}^{T - 2} \langle \laginf, \bxs{t + 1} \rangle - \frac{1}{T} \E\sum_{t = 0}^{T - 2} \langle \laginf, \bYFLT{t} \cdot \PmatN \rangle - \frac{1}{T} \E \langle \laginf, \PmatN \cdot \bYFLT{T - 1} \rangle\\
    \leq &\frac{1}{T} \left| \E\sum_{t = 0}^{T - 2} \E \left[ \langle \laginf, \bxs{t + 1} \rangle -  \langle \laginf, \bYFLT{t} \cdot \PmatN \rangle| \bs(t), \bANN(t)\right] \right| + O\left( \frac{1}{T}\right)
\end{align*}
Now we consider each term in the sum.
\begin{align*}
    &\E \left[ \langle \laginf, \bxs{t + 1} \rangle -  \langle \laginf, \bYFLT{t} \cdot \PmatN \rangle| \bs(t), \bANN(t)\right] = (*)
\end{align*}

Let $[\bANN(t)]$ denote the $\Real{|\sspace\|}$ dimensional random vector produced using randomized rounding. Recall, our algorithm picks $[\bANN(t)]$ as follows,
\begin{equation*}
\bxs{t} \xrightarrow[\eqref{EQ::FLUIDLP1}-\eqref{EQ::FLUIDLP4}]{\text{$\THor(\bxs{t})$-horizon LP}} \bYFLT{t} \xrightarrow[\text{rounding}]{\text{Randomized }} \bANN(t) \xrightarrow[]{\eqref{EQ:MKEVOL}} \bs(t + 1).
\end{equation*}

This process implies that the fluid actions $\bYFLT{t}$ are a deterministic function of $\bxs{t}$. Hence,
\begin{align*}
   (*) =& \langle \laginf, \E \left[\bxs{t + 1}| \bs(t), \bANN(t)\right] - \bYFLT{t} \cdot \PmatN \rangle\\
   \leq& \|\laginf\|_{\infty} \frac{1}{N} \E\|[\bANN(t)] - [\bYFLU{t}]\|_1 \leq \|\laginf\|_{\infty} \left(\frac{\alpha N - \lfloor \alpha N \rfloor}{N} \right)
\end{align*}
Plugging this value back into the sum gives us the inequality (a).

\subsection{Deriving Inequality (b)}
We set the value of $\THor$ so that, $\CostT{\bYFLT{t} \cdot \PmatN}{\THor + 1} - \CostT{\bYFLT{t} \cdot \PmatN}{\THor} < \epsilon$. The existence of such a $\THor$ follows from inequality \eqref{EQN::INFCOST3}, Lemma \ref{LEM::INFCOST} and the dynamic programming principle, 
\begin{align}\label{EQN::INEQb1}
    \rotcost{c}{\bxs{t}, \ufl{t}} = \CostT{\bxs{t}}{\THor} - \CostT{\bYFLT{t} \cdot \PmatN}{\hat\THor - 1} \nonumber\\
    \leq \epsilon + \CostT{\bxs{t}}{\THor} - \CostT{\bYFLT{t} \cdot \PmatN}{\THor}
\end{align}
 Denoting this time by $\THor(t + 1)$ we get: 
\begin{align}\label{EQN::INEQb2}
    \frac{1}{T}\sum_{t = 0}^{T - 1}\rotcost{c}{\bxs{t}, \ufl{t}} \leq& \epsilon + \frac{1}{T}\sum_{t = 0}^{T - 2}\CostT{\bxs{t + 1}}{\THor} - \frac{1}{T}\sum_{t = 0}^{T - 2}\CostT{\bYFLT{t} \cdot \PmatN}{\THor} + O\left(\frac{1}{T}\right)
\end{align} 

\subsection{Deriving Inequality (d)}

We will use our notation from \ref{APP:JENGAP} to show our result. We will fix $\THor$ for the purposes of what is to follow. Consider the following terms,
\begin{align*}
    &\E \left[\VLPROB{\THor}{\bYFLT{t} \cdot \PmatN} - \VLPROB{\THor}{\bxs{t + 1}} | \bxs{t}, \bANN(t) \right]\\
     =& \VLPROB{\THor}{\bYFLT{t} \cdot \PmatN} - \E_{\Tilde{X}(t + 1) \sim \Phi(\bxs{t}, \bANN(t))} \left[\VLPROB{\THor}{\Tilde{X}(t + 1)} | \bxs{t}, \bANN(t)\right]\\
     =& \VLPROB{\THor}{\bYFLT{t} \cdot \PmatN} - \E_{X(t + 1) \sim \Phi(\bxs{t}, \ufl{t}} \left[\VLPROB{\THor}{X(t + 1)} | \bxs{t}, \bANN(t)\right]\\
     &+ \E_{X(t + 1) \sim \Phi(\bxs{t}, \ufl{t}} \left[\VLPROB{\THor}{X(t + 1)} | \bxs{t}, \bANN(t)\right]\\
     &- \E_{\Tilde{X}(t + 1) \sim \Phi(\bxs{t}, \bANN(t)} \left[\VLPROB{\THor}{\Tilde{X}(t + 1)} | \bxs{t}, \bANN(t)\right]
\end{align*}
The first equality follows since $\bYFLT{t}$ is $\bxs{t}$ measurable and the second equality follows from definition. Now note, one can couple the states of $X(t + 1) := \{X^{(1)}(t + 1), X^{(2)}( t+ 1) \dots X^{(N)}( t+ 1)\}$ and $\Tilde{X}(t + 1) := \{\Tilde{X}^{(1)}(t + 1), \Tilde{X}^{(2)}(t + 1) \dots\Tilde{X}^{(N)}(t + 1)\}$ as follows: the $i^{\text{th}}$ state of the system evolve to same next state if they take the same actions $[\ufl{t}]^i, \bANN^{(i)}(t)$. Due to randomized rounding, the marginals on the actions $\bANN^{(i)}(t)$ are the same as $[\ufl{t}]^{(i)}$. Finally, by Lemma (\ref{LEM::VALT}), 
\begin{align*}
    &\E_{X(t + 1) \sim \Phi(\bxs{t}, \ufl{t})} \left[\VLPROB{\THor}{X(t + 1)} | \bxs{t}\right] \\
     &- \E_{\Tilde{X}(t + 1) \sim \Phi(\bxs{t}, \bANN(t))} \left[\VLPROB{\THor}{\Tilde{X}(t + 1)} | \bxs{t}, \bANN(t)\right] \\
     &\leq \left[\|\mu\|_{\infty} + \left(1 + \frac{k}{\syncconst{k}} \right)\sum_{l = 1}^{\THor}\left(1 - \frac{\syncconst{k}}{k} \right)^{l - 1}\right]\frac{\alpha N - \lfloor \alpha N \rfloor}{N}
\end{align*}

We may now substitute this value of the horizon into the Jensen gap lemma \ref{LEM::JENGAP} to obtain the following bound,
\begin{align}
    \THor(\epsilon)\Bigg(\Bigg(1 + \frac{k}{\syncconst{k}} \Bigg)\frac{k}{\syncconst{k}} & + \|\laginf\|_{\infty} + 1\Bigg)\sqrt{\frac{\pi}{N}} \nonumber\\
    +& \left[\|\mu\|_{\infty} + \left(1 + \frac{k}{\syncconst{k}} \right)\sum_{l = 1}^{\THor(\epsilon)}\left(1 - \frac{\syncconst{k}}{k} \right)^{l - 1}\right] \frac{\alpha N - \lfloor \alpha N \rfloor}{N}.
\end{align}
When combined with the bound on randomized rounding completes the proof.

     

\subsection{Sketch of Proof of Corollary \ref{COR::GENMDP}}\label{APP:GENMDP}
\begin{proof}
    Our proof sketch will simply draw parallels between the different components used in the Restless bandit setting and the more general Heterogeneous Weakly Coupled Setting. 
    \begin{enumerate}
        \item Firstly, it should be noted that one can define a gain and bias function $\gstarWMDP$ similar to the RMAB problem as a solution to the fixed point problem:
    \begin{align}\label{EQ:FXWMDP}
    \gstarWMDP :=& \max_{\bYN}\left[ \frac{1}{N} \sum_{n = 1}^N \Rvec{n} \cdot \bYinf{(n)} \right] \nonumber\\
    &\text{such that,} \nonumber \\
    &\sum_{a \in \mcl{A}} \bYWMDPN{(n)} = \sum_{{a \in \mcl{A}}} \bYWMDPN{(n)} \cdot \Pmat{n}{a}   \hspace{0.2 in } \forall n \nonumber\\
    &\sum_{n = 1}^{N} \sum_{(s,a)}\bYWMDPN{n}\constrWMDP{(n),(e)}{s}{a} \leq N \budgetk{e} \hspace{0.2 in} \text{for all $e$}
\end{align}
    \item Next, one can use Assumption \ref{AS:SA2} and \ref{AS:SA3} to show that a parallel notion of a bias function, $\hstarWMDP{\cdot}$, defined by the point-wise limit:
    \[
    \hstarWMDP{\bx} := \lim_{\THor \to \infty}\THor \gstarWMDP - \VLWMDP{\THor}{\bx}
    \]
    is Lipschitz continuous and hence, bounded in the domain of $\bx$. 
    \item Assumption \ref{AS:SA2} and \ref{AS:SA3} ensure that the constraint set is non-empty. Therefore, the linear program, \eqref{EQ:FXWMDP} has a nonempty set that satisfies the fixed-point constraint and hence, satisfies strong duality. Let $\bYstarN$ be the fixed point solution to \eqref{EQ:FXWMDP}. If $\laginf$ is the Lagrange multipliers corresponding to the fixed point constraint, one can construct a rotated cost function,
    \begin{align*}
        0 = \gstar -  \frac{1}{N} \sum_{n = 1}^{N}\Rvec{n} \cdot \bYstar{n} \leq \gstar - \frac{1}{N} \sum_{n = 1}^{N}\Rvec{n} \cdot \bYinf{n} + \langle \laginf, (\bYZN + \bYUN) - \bYN \cdot \PmatN \rangle
    \end{align*}
    and recover the dissipativity result in the more general setting. The monotonicity of the equivalent surrogate loss function follows.
    
    \item Since, we cannot directly use randomized rounding to solve our problem, we solve the fluid problem and we construct a feasible policy for the fluid problem exactly as was done in \ref{APP:JENGAP} for the proof of Lemma \ref{LEM::JENGAP}.
    
    \item Since, we assume that there exists a constant $b$ independent of $N$ such that our constraints for each arm are bounded, $0 \leq \constrWMDP{(n), (e)}{\cdot}{\cdot} \leq b$, Mcdiarmid's inequality continues to hold ensuring much of the mechanisms of our proof of Jensen's Gap remains unchanged in this case. The previous bounds on the difference between the fluid solution and the expected value of the fluid solution will continue to hold in this case.
    \end{enumerate}
    Given these parallels, one can assemble the last steps of the proof exactly as we have done in Section \ref{SUBSEC::STEPSOFPROOF} to complete our result. 
\end{proof}
\section{Proof of Theorem \ref{THM::ALMOSTOPT2}}\label{APP::THM2}
We begin with a parallel definition of the surrogate cost function, Definition \ref{DEF::SURCST} when our time horizon $\hat{\THor}(\bx, t)$ is set to $B/\epsilon - t \mod B/\epsilon $ as required by our finite horizon variant. Given any initial joint state $\bs(0) = \bs$, let $\bs(t)$ be the evolution of the states over a time interval $t \in \{0, 1 \dots \THor - 1\}$ under the modified LP-update policy. We denote by $\bxs{(t)}$ its one-hot coded state vector. Next, we let $\bualg{t}$ denote the policy when computing a $\hat{\THor}(\bx, t) = B/\epsilon - t \mod B\epsilon$ horizon dynamic program. $\bAalg{t}$ then denotes the actions taken by Algorithm \ref{algo::MPC} after randomized rounding. The new surrogate loss function for any horizon $\THor$ is set to,
\begin{equation}\label{EQN::COSTFIN}
    \hat{C}_{\THor}(\bxs{\bs}) := \sum_{t = 0}^{\THor - 1} \rotcost{c}{\bxs{t}, \bualg{t}}
\end{equation}

\begin{lemma}\label{LEM::APBND}
    The cost \eqref{EQN::COSTFIN} for any $\THor < B/\epsilon$ can be bounded as follows:
    \[
    \hat{C}_{\THor}(\bxs{\bs}) \leq B + \frac{\THor (\THor - 1)}{4}\left[\left(1 + \frac{k}{\syncconst{k}} \right)\frac{k}{\syncconst{k}}  + \|\laginf\|_{\infty} + 1\right]\left( \sqrt{\frac{\pi}{N}} \right)
    \]
\end{lemma}
\begin{proof}
    From the dynamic programming principle we have,
    \[
    \rotcost{c}{\bxs{t}, \bualg{t}} = \CostT{\bxs{t}}{\hat{\THor}(\bx, t)} - \CostT{\Phi(\bxs{t}, \bualg{t})}{\hat{\THor}(\bx, t) - 1}.
    \]
    It then follows that,
    \begin{align*}
       &\hat{C}_{\THor}(\bxs{\bs}) =  \sum_{t = 0}^{\THor - 1} \rotcost{c}{\bxs{t}, \bualg{t}}\\
       =&\CostT{\bxs{t}}{\THor} + \sum_{t = 1}^{\THor - 1}\CostT{\bxs{t + 1}}{\hat{\THor}(\bx, t) - 1} - \CostT{\Phi(\bxs{t}, \bualg{t})}{\hat{\THor}(\bx, t) - 1}.
    \end{align*}
    Now one invokes the Jensen's gap lemma, Lemma \ref{LEM::JENGAP} on each term in the summation and we note that $\CostT{\bxs{t}}{t} \leq B$ which yields the required result.
\end{proof}
Consider a time horizon $T$, fix $\epsilon$ and recall that $B$ is an upper bound on $\CostT{\bx}{\infty}$ for any $\bx$. We split this time horizon into segments of length $\frac{T}{(B/\epsilon)}$. We are now ready to complete our proof.

\begin{proof}[Proof of Theorem \ref{THM::ALMOSTOPT2}\\]

Using the same steps as those that brought us inequality $(a)$ for an interval of length $T$ from Section \ref{SUBSEC::STEPSOFPROOF}
    \begin{align*}
        \gstar - \VLPF \leq& \frac{1}{T} \left[ \sum_{t = 0}^{T - 1} \rotcost{c}{\bxs{t}, \bualg{t} (\bxs{t}}\right] + \left(1 + \|\laginf\|_{\infty} \right)\frac{\alpha N - \lfloor \alpha N \rfloor}{N} + \mcl{O}\left(\frac{1}{T} \right)\\
        =&\frac{1}{T} \times \frac{T}{B/\epsilon}\left[ \sum_{t = 0}^{B/\epsilon - 1} \rotcost{c}{\bxs{t}, \bu_{t}(\bxs{t}}\right] + \left(1 + \|\laginf\|_{\infty} \right)\frac{\alpha N - \lfloor \alpha N \rfloor}{N} + \mcl{O}\left(\frac{1}{T} \right)\\
        =& \frac{1}{T} \times \frac{T}{B/\epsilon}\left[ \hat{C}_{B/\epsilon}(\bxs{\bs})\right] + \left(1 + \|\laginf\|_{\infty} \right)\frac{\alpha N - \lfloor \alpha N \rfloor}{N} + \mcl{O}\left(\frac{1}{T} \right)\\
        \leq& \frac{1}{T} \times \frac{T}{B/\epsilon}\left(B +  \frac{B/\epsilon (B/\epsilon - 1)}{4}\left[\left(1 + \frac{k}{\syncconst{k}} \right)\frac{k}{\syncconst{k}}  + \|\laginf\|_{\infty} + 1\right]\sqrt{\frac{\pi}{N}}\right)\\
        &+ \left(1 + \|\laginf\|_{\infty} \right)\frac{\alpha N - \lfloor \alpha N \rfloor}{N} + \mcl{O}\left(\frac{1}{T} \right)\\
        =&  \epsilon + \frac{B}{4\epsilon}\left[\left(1 + \frac{k}{\syncconst{k}} \right)\frac{k}{\syncconst{k}}  + \|\laginf\|_{\infty} + 1\right]\sqrt{\frac{\pi}{N}} + \left(1 + \|\laginf\|_{\infty} \right)\frac{\alpha N - \lfloor \alpha N \rfloor}{N} + \mcl{O}\left(\frac{1}{T} \right)
    \end{align*}
    Where the inequality follows from Lemma \ref{LEM::APBND}. This completes the proof. 
\end{proof}

\section{Proof of Proposition \ref{PROP:WHITLUPDATE}\\}\label{APP::PROOFWHITLUPDATE}
For the purposes of this problem we will replace the action constraint $\sum_{n = 1}^{N} a_n \leq \budget N$ to an equality $\sum_{n = 1}^{N} a_n = \budget N$. One way to do so is by following the arguments proposed in \citet{GGY23b}: create a collection of $\alpha N$ dummy arms that have a single state and provide no rewards at all. Whenever the optimal action demands that less than $\budget N$ arms are pulled, the remaining arms can be selected from this pool of dummy arms. Note, such a system does not differ in any way from our own system in terms of state evolution or rewards. We do this to make the next set of arguments easier to follow.

We consider the following primal problem (from \cite{puterman2014markov}),
\begin{align}\label{EQ:PUTERMAN}
    \min_{g, \laginf_{rel}} g &+ \blamT_{rel} \budget \\
    &\text{s.t.} \nonumber\\
    g + \laginf_{rel,s_i}^{(i)} -& \sum_{s' \in \sspaceN{i}}P^{(i)}_{s'|s_i, a}\laginf_{rel,s'}^{(i)} \geq r^{(i)}(s_i, a) - a \blamT_{rel} \text{ for all } (s_i, a).
\end{align}
The dual LP to the primal \eqref{EQ:PUTERMAN} is given by, 
\begin{align*}
    \max_{\bYN \in \mcl{Y}}& \frac{1}{N}\sum_{n = 1}^{N} \Rvec{n} \cdot \bYinf{(n)} \\
    &\text{such that,} \nonumber \\
    &\bYZinf{(n)} + \bYUinf{(n)} =  \bYZinf{(n)} \cdot \Pmat{n}{0} + \bYUinf{(n)} \cdot \Pmat{n}{1}  \hspace{0.2 in } \forall n \\
    &\sum_{n = 1}^{N + \alpha N} \| \bYUinf{n}\|_1 = N. \alpha \hspace{0.2 in }
\end{align*}
where we choose $\blamT_{rel}$ so that the two problems are equivalent. Further, the solution to the fixed point problem can be found by computing the two variables $(\gstar, \laginf_{rel})$ that satisfy the fixed point equation \eqref{EQ::LPPRIO}, 
\begin{equation*}
    \gstar + \laginf_{rel, s_i}^{(i)} = \max_{a \in \{0, 1\}}\{r(s_i, a) - a \blamT_{rel} + \sum_{s' \in \sspaceN{i}}P^{(i)}_{s'|s_i, a}\laginf_{rel, s'}^{(i)}\}.
\end{equation*}
and $\laginf_{rel}$ is a Lagrange multiplier corresponding to the fixed point constraint. Consider now the LP-update policy when $\THor$ is set to $1$.  
\begin{align*}
    \VLPROB{1}{\bx} =&
        \max_{\{\bYN |\bYN_{\cdot, 1} + \bYN_{\cdot, 0} = \bx\}}\frac{1}{N} \sum_{n = 1}^N  \Rvec{n} \cdot \bYinf{(n)}  + \langle \laginf,\bYN_{\cdot, 1} \cdot \PmatNN{1} + \bYN_{\cdot, 0} \cdot \PmatNN{0}\rangle\\
        \text{s.t.}&\\
        \sum_{n = 1}^{N + \alpha N} \| \bYUinf{n}\|_1 =& \budget N.
\end{align*}
There is at most one randomization for a $1$ horizon Linear program, see for example ~\cite{altman1999constrained}, hence, this LP is solved by constructing a priority policy. Let $\blamstarT(\bx)$ be the lagrange multiplier, then an arm $i$ in state $s_i$ will have higher priority then arm $j$ in state $s_j$ when,
\begin{align*}
    &r^{(i)}(s_i, 1) -  \blamstarT({\bx}) + \sum_{s' \in \sspaceN{i}}P^{(i)}_{s'|s_i, 1}\laginf_{ s'}^{(i)} - r^{(i)}(s_i, 0) - \sum_{s' \in \sspaceN{i}}P^{(i)}_{s'|s_i, 0}\laginf_{ s'}^{(i)} \\
    \geq& r^{(j)}(s_j, 1) -  \blamstarT({\bx}) + \sum_{s' \in \sspaceN{i}}P^{(i)}_{s'|s_j, 1}\laginf_{ s'}^{(j)} - r^{(j)}(s_j, 0) - \sum_{s' \in \sspaceN{i}}P^{(i)}_{s'|s_j, 0}\laginf_{ s'}^{(j)}\\
    \iff& r^{(i)}(s_i, 1)  + \sum_{s' \in \sspaceN{i}}P^{(i)}_{s'|s_i, 1}\laginf_{ s'}^{(i)} - r^{(i)}(s_i, 0) - \sum_{s' \in \sspaceN{i}}P^{(i)}_{s'|s_i, 0}\laginf_{ s'}^{(i)} \\
    \geq& r^{(j)}(s_j, 1) + \sum_{s' \in \sspaceN{i}}P^{(i)}_{s'|s_j, 1}\laginf_{ s'}^{(j)} - r^{(j)}(s_j, 0) - \sum_{s' \in \sspaceN{i}}P^{(i)}_{s'|s_j, 0}\laginf_{ s'}^{(j)}
    \end{align*}
     From Definition \ref{DEF:LAGMUL}, $\laginf$ is another collection of multipliers for the same constraint. It now follows that $\laginf_{rel} - \laginf = \Bar{c}\Onne$ for some constant $\Bar{c}$ i.e, the two lagrange multipliers differ by a constant. 
    \begin{align*}
    \iff& r^{(i)}(s_i, 1)  + \sum_{s' \in \sspaceN{i}}P^{(i)}_{s'|s_i, 1}(\laginf_{ s'}^{(i)} + \Bar{c}) - r^{(i)}(s_i, 0) - \sum_{s' \in \sspaceN{i}}P^{(i)}_{s'|s_i, 0}(\laginf_{ s'}^{(i)} + \Bar{c}) \\
    \geq& r^{(j)}(s_j, 1) + \sum_{s' \in \sspaceN{i}}P^{(i)}_{s'|s_j, 1}(\laginf_{ s'}^{(j)} + \Bar{c}) - r^{(j)}(s_j, 0) - \sum_{s' \in \sspaceN{i}}P^{(i)}_{s'|s_j, 0}(\laginf_{ s'}^{(j)} + \Bar{c})\\
    \iff& r^{(i)}(s_i, 1)  + \sum_{s' \in \sspaceN{i}}P^{(i)}_{s'|s_i, 1}\laginf_{rel, s'}^{(i)} - r^{(i)}(s_i, 0) - \sum_{s' \in \sspaceN{i}}P^{(i)}_{s'|s_i, 0}\laginf_{rel,s'}^{(i)}  
    \end{align*}
    \begin{align*}
    \geq& r^{(j)}(s_j, 1) + \sum_{s' \in \sspaceN{i}}P^{(i)}_{s'|s_j, 1}\laginf_{rel, s'}^{(j)} - r^{(j)}(s_j, 0) - \sum_{s' \in \sspaceN{i}}P^{(i)}_{s'|s_j, 0}\laginf_{rel,s'}^{(j)}\\
    \iff& r^{(i)}(s_i, 1) - \blamT_{rel} + \sum_{s' \in \sspaceN{i}}P^{(i)}_{s'|s_i, 1}\laginf_{rel, s'}^{(i)} - r^{(i)}(s_i, 0) - \sum_{s' \in \sspaceN{i}}P^{(i)}_{s'|s_i, 0}\laginf_{rel,s'}^{(i)} \\
    \geq& r^{(j)}(s_j, 1) - \blamT_{rel} + \sum_{s' \in \sspaceN{i}}P^{(i)}_{s'|s_j, 1}\laginf_{rel, s'}^{(j)} - r^{(j)}(s_j, 0) - \sum_{s' \in \sspaceN{i}}P^{(i)}_{s'|s_j, 0}\laginf_{rel,s'}^{(j)}\\
    \iff& \indxlp{i} \geq \indxlp{j}
\end{align*}
 Therefore, the LP-priority policy solves the LP-update policy when $\THor$ is set to $1$, completing the proof. 

\section{Randomized Rounding Procedure\\}\label{APP:RR}

Before describing the randomized rounding procedure, let us consider the following problem: we are given a vector of $n$ probabilities $z_1\dots z_N\in[0,1]$ and a budget $Z\in\{1\dots N\}$ such that $\sum_{i=1}^N z_i\le Z$. The randomized rounding problem is to find a coupling of $n$ Bernoulli random variables $Z_1\dots Z_N$ such that 
\begin{align}
    \label{eq:randomized_rounding}
    \Prob{Z_i}=1-\Prob{Z_i=0}=z_i\text{ and $\sum_{i=1}^N Z_i \le Z$ (almost surely).}
\end{align}
This problem can be solved using the algorithm described in Section~5.2.3 of \cite{ioannidis2016adaptive} as follows. Let $\bar{Z}=\lceil\sum_{i=1}^N z_i\rceil$ and define $z_{N+1}=\bar{Z} - \sum_{i=1}^N z_i$. Using the procedure described in the aforementioned paper, there exists a coupling of $k+1$ Bernoulli random variables $Z_1\dots Z_{N+1}$ such that $\sum_{i=1}^{N+1} Z_i = \bar{Z}$. By construction, $\bar{Z}\le Z$. Hence, the variables $Z_1\dots Z_N$ satisfy \eqref{eq:randomized_rounding}.

Adapting this procedure to our case requires minimal changes. In the randomized rounding problem presented in Section~\ref{SEC::IPL}, we are given a collection of variables $y_{1}\dots y_{N}\in[0,1]^N$ that satisfies $\sum_{n=1}^N y_{n}\le \alpha N$ and we want to find an action vector $A_{1}\dots A_N$ that satisfies $\sum_{n=1}^N A_n\le \alpha N$ (almost surely) and such that $\sum_{i=1}^N| |\E A_n-y_n|\le \alpha N - \floor{\alpha N}$. The only difference with the problem presented above is that $\alpha N$ might not be an integer. Our idea is therefore to the procedure of \cite{ioannidis2016adaptive} with the rescaled variables $z_n=y_n\floor{\alpha N}/(\alpha N)$. This is summarized in Algorithm~\ref{algo::rounding}.
\begin{algorithm}
	\begin{algorithmic}[1]
    \caption{Randomized rounding}
    \label{algo::rounding}
        \State \textbf{Input}: Vector $\bYN\in\sspace$ such that $\sum_n y_n\le\alpha N$ and $y_n\in[0,1]$. 
        \State Let $z_n = \frac{\floor{\alpha N}}{\alpha N} y_n$. 
        \State Use the procedure of \cite{ioannidis2016adaptive} to construct a collection of $N$ Bernoulli random variables $A_1\dots A_N$ such that $\E[A_n]=z_n$ and $\sum_{n=1}^N A_n\le \floor{\alpha N}$.
        \Ensure $\sum_{n=1}^N A_n\le\alpha N$ and $\sum_{n=1}^N |\E{A_n}-y_n|\le \alpha N - \floor{\alpha N}$.
    \end{algorithmic}
\end{algorithm}

To prove the correctness of Algorithm~\ref{algo::rounding}, we first remark that, by construction, $\sum_{n=1}^N A_n \le\floor{\alpha N}\le\alpha N$. For the second part, we have:
\begin{align*}
    \sum_{n=1}^N |\E{A_n}-y_n| &= \sum_{n=1}^N |z_n-y_n|\\
    &= \sum_{n=1}^N \left|\frac{\floor{\alpha N}}{\alpha N} y_n - y_n\right|\\
    &= \sum_{n=1}^N y_n \left(1 - \frac{\floor{\alpha N}}{\alpha N}\right)\\
    &\le \alpha N - \floor{\alpha N},
\end{align*}
where the last line is simply because $\sum_{n=1}^N y_n\le\alpha N$.

 \section{Experimental Details\\}\label{APP:EXP}

In this section we provide a details for our simulations results. We provide the exact matrices and budgets used in our counter examples for completeness. All our simulations are implemented in \texttt{Python}, by using \texttt{numpy} for the random generators and array manipulation and \texttt{pulp} to solve the linear programs. All simulations are run on a personal laptop (HP ENVY 2023). 

 \subsection{Example \cite{HXCW23}}
This example has $8$ dimensions and its parameters are:
\begin{align*}
    P_{ 0 }=\left(
    \begin{array}{cccccccc}
       1 & & & & & & & \\
       1 & & & & & & & \\
        &0.48 &0.52 & & & & & \\
        & &0.47 &0.53 & & & & \\
        & & & &0.9 &0.1 & & \\
        & & & & &0.9 &0.1 & \\
        & & & & & &0.9 &0.1 \\
       0.1 & & & & & & &0.9 \\
    \end{array}\right)
\end{align*}
\begin{align*}
    P_{ 1 }=\left(
    \begin{array}{cccccccc}
       0.9 &0.1 & & & & & & \\
        &0.9 &0.1 & & & & & \\
        & &0.9 &0.1 & & & & \\
        & & &0.9 &0.1 & & & \\
        & & &0.46 &0.54 & & & \\
        & & & &0.45 &0.55 & & \\
        & & & & &0.44 &0.56 & \\
        & & & & & &0.43 &0.57 \\
    \end{array}\right)
\end{align*}

$r_{ 1 } = 0 $ for all state and $r_{ 0 } =  (0, 0,0 ,0 ,0 ,0 ,0 ,0.1)$.

The budget is set to $\alpha=0.5$ (which is the parameter used in \cite{HXCW23}).

\subsection{Example \cite{chen:tel-04068056}}

This example has $|\sspace|=3$ dimensions and its parameters are:

\begin{align*}
    P_{ 0 }=\left(
    \begin{array}{ccc}
       0.022 &0.102 &0.875 \\
       0.034 &0.172 &0.794 \\
       0.523 &0.455 &0.022 \\
    \end{array}\right)\qquad 
    P_{ 1 }=\left(
    \begin{array}{ccc}
       0.149 &0.304 &0.547 \\
       0.568 &0.411 &0.020 \\
       0.253 &0.273 &0.474 \\
    \end{array}\right)
\end{align*}

$r_{ 1 } = ( 0.374, 0.117, 0.079) $ and $r_0=0$ for all the states. 

The budget is set to $\alpha=0.4$ (which was the parameter used in \cite{chen:tel-04068056}).


\end{document}